\documentclass[a4paper, 10pt]{amsart}

\usepackage[margin=0.75in,headheight=15pt,footskip=25pt]{geometry}

\usepackage{cite}

\usepackage{amsmath}
\usepackage{amssymb}
\usepackage{amsthm}
\usepackage{subcaption}

\usepackage{tikz}
\usepackage{pgfplots}
\usetikzlibrary{calc,patterns,angles,quotes}
\usepackage{enumerate}

\usepackage{hyperref}

\hypersetup{colorlinks=true,allcolors=blue}
\usepackage[usenames,dvipsnames]{pstricks}
\usepackage{siunitx,booktabs}
\def\doublewidetilde#1{\widetilde{\widetilde{#1\mkern0mu}}}

\theoremstyle{plain}
\newtheoremstyle{spacedplain}
  {8pt}      
  {8pt}      
  {\itshape} 
  {}         
  {\bfseries}
  {.}        
  {0.5em}    
  {}

\newtheoremstyle{spacedremark}
  {8pt}
  {8pt}
  {}         
  {}
  {\bfseries}
  {.}
  {0.5em}
  {}

\theoremstyle{spacedplain}
\newtheorem{theorem}{Theorem}[section]
\newtheorem{lemma}{Lemma}[section]

\newtheorem{assumption}{Assumption}[section]
\newtheorem{remark}{Remark}[section]
\newtheorem{definition}{Definition}[section]

\usepackage{algorithm}
\usepackage{algorithmicx} 
\usepackage[skip=6pt]{caption} 

\makeatletter
\def\@floatboxreset{\reset@font\raggedright} 
\makeatother

\numberwithin{figure}{section}
\numberwithin{equation}{section}

\usepackage{graphicx}%
\usepackage{multirow}%

\usepackage{mathrsfs}%
\usepackage[title]{appendix}%
\usepackage{xcolor}%
\usepackage{textcomp}%
\usepackage{manyfoot}%
\usepackage{booktabs}%

\usepackage{algpseudocode}%
\usepackage{listings}%

\usepackage{amssymb}
\usepackage{subcaption}

\usepackage{tikz}
\usepackage{pgfplots}
\usetikzlibrary{calc,patterns,angles,quotes}
\usepackage{enumerate}

\usepackage{hyperref}

\hypersetup{colorlinks=true,allcolors=blue}
\usepackage[usenames,dvipsnames]{pstricks}
\usepackage{siunitx,booktabs}
\usepackage{url}

\numberwithin{figure}{section}
\numberwithin{equation}{section}

\newcommand{\tribar}{\vert\kern-0.25ex\vert\kern-0.25ex\vert}

\newcommand{\M}{\bold{M}}
\renewcommand{\S}{\bold{S}}
\newcommand{\Q}{\bold{T}}
\newcommand{\D}{\bold{D}}

\newcommand{\CC}{\mathcal{C}}

\newcommand{\ddiv}{\text{div}\,}
\newcommand{\Gn}{\mathcal{G}^n}
\newcommand{\F}{\mathcal{F}}
\newcommand{\Gno}{\mathcal{G}^{n-1}}
\newcommand{\Th}{\mathcal{T}_h}
\newcommand{\Sh}{\mathcal{S}_{h}}

\newcommand{\N}{N}
\newcommand{\Nh}{\mathcal{N}_h}
\newcommand{\Nhb}{\mathcal{N}_h^{\partial}}
\newcommand{\Ehh}{\mathcal{E}_h}
\newcommand{\Ehb}{\mathcal{E}_h^{\partial}}
\newcommand{\NT}{N_0}
\newcommand{\Zh}{\mathcal{Z}_h}

\newcommand{\del}{\boldsymbol{\delta}}

\newcommand{\Rn}{\mathbb{R}^{\N}}
\newcommand{\R}{\mathbb{R}}
\newcommand{\al}{\boldsymbol{\alpha}}

\newcommand{\bpsi}{\boldsymbol{\psi}}
\newcommand{\be}{\boldsymbol{\beta}}
\newcommand{\ga}{\boldsymbol{\gamma}}

\newcommand{\Uad}{\mathbb{U}_{\mathsf{ad}}}
\newcommand{\Uadk}{\mathbb{U}_{\mathsf{ad},k}}

\newcommand{\Xh}{\mathcal{X}_h}

\newcommand{\bff}{\bold{f}}
\newcommand{\bfb}{\bold{b}}
\newcommand{\bfr}{\bold{r}}
\newcommand{\zet}{\boldsymbol{\zeta}}
\newcommand{\bfg}{\bold{g}}

\newcommand{\bv}{\bold{v}}
\newcommand{\varthet}{\boldsymbol{\vartheta}}

\newcommand{\bU}{\bold{U}}
\newcommand{\bu}{\bold{u}}
\newcommand{\x}{\boldsymbol{x}}
 
\usepackage{stmaryrd}

\newcommand{\logLogSlopeTriangle}[5]
{

    \pgfplotsextra
    {
        \pgfkeysgetvalue{/pgfplots/xmin}{\xmin}
        \pgfkeysgetvalue{/pgfplots/xmax}{\xmax}
        \pgfkeysgetvalue{/pgfplots/ymin}{\ymin}
        \pgfkeysgetvalue{/pgfplots/ymax}{\ymax}

        \pgfmathsetmacro{\xArel}{#1}
        \pgfmathsetmacro{\yArel}{#3}
        \pgfmathsetmacro{\xBrel}{#1-#2}
        \pgfmathsetmacro{\yBrel}{\yArel}
        \pgfmathsetmacro{\xCrel}{\xArel}

        \pgfmathsetmacro{\lnxB}{\xmin*(1-(#1-#2))+\xmax*(#1-#2)} 
        \pgfmathsetmacro{\lnxA}{\xmin*(1-#1)+\xmax*#1} 
        \pgfmathsetmacro{\lnyA}{\ymin*(1-#3)+\ymax*#3} 
        \pgfmathsetmacro{\lnyC}{\lnyA+#4*(\lnxA-\lnxB)}
        \pgfmathsetmacro{\yCrel}{\lnyC-\ymin)/(\ymax-\ymin)} 

        \coordinate (A) at (rel axis cs:\xArel,\yArel);
        \coordinate (B) at (rel axis cs:\xBrel,\yBrel);
        \coordinate (C) at (rel axis cs:\xCrel,\yCrel);

        \draw[#5]   (A)-- node[pos=0.5,anchor=north] {1}
                    (B)-- 
                    (C)-- node[pos=0.5,anchor=west] {#4}
                    cycle;
    }
}
\newcommand{\invLogSlopeTriangle}[5]
{

    \pgfplotsextra
    {
        \pgfkeysgetvalue{/pgfplots/xmin}{\xmin}
        \pgfkeysgetvalue{/pgfplots/xmax}{\xmax}
        \pgfkeysgetvalue{/pgfplots/ymin}{\ymin}
        \pgfkeysgetvalue{/pgfplots/ymax}{\ymax}

        \pgfmathsetmacro{\xArel}{#1}
        \pgfmathsetmacro{\yArel}{#3}
        \pgfmathsetmacro{\xBrel}{#1+#2}
        \pgfmathsetmacro{\yBrel}{\yArel}
        \pgfmathsetmacro{\xCrel}{\xArel}
\pgfmathsetmacro{\lnxB}{\xmin*(1-(#1-#2))+\xmax*(#1-#2)}
\pgfmathsetmacro{\lnxA}{\xmin*(1-#1)+\xmax*#1}
\pgfmathsetmacro{\lnyA}{\ymin*(1-#3)+\ymax*#3}
\pgfmathsetmacro{\lnyC}{\lnyA + #4*(\lnxA - \lnxB)}
\pgfmathsetmacro{\yCrel}{(\lnyC-\ymin)/(\ymax-\ymin)}

        \coordinate (A) at (rel axis cs:\xArel,\yArel);
        \coordinate (B) at (rel axis cs:\xBrel,\yBrel);
        \coordinate (C) at (rel axis cs:\xCrel,\yCrel);

\pgfmathtruncatemacro{\absSlopeInt}{abs(#4)}

\draw[#5] (A) -- node[pos=0.5,anchor=south] {1}        
          (B) --
		  (C) -- node[midway, above, allow upside down=false, anchor=east] {\absSlopeInt}
          cycle;

    }
}

\begin{document}

\title[Algebraic flux correction scheme for a parabolic optimal control problem]{Numerical analysis of a stabilized scheme for an optimal control problem governed by a parabolic convection--diffusion equation}

\author{Christos Pervolianakis}
\address{Institut für Mathematik, Friedrich-Schiller-Universität Jena, 07743, Jena, Germany}
\email{christos.pervolianakis@uni-jena.de}


\date{\today}

\begin{abstract}
We consider an optimal control problem on a bounded domain $\Omega\subset\mathbb{R}^2,$ governed by a parabolic convection--diffusion--reaction equation with pointwise control constraints. We follow the optimize--then--discretize approach, in which the state and co-state variables are discretized using the piecewise linear finite element method. For stabilization, we apply the algebraic flux correction method. Temporal discretization is performed using the backward Euler method. The discrete control variable is obtained by projecting the discretized adjoint state onto the set of admissible controls. 
The resulting stabilized fully--discrete scheme is nonlinear and a fixed point argument is used to prove its existence and uniqueness under a mild condition between the time step $k$ and the mesh size $h,$ e.g., $k = \mathcal{O}(h).$ Furthermore, assuming sufficient regularity of the exact solution, we derive error estimates in the $L^{2}$ and energy norms with respect to the spatial variable, and in the $\ell^\infty$ norm with respect to time for the state and co-state variables. For the control variable, we also derive an $L^{2}$-norm error estimate with respect to space and an $\ell^\infty$-norm estimate in time. 
Finally, we present numerical experiments that validate the the order of convergence of the stabilized fully--discrete scheme based on the algebraic flux correction method. We also test the stabilized fully--discrete scheme in optimal control problems that governed by a convection--dominant equation where the solution possesses interior layers.
\end{abstract}


\maketitle

\section{Introduction}
Let $\Omega\subset \R^2,$ a bounded convex polygonal domain with Lipschitz boundary $\partial\Omega$. We shall consider the following optimal control problem:
\begin{align}\label{control_problem}
\min_{u\in \Uad} J(y,u) = \min_{u\in \Uad} \frac{1}{2}\int_0^T\int_\Omega  (y(t) - y_d(t))^2 \,d\x\,dt + \frac{\lambda}{2}\,\int_0^T\int_\Omega u(t)^2\,d\x\,dt,
\end{align}
subject to the state equation governed by the parabolic convection--diffusion--reaction equation,
\begin{equation}\label{conv_diff}
\begin{cases}
y_t - \mu\,\Delta\,y + \bfb \cdot \nabla y + \sigma\,y =  Bu + G, & \text{in }{{\Omega}}\times [0,T],\\
y = 0, & \text{on }\partial {\Omega}\times[0,T],\\
y(\cdot,0)  = y^0, & \text {in }{{\Omega}},
\end{cases}
\end{equation}
with the regularization parameter $\lambda > 0,$ the diffusion coefficient $\mu>0$ and the non-negative reaction coefficient $\sigma \geq 0.$ The operator $B$ is bounded from $L^{2}(0,T;L^{2})$ to $L^{2}(0,T;L^{2}).$ In addition, the right hand side function $G\in L^{2},$ and the velocity field $\bfb := \bfb(\x,t)\in [W^{1}_{\infty}(\Omega)]^2,\,\forall\,t\in[0,T],$ and we assume that the following coercivity condition holds
\begin{align}\label{coercivity_cond}
\sigma - \frac{1}{2}\ddiv\bfb \geq \sigma_0 > 0.
\end{align} 
The set of admissible controls is defined by
\begin{align*}
\Uad = \left\lbrace u\in L^{2}(0,T; L^{2}(\Omega))\,:\, \mathsf{u_a} \leq u(\x,t) \leq \mathsf{u_b},\;\;\text{a.e in }\Omega\times (0,T]\right\rbrace,
\end{align*}
with $\mathsf{u_a},\,\mathsf{u_b}\in\mathbb{R}$ such that $\mathsf{u_a} < \mathsf{u_b}.$ \par

Many physical or chemical optimization processes can be modeled by optimal control problems that are governed by (non--)linear--convection--diffusion--reaction equations, see, e.g., the optimal control of Keller--Segel equations \cite{ryu2001}, where optimal control strategies for chemotactic movement are presented. For more applications and a systematic introduction of finite elements for optimal control problems, we refer to the book \cite{hinze2009book}. \par

When discretizing optimal control problems the "optimize-then-discretize" and the "discretize-then-optimize" approaches are usually considered, see, e.g., \cite{hinze2005}. More specifically, in the "optimize-then-discretize" approach, the optimality conditions are derived in continuous level while in "discretize-then-optimize" approach, the optimal control problem is first discretized by a numerical method, e.g., the finite element method, and then, the optimality conditions are derived. In this work, we focus on the discretize-then-optimize-approach, since we will use a nonlinear stabilization method, and more specifically the algebraic flux correction method. 

There are many works considering optimal control problems that are governed by elliptic or parabolic equations, where their discretization is based on the finite element method. To name a few, for optimal control problems governed by elliptic equations, we refer to \cite{arada2002,casas2005,meyer2004,becker2007b,fu2009} and for parabolic equations \cite{becker2007,meidner2008,meidner2008b,springer2014,vondaniels2015,meidner2011, neitzel2012,zhang2023}. For the optimal control problems governed by elliptic equations, an error estimate of order $\mathcal{O}(h)$ was shown for the control variable using constant discretizations. The estimate can be refined to $\mathcal{O}(h^{3/2})$ in case where the control variable is discretized by (bi--/tri--)linear $H^{1}-$conforming element, see, e.g., \cite{becker2007b,yan2009}. When the control variable is not discretized and implicitly utilize the optimality conditions by discretizing the state and the co-state variables, then a $\mathcal{O}(h^2)$ can be proved, see \cite{hinze2005}. \par

It is known that, convection--diffusion equations may possess layers, i.e., small regions where the solution has a large gradient, and therefore the standard finite element method may fail to approximate the exact solution due to these layers. Consequently, nonlinear discretizations have been developed which preserve the maximum principle of the solution and accurately determine the position of the layers, e.g., the streamline--upwind Petrov--Galerkin method \cite{brooks1982}, the so-called SOLD methods \cite{john2007}, the least-square stabilization of the gradient jumps across element edges for an optimal control problem governed by an elliptic convection--diffusion equation \cite{yan2009} as well as the Algebraic Flux Correction schemes (AFC), which have gained a systematical attention in the last two decades, see, e.g., \cite{kuzmin2005,john2021,kuzmin2010book,barrenechea2025, barrenechea2016, barrenechea2018, barrenechea2017b, chatzipantelidis2022, barrenechea2024, jha2021} and the references therein. \par

Many stabilized numerical methods based on the finite element method, has been proposed for the optimal control problem that is governed by an elliptic convection--diffusion--reaction. These include the streamline upwind Petrov
Galerkin (SUPG) method, see, \cite{collis2002}, where the authors compare two different approaches for the discretization of the optimal control problem; the "optimize-then-discretize" and the "discretize-then-optimize". In the former, the necessary optimality conditions are derived at the continuous level, while in the latter, the discretization method is applied where afterwards one has to solve a finite dimensional optimization problem. The authors showed that in the context of the SUPG, these two approaches are not generally the same. Further, the authors in \cite{becker2007a} proposed and analyzed a numerical scheme where the optimality conditions were discretized based on local projections, the so--called (LPS) method, a method that uses symmetrical penalty terms, see, e.g, \cite{guermond1999} in the context of the convection--diffusion equation. Contributions on the optimal control problem that is governed by an elliptic convection--diffusion equation have been also proposed using an edge stabilization (ES) method, see, e.g., \cite{yan2009}; also in the context of the convection--diffusion equation in \cite{burman2004}. We also refer to the survey \cite{weng2015}. Recently, the AFC schemes have been proposed in elliptic optimal control problems, see, e.g., \cite{baumgartner2022,baumgartner2025}. \par

In this work, our aim is to develop and analyze a numerical scheme for \eqref{control_problem}--\eqref{conv_diff}, whose solution satisfies the discrete analogue of the maximum principle. To enforce this property, we use the  algebraic flux correction method.

\subsection{Optimality conditions}

Throughout the paper, we use the standard notation for Lebesgue  and Sobolev spaces, namely we denote $W^{m}_p=W^{m}_p(\Omega)$, $H^{m}=W^{m}_{2}$, $L^{p}=L^{p}(\Omega)$, and  with $\|\cdot\|_{m,p}=\|\cdot\|_{W^{m}_p}$, $\|\cdot\|_{m}=\|\cdot\|_{H^{m}}$,  $\|\cdot\|_{L^{p}}=\|\cdot\|_{L^{p}(\Omega)}$, for $m\in\mathbb{N}$ and $p\in[1,\infty]$, the corresponding norms. In addition, we let $L^{r}(0,T; W^{m}_p(\Omega)),$ be the Banach space of all $L^{r}$ integral functions defined from $[0,T]$ into $W^{m}_{p},$ equipped with norm
\begin{align*}
\|v\|_{L^{r}(0,T; W^{m}_p)} = \left(\int_0^T \|v\|_{m,p}^r\,ds\right)^{1/r},\;\;1\leq p < \infty.
\end{align*}

The basis for the methods studied is the variational formulation of the model problem, to find functions $(y,u)\in L^{\infty}(0,T;H^{1}_0)\times \Uad$ such that
\begin{align}\label{weak_u}
\min_{u\in \Uad} J(y,u) = \min_{u\in \Uad} \frac{1}{2}\int_0^T \left( \|y(t) - y_d(t)\|_{L^{2}}^2 + \lambda\,\|u(t)\|^2_{L^{2}}\right)\,dt,
\end{align}
subject to the weak formulation to the state equation,
\begin{align}
(y_t(t), v) + \alpha(y(t), v) = (Bu(t) + G(t),v), \;\;\;\forall\, v\in H^{1}_0,\label{weak_y}
\end{align}
where the bilinear form $\alpha(\cdot,\cdot)\,:\,H^{1}_0\times H^{1}_0 \to \mathbb{R}$ is defined as $\alpha(v,w) := (\mu\nabla v, \nabla w) + (\bfb\cdot\nabla v,w) + (\sigma v, w)$ for all $v,\,w\in H^{1}_0.$ The bilinear for is bounded and coercive in $H^{1}_0,$ due to \eqref{coercivity_cond}. We define the energy norm, $\tribar v \tribar^2 := \mu\|\nabla v\|_{L^{2}}^2 + \sigma\|v\|_{L^{2}}^2,\,v\in H^{1},$ in which we will derive a priori error estimates. We note that, $\|\nabla v\|_{L^{2}} \leq \mu^{-1/2}\tribar v \tribar,\,v\in H^{1}.$

It is known, see \cite[Chapter III, Theorem 2.1]{lions1971}, that the optimal control problem \eqref{weak_u}--\eqref{weak_y} has a unique solution $(y,u),$ if and only if there exists a co-state variable $\overline p$ such that $(\overline y,\overline p, \overline u)\in H^{1}_0 \times H^{1}_0 \times \Uad,$ satisfies the following optimality conditions for $t\in[0,T],$
\begin{alignat}{2}
(\overline y_t , v) + \alpha(\overline y, v) 
  &\;=\; (B\overline{u} + G,v) 
  &\quad& \forall\, v\in H^{1}_0, \label{optimal_weak_y} \\
- (\overline p_t, v) + \alpha(v, \overline p) 
  &\;=\; (\overline y - y_d ,v) 
  &\quad& \forall\, v\in H^{1}_0, \label{optimal_weak_p} \\
(\lambda\, \overline u +  B^*\overline p ,  \widehat u - \overline u) 
  &\;\geq\; 0 
  &\quad& \forall\,\widehat u\in\Uad, \label{optimal_weak_cond} \\
\overline y(\x,0) & = y_0(\x),\,\overline p(\x,T) = 0, 
  &\quad&  \notag
\end{alignat}
where $B^*$ is the adjoint operator of $B.$ We now introduce the projection operator $\Pi_{[\mathsf{a},\mathsf{b}]}\,:\,L^{1} \to \Uad$ as 
\begin{align}\label{operator_Pi}
\Pi_{[\mathsf{a},\mathsf{b}]}(v) : = \max\left\lbrace \mathsf{a}, \min\{\mathsf{b},v\}\right\rbrace,\;\;\;\text{a.e in } \Omega.
\end{align}
In view of \cite[Lemma 2.2]{meyer2004}, we can express the \eqref{optimal_weak_cond}, in the following equivalent form,
\begin{align}\label{optimal_weak_u2}
\overline u(t) := \Pi_{[\mathsf{u_a},\mathsf{u_b}]}(-\lambda^{-1}B^* \overline p(t)),\;\;\;t\in[0,T].
\end{align}

\subsection{Finite element method}

The finite element methods studied are based on triangulations $\Th =\bigcup_{i=1}^{\N}K_i,\,\N\geq1,$ of $\Omega,$ where $h = \max_{1\leq i\leq \N}\text{diam}(K_i).$ We use the finite element spaces
\begin{align}
\Sh &: = \left\lbrace \chi \in\mathcal{C}(\overline{{\Omega}})\,:\,\chi|_{K} \in \mathbb{P}_1,\;\forall\;K\in \Th,\,\;\text{such that}\;\,\chi = 0\;\;\text{on}\;\;\partial\Omega\right\rbrace,\label{fem_space}
\end{align} 
The semi--discrete approximation of the variational problem \eqref{weak_u}--\eqref{weak_y}, may be written as follows: Find functions $(y_h,u_h)\in \Sh \times \Uad$, with $y_h(0) = y^0_h\in \Sh$, such that
\begin{align}\label{semi_weak_u}
\min_{u_h\in \Uad} J(y_h,u_h) = \min_{u_h\in \Uad} \frac{1}{2}\int_0^T \left( \|y_h - y_d\|_{L^{2}}^2 + \|u_h\|^2_{L^{2}}\right)\,dt,
\end{align}
subject to 
\begin{align}
(y_{h,t}, \chi) + \alpha(y_h, \chi)  = (Bu_h + G, \chi), \;\;\;\forall\, \chi\in \Sh.\label{semi_weak_y}
\end{align}
The optimal control problem \eqref{semi_weak_u}--\eqref{semi_weak_y} has a unique solution $(y_h, u_h)\in \Sh \times \Uad$ if and only if there exists a co-state variable $p_h\in \Sh$ such that the following optimality conditions are satisfied.
\begin{alignat}{2}
(y_{h,t}, \chi) + \alpha(y_h, \chi) 
  &\;=\; (Bu_h + G,\chi) 
  &\quad& \forall\, \chi\in \Sh, \label{optimal_semi_weak_y} \\
- (p_{h,t}, \chi) + \alpha(\chi,  p_h) 
  &\;=\; (y_h - y_d, \chi) 
  &\quad& \forall\, \chi\in \Sh, \label{optimal_semi_weak_p} \\
(\lambda\,u_h +  B^*p_h, \psi - u_h) 
  &\;\geq\; 0 
  &\quad& \forall\,\psi\in\Uad, \label{optimal_semi_weak_u} \\
y_h(0) & = y_h^0,\;\;p_h(T) = 0. 
  &\quad& \label{optimal_semi_weak_cond}
\end{alignat}

Similarly to \eqref{optimal_weak_u2}, one can express \eqref{optimal_semi_weak_u} by using \eqref{operator_Pi}, as
\begin{align}\label{optimal_semi_weak_u2}
u_h(t) := \Pi_{[\mathsf{u_a},\mathsf{u_b}]}(-\lambda^{-1} B^*p_h(t)),\;\;\;t\in[0,T].
\end{align}
The equation \eqref{optimal_semi_weak_u2} means that the control variable $u_h\in\Uad$ is the projection of the finite element function onto the admissible space $\Uad.$ Note that $u_h$ in general does not belong to $\Sh,$ see, e.g., its definition \eqref{optimal_semi_weak_u2}. While it is not directly discretized, it depends on $h,$ indirectly due to its definition, thus, we use the subscript $h$ on it. \par

In this work, we will stabilize the semi--discrete scheme \eqref{optimal_semi_weak_y}--\eqref{optimal_semi_weak_u2} using the algebraic flux correction method. To accomplish this, we need first to write the matrix formulation of the latter semi--discrete scheme. Let $\Zh = \lbrace Z_j\rbrace_{j=1}^{\N}$ be the set of nodes in the triangulation $\Th$ and $\lbrace \phi_j \rbrace_{j=1}^{\N}\subset \Sh$ the corresponding nodal basis, with $\phi_j(Z_i)=\delta_{ij}.$ Then, we may write $y_h(t)=\sum_{j=1}^{\N}\alpha_j(t)\phi_j,\,p_h(t)=\sum_{j=1}^{\N}\beta_j(t)\phi_j$, with $y_h^0=\sum_{j=1}^{\N}\alpha_j^0\phi_j,\,p_h^{\NT}=\sum_{j=1}^{\N}\beta_j^{\NT}\phi_j.$ The semi--discrete problem  \eqref{optimal_semi_weak_y}--\eqref{optimal_semi_weak_cond} can then be expressed, with 
$\al = \al(t),\,\be = \be(t)$ and $\al = (\alpha_1, \dots, \alpha_{\N})^T,\,\be = (\beta_1, \dots, \beta_{\N})^T,$ as follows,
\begin{align}
\M\al^\prime(t) + (\mu\,\S + \Q + \sigma\,\M ) \al(t) & = \bfr(Bu_h) + \bfr(G), \;\qquad\text{ for } t\in[0,T], \text{ with }\al(0)=\al^0,\label{matrix_y}\\
-\M\be^\prime(t) + (\mu\,\S - \Q + \sigma\,\M ) \be(t) & = \M\,\al(t) - \bfr(y_d), \qquad\,\text{ for } t\in[0,T], \text{ with }\be(T)=\mathbf{0},\label{matrix_p}
\end{align}
where $\mathbf{0}$ the zero vector and for an $s\in L^{2},$ the vector $\bfr(s) = (r_i(s)),$ with elements $r_i(s) = (s, \phi_i),\,i=1,\ldots,\N.$

 The matrix $\M=(m_{ij})$ with elements  $m_{ij}=(\phi_i,\phi_j)$ and $\S=(s_{ij})$ with elements $s_{ij}=(\nabla\phi_i,\nabla\phi_j)$ are the usual mass and stiffness matrix, respectively. The matrix due to the convection term in \eqref{optimal_semi_weak_y} and in \eqref{optimal_semi_weak_p} is $\Q=(\tau_{ij})$, and its elements are defined as follows,
\begin{align}\label{T_def}
\tau_{ij} = (\bfb \cdot \nabla \phi_j ,  \phi_i), \quad \text{for } i,j=1,\dots,\N.
\end{align}
The mass matrix $\M$ and the stiffness matrix $\S$, are both symmetric and positive definite while on the other hand the matrix due to the convection term is not symmetric. \par

\subsection{Contributions and outline of the paper}

In this paper, our aim is to analyze a fully--discrete scheme for the approximation of the optimal control problem governed by a parabolic convection--diffusion equation \eqref{control_problem}--\eqref{conv_diff}, by enforcing it to satisfy the discrete maximum principle. We follow the optimize--then--discretize approach, where for the state and the co-state variables, we consider the piecewise finite element method alongside with the algebraic flux correction method for their stabilization in space and for the temporal discretization, we use the backward Euler method. For the control discretization, we use the variational discretization proposed in \cite{hinze2005}, where the control variable is not directly discretized by finite element method, rather than, it is a projection of a finite element function that belongs to $\Uad.$ \par

The resulting stabilized fully--discrete scheme with the AFC method is nonlinear, thus, we employ a fixed point argument to prove, under a mild condition between the time step and the mesh step $k,\,h,$ respectively, e.g., for $k = \mathcal{O}(h),$ its existence and uniqueness.
\par

The fully--discrete schemes we consider approximate $\overline y^n$ by $Y^n \in \Sh$ where $\overline y^n=\overline y (\cdot,t^n)$, solution of \eqref{optimal_weak_y} and $t^n=nk$, $n=0,\dots, \NT$ and $\NT\in\mathbb{N}$, $\NT\ge1$, $k=T/\NT$. Similarly, approximate $\overline p^n$ by $P^n \in \Sh$ where $\overline p^n=\overline p(\cdot,t^n)$ solution of \eqref{optimal_weak_p} and $\overline u^n$ by $U^n \in \Uad$ where $\overline u^n=\overline u (\cdot,t^n)$, solution of \eqref{optimal_weak_cond}. Under sufficient regular solution, namely $\overline y,\,\overline p \in L^{\infty}(0,T;H^{2})\cap H^{2}(0,T;L^{2})$ with $\overline y_t,\,\overline p_t \in L^{2}(0,T;H^{2})\cap L^{\infty}(0,T;H^{1}_0),\;\overline u\in L^{2}(0,T;H^{1})\cap H^{1}(0,T;L^{2}),$ and small $k,\,h,$ with $k=\mathcal{O}(h),$ we derive error estimates for the AFC scheme 
\begin{align*}
\max_{0\leq n\leq \NT} \left( \|Y^n-\overline y^n\|_{L^{2}} + \|P^n-\overline p^n\|_{L^{2}} + \|U^n-\overline u^n\|_{L^{2}} \right) + k\sum_{n=0}^{\NT}\left( \tribar Y^n-\overline y^n \tribar^2 + \tribar P^n-\overline p^n \tribar^2 \right) & \leq C (k+h^2),\\
\max_{0\leq n\leq \NT} \left( \tribar Y^n-\overline y^n \tribar + \tribar P^n-\overline p^n \tribar \right)  & \leq C (k+h),
\end{align*}
where the constant $C$ is independent of $k,\,h,$ and depends on the norms of the solution as well as the parameters $\lambda,\,\mu^\nu,\,\nu<0.$ \par
The paper is organized as follows: In Section \ref{section:AFC}, we introduce the semi--discrete AFC scheme for the discretization of the optimal control problem governed by a parabolic convection--diffusion equation. Further, we recall some auxiliary results for the stabilization terms from \cite{chatzipantelidis2022}, that we will employ in the analysis. In Section \ref{section:fully_discrete}, we discretize semi--discrete scheme using the backward Euler method. For a smooth solution of \eqref{optimal_weak_y}--\eqref{optimal_weak_u2}, we derive error estimates in $L^{2}$ and in energy norm in spatial variable and $\ell^\infty$ in time. In Section \ref{section:fixed_point}, we define a fixed point scheme for the fully--discrete scheme and we prove that the resulting nonlinear fully--discrete scheme has a unique solution. Finally, in Section \ref{section:numerical_results}, we test the stabilized fully--discrete scheme on optimal control problems governed by a convection-dominated equation, where the solution exhibits layers, and we compute the experimental order of convergence which validates our theoretical results.

\section{Algebraic flux correction}\label{section:AFC}

\subsection{Mesh assumptions}

We consider a  family of regular triangulations $\Th$ of the domain $\overline\Omega\subset\R^2$. We will assume that the family $\Th$ satisfies the following assumption.
 
\begin{assumption}\label{mesh-assumption}
Let $\Th =\bigcup_{i=1}^{\N}K_i,\,\N\geq1,$ be a  family of regular triangulations  of $\overline\Omega$ such that any edge of any $K$ is either a subset of the boundary $\partial\Omega$ or an edge of another $K \in \Th$, and in addition
\begin{enumerate}
\item  $\Th$ is shape regular, i.e, there exists a constant $\varpi >0,$ independent of $K$ and $\Th,$ such that 
\begin{equation}\label{shape_regularity}
\frac{h_K}{\varrho_K} \leq \varpi,\quad \forall K\in\Th,
\end{equation}
where $h_K$ the longest edge of $K,\,\varrho_K=\text{diam}(B_K)$, and $B_K$ is the inscribed ball in $K$.
\item The  family of triangulations $\Th$ is quasi-uniform, i.e., there exists constant $\varrho>0$ such that
\begin{align}\label{quasi-uniformity}
\frac{\max_{K\in\Th}h_K}{\min_{K\in\Th}h_K} \leq \varrho,\quad\forall K\in\Th.
\end{align}
\item
We assume that $\Th$ satisfies an acute condition, i.e., all interior angles of a triangle $K\in\Th$ are less than $\pi/2$. 
\end{enumerate}
\end{assumption}

\begin{remark}
The third part of the Assumption \ref{mesh-assumption} implies that $s_{ij}\le0,\,j\neq i$ and $s_{ii}>0,$ see, e.g., {\normalfont\cite{draganescu2004}}. A motivation why we impose such an assumption, it will explained in Remark \ref{remark:optimal2}.
\end{remark}

Since $\Th$ satisfies \eqref{quasi-uniformity}, we have for all $\chi\in \Sh,$ cf., e.g., \cite[Chapter 4]{brenner2008}, 
\begin{equation}\label{eq:inverse_estimate}
\|\chi\|_{L^{\infty}} + \|\nabla \chi\|_{L^{2}}  \le Ch^{-1}\|\chi\|_{L^{2}}
\text{ and } \|\nabla \chi\|_{L^{\infty}}  \le Ch^{-1}\|\nabla\chi\|_{L^{2}}.
\end{equation}

\begin{figure}
\centering
\includegraphics[scale=0.7]{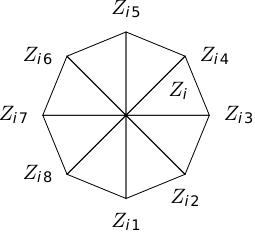} 
\caption{The subdomain $\omega_i$ of the triangulation $\Th.$}\label{fig:patches}
\end{figure}

Let $\Nh$ be the the indices of all the nodes of $\Th$, i.e., $\Nh:=\{ i: Z_i \text{ a node of the triangulation } \Th\}$ which  can be splitted into the indices of the internal nodes, $\Nh^0$,  and  the indices of the nodes  on the boundary  $\partial\Omega$,  $\Nhb$, i.e.  $\Nh := \Nh^0\cup \Nhb$.  Also let $\Ehh$ be the set of all edges of the triangulation $ \Th.$ Similarly, we split this  set into  the internal edges, $\Ehh^0$  and the edges on the boundary $\partial\Omega$, $\Ehb$, i.e. $\Ehh := \Ehh^{0}\cup\Ehb$. We denote  $\omega_e$ the collection of triangles with a common edge $e\in\Ehh$,  see Fig. \ref{fig:patches}, and $\omega_i$, $i\in \Nh$, the collection of triangles with a common vertex $Z_i$,  i.e. $\omega_i = \cup_{Z_i\in K}\overline{K},$ see Fig. \ref{fig:patches}. The sets $\Nh(\omega)$ and $\Ehh(\omega)$ contain the vertices and  the edges, respectively, of a  subset of $\omega\subset\Th$ and $\Zh^i$ the set of nodes adjacent to $Z_i$, $\Zh^i:=\{ j: Z_j\in \Zh, \text{adjacent to }Z_i\}$. Using the fact that $\Th$ is shape regular, there exists a constant $C_\varpi$, independent of $h$,  such that the number of vertices in $\Zh^i$  is less than  $C_\varpi$, for $i=1,\dots, \N$. Also  $e_{ij}\in\Ehh$ denotes an edge of $\Th$ with endpoints $Z_i$, $Z_j\in \Zh$.

\subsection{Auxiliary results}

Recall that the bilinear form  $(\cdot,\cdot)_h,$ see, e.g., \cite[Chapter 15]{thomee2006}, is an inner product in $\Sh$ that approximates $(\cdot,\cdot)$ and is defined by
\begin{equation}\label{quadrature}
(\psi,\chi)_h = \sum_{K\in\Th}Q_h^K(\psi\chi),\ \text{ with }Q_h^K(g) = \frac{1}{3}|K|\sum_{j=1}^3g(Z_j^K)\approx \int_K g\,dx,
\end{equation}
with $\{Z_j^K\}_{j=1}^3$ the vertices of a triangle $K\in\Th.$ For the inner product $(\cdot,\cdot)_h$ introduced in \eqref{quadrature}, the following holds.

It is known that $(\cdot,\cdot)_h$ induces an equivalent norm to $\|\cdot\|$ on $\Sh$. In particular, exists constants $C_1,\,C_2$ independent on $h$, such that
\begin{equation}\label{mass_lump_equivalence}
C_1\|\chi\|_h  \leq \|\chi\|  \leq C_2\|\chi\|_{h},\ \text{ with }\ \|\chi\|_h = (\chi, \chi)_h^{1/2},\quad\forall\chi\in \Sh.
\end{equation}

\begin{lemma}{\normalfont{\cite[Lemma 2.3]{chatzipantelidis2012}}}\label{lemma:mass_lump_error}
Let $\varepsilon_{h}(\chi ,\psi): = (\chi , \psi) - (\chi , \psi)_h$. Then,
\begin{equation}
\vert \varepsilon_{h}(\chi ,\psi)\vert\leq Ch^{i+j}\|\nabla^i\chi \|_{L^{2}}\|\nabla^j\psi\|_{L^{2}},\;\;\;\forall\,\chi ,\psi\in \Sh,\;\;\;\text{and}\;\;\;i,j=0,1,\nonumber
\end{equation}
where the constant $C$ is independent of $h.$
\end{lemma}

We recall a modified $L^{2}-$projection operator $\widetilde P_h \,:\,L^{2} \to \Sh,$ see, e.g., \cite[Lemma 2.2]{bartels2022}, defined by
\begin{align}\label{L2_modf_projection}
(\widetilde P_hv, \chi)_h =  ( v,  \chi),\quad\forall \chi\in \Sh.
\end{align}
For the latter modified  $L^{2}-$projection, similar to $P_h$ properties can be proved, see, e.g., \cite[Lemma 2.2]{bartels2022}. More specifically, there exists a constant $C,$ such that for a $v\in H^{1}_0,$ holds that
\begin{align*}
\|\nabla \widetilde P_hv\|_{L^{2}} + h^{-1}\|v - \widetilde P_hv\|_{L^{2}} \leq C\|\nabla v\|_{L^{2}}.
\end{align*}

\subsection{A stabilized semi--discrete scheme}
 
We replace the composite mass matrix $\M$ by the corresponding lumped mass matrix $\M_L$ and the negative off-diagonal elements of $\Q$ are cured by an artificial diffusion operator $\D=(d_{ij})$ so that the off-diagonal elements of $\Q+\D$ be non-positive, elementwise. Similarly, we also define an artificial diffusion $\widehat \D = (\widehat d_{ij})$ so that the off-diagonal elements $\Q+\widehat\D$ be non-negative, elementwise. 

\begin{remark}\label{remark:acute}
In order to ensure that the matrix $\mu\S + \Q + \D$ in \eqref{matrix_y}, as well as the matrix $\mu\S - \Q - \widehat \D$ in \eqref{matrix_p}, each have non-positive off-diagonal entries and positive diagonal entries, it suffices to construct the matrices $\D$ and $\widehat \D$ such that the matrices $\Q + \D$ and $-\Q - \widehat \D$ possess these properties, since $\Th$ satisfies an acute condition, see Assumption \ref{mesh-assumption}. We note that the third condition of Assumption \ref{mesh-assumption} is not a necessary condition to preserve the discrete maximum principle.
\end{remark}

Since, we would like not to add mass in our scheme, the matrices $\D,\,\widehat \D,$ must be symmetric with zero row and column sums, cf. \cite{kuzmin2002},  which is true if  $\D=(d_{ij})_{i,j=1}^{\N}$ is defined by
\begin{equation}\label{D_def}
d_{ij} : = \min\{ - \tau_{ij},0, -\tau_{ji}\} = d_{ji} \le 0,\quad\forall j\neq i\ 
\text{ and }\ 
d_{ii} : = -\sum_{j\neq i}d_{ij}.
\end{equation}
In a similar way, we define the matrix $\widehat \D=(\widehat d_{ij})_{i,j=1}^{\N},$ with elements,
\begin{equation}\label{D_widehat_def}
\widehat d_{ij} : = \max\{ - \tau_{ij},0, - \tau_{ji}\}= \widehat d_{ji} \ge 0,\quad\forall j\neq i\ 
\text{ and }\ 
\widehat d_{ii} : = -\sum_{j\neq i}\widehat d_{ij}.
\end{equation}

In view of the above definitions, we can derive for a $\bfb \in L^{\infty}$, \eqref{T_def} and the fact that the triangulation $\Th$ is shape regular, i.e., \eqref{shape_regularity},  that there exists a constant $C(\varpi)>0,$ independent of $h$, such that
\begin{equation}\label{est_dij_2D_L_infty}
\begin{aligned}
\vert d_{ij}\vert  \leq \vert\tau_{ij}\vert + \vert\tau_{ji}\vert & \leq \|\bfb\|_{L^{\infty}}\sum_{K\in \omega_i} (\|\nabla \phi_i\|_{L^{2}(K)}\|\phi_j\|_{L^{2}(K)} + \|\nabla\phi_j\|_{L^{2}(K)}\|\phi_i\|_{L^{2}(K)})\\
& \leq C \|\bfb\|_{L^{\infty}} \sum_{K\in \omega_i}h_K \leq C(\varpi)\,\|\bfb\|_{L^{\infty}}\,h.
\end{aligned}
\end{equation}
A similar estimate holds also for \eqref{D_widehat_def}, i.e., 
\begin{align}\label{est_dij_2D_L_infty_hat}
\vert \widehat d_{ij}\vert \leq C(\varpi)\,\|\bfb\|_{L^{\infty}}\,h.
\end{align}

The resulting system for the approximation of \eqref{control_problem}--\eqref{conv_diff} is expressed as follows, we seek $\al(t),\,\be(t),\,\ga(t)\in\Rn$ such that, for $t\in[0,T]$,
\begin{align}
\M_L\al^\prime(t) + (\mu\,\S + \Q + \D + \sigma\,\M_L)  \al(t) & = \bfr(Bu_h) + \bfr(G), \qquad\text{ for } t\in[0,T], \text{ with }\al(0)=\al^0,\label{LOW_matrix_y}\\
-\M_L \be^\prime(t) + (\mu\,\S - \Q - \widehat \D  + \sigma\,\M_L) \be(t) & = \M\,\al(t) - \bfr(y_d), \;\;\;\quad\text{ for } t\in[0,T], \text{ with }\be(T)=\mathbf{0}.\label{LOW_matrix_p}
\end{align}
We will use the latter ODE as a basis for the construction of the algebraic flux correction scheme (AFC scheme). This procedure is commonly used in scalar conservation laws, see, e.g., \cite[Chapter 4]{kuzmin2010book}, \cite{kuzmin2004}. To formulate the variational formulation of \eqref{LOW_matrix_y}--\eqref{LOW_matrix_p}, recall the stabilization term that also used in \cite{barrenechea2016}.\par

Define for a function $s\in \Sh,$ its nodal values as $s_i = s(Z_i),\;i=1,\ldots,\N.$ The bilinear form $d_{\D,h}(\cdot,\cdot):\CC \times \CC \to \R,$  is defined by
\begin{equation}\label{stab_term}
d_{\D,h}(v,z) := \sum_{i,j=1}^{\N}\,d_{ij}(v_i - v_j)z_i = \sum_{i<j}\,d_{ij}(v_i - v_j)(z_i - z_j),\quad\forall v,\,z\in{\CC},
\end{equation}
where the last equality is due to the symmetry of matrix $\D,$ see, e.g., \cite{barrenechea2018}. In a similar way, we can define the stabilization term due to the artificial diffusion operator $\widehat \D,$ the $d_{\widehat\D,h}.$\par

We define the bilinear form $\alpha_h(\cdot,\cdot)\,:\,\Sh\times\Sh \to\mathbb{R},$ as $\alpha_h(\psi,\chi):= (\mu\nabla \psi, \nabla \chi) + (\bfb\cdot \nabla \psi, \chi) + (\sigma\psi,\chi)_h$ for all $\psi,\,\chi\in\Sh.$ We define the mesh--dependent energy norm, $\tribar v \tribar^2_h := \mu\|\nabla v\|_{L^{2}}^2 + \sigma\|v\|_h^2,\,v\in H^{1}.$ The latter bilinear form is bounded and coercive on $\Sh,$ see the condition \eqref{coercivity_cond} as well as Lemma \ref{lemma:mass_lump_error}.
\begin{remark}\label{remark:h_norm}
Since the norm $\|\cdot\|_h$ is equivalent to $\|\cdot\|_{L^{2}}$ on $\Sh,$ see \eqref{mass_lump_equivalence}, it follows that $\tribar  \cdot \tribar_h$ is equivalent $\tribar \cdot \tribar$ on $\Sh.$ In particular, it holds 
\begin{align}\label{h_norm_equivalence}
\widetilde C_1\tribar \chi \tribar_h   \leq \tribar \chi \tribar \leq \widetilde C_2 \tribar \chi \tribar_h,\quad\forall\chi\in \Sh \text{ with } \widetilde C_1 = \min\{1, C_1\},\,\widetilde C_2 = \max\{1, C_2\}.
\end{align}
\end{remark}
Following the analysis in \cite{barrenechea2016}, the low order scheme \eqref{LOW_matrix_y}--\eqref{LOW_matrix_p}, may be formulated as follows.
We seek $(y_h,p_h,u_h)\in \Sh \times \Sh \times \Uad$, such that
\begin{align}
(y_{h,t}, \chi)_h + \alpha_h( y_h , \chi) - \widehat d_{\D,h}(y_h;y_h,\chi)  & = (Bu_h  + G, \chi),\label{optimal_LOW_semi_weak_y}
\end{align}
for all $\chi \in\Sh$ with $y_h(0)=y_h^0,$ and
\begin{align}
- (p_{h,t}, \chi)_h + \alpha_h(\chi,p_h) +   d_{\widehat\D,h}(p_h;p_h,\chi) & = (y_h - y_d,\chi),\label{optimal_LOW_semi_weak_p}
\end{align}
for all $\chi \in\Sh$ with $p_h(T)=0.$ The control variable is defined as
\begin{align}
(\lambda\, u_h +  B^*p_h, \psi - u_h) & \geq 0, \;\;\;\psi\in\Uad.\label{optimal_LOW_semi_weak_u}
\end{align}
Similarly to \eqref{optimal_weak_u2}, one can express \eqref{optimal_LOW_semi_weak_u} by using \eqref{operator_Pi}, as
\begin{align}\label{optimal_LOW_semi_weak_u2}
u_h(t) := \Pi_{[\mathsf{u_a},\mathsf{u_b}]}(-\lambda^{-1} B^*p_h(t)),\;\;\;t\in[0,T].
\end{align}

It is known, that the latter scheme may pollute the accuracy of the approximation. A common remedy is to return to the system the elements that are harmless to the discrete maximum principle, see, e.g., \cite{kuzmin2010book, kuzmin2004} and the references therein. This approach is called algebraic flux correction scheme, which involves the decomposition of this error into internodal fluxes, which can be used to restore the optimal accuracy in regions where the solution is well resolved and no modifications of the standard FEM are required.  Algebraic flux corrections schemes were proposed in \cite{kuzmin2002} and have attained a lot of attention by many authors, see, e.g.,\cite{kuzmin2005, john2021, kuzmin2010book, barrenechea2025, barrenechea2016, barrenechea2018, barrenechea2017b, chatzipantelidis2022, barrenechea2024, kuzmin2004, kuzmin2002}
.\par

In what follows, we will describe the application of the AFC method to the low-order scheme \eqref{optimal_LOW_semi_weak_y}--\eqref{optimal_LOW_semi_weak_u2}, by recalling arguments that have already presented extensively in the literature in the context of the convection--diffusion equations, see, e.g., \cite{chatzipantelidis2022,barrenechea2016, barrenechea2018, barrenechea2017b,barrenechea2025}.\par

Let $\bff=(f_1,\ldots,f_\N)^T$  denote  the error of inserting the artificial diffusion operator $\D$ in \eqref{LOW_matrix_y}, i.e., $\bff(\al) =   \D\al.$ Similarly, the function $\widehat\bff=(\widehat f_1,\ldots,\widehat f_\N)^T$ denote the error of inserting the artificial diffusion operator $\widehat \D$ in \eqref{LOW_matrix_p}, i.e., $\widehat \bff(\be) =  -\widehat \D\be.$ 
In addition, we define the functions $\bfg=(g_1,\ldots,g_\N)^T,\,\widehat\bfg=(\widehat g_1,\ldots,\widehat g_\N)^T,$ such that $\bfg(\al) = (\M_L - \M)\al,\,\widehat\bfg(\be) = (\M_L - \M)\be,$ where $\bfg(\al^\prime),\,\widehat\bfg(\be^\prime)$
denotes the error due to the replacement of the mass matrix $\M$ by its mass lumped $\M_L.$ Here, $\al',\,\be^\prime$ denotes the derivative with respect to temporal variable of $\al = \al(t),\,\be = \be(t),$ respectively. Note that, we use four different correction factors, two for the error due to the matrices $\D,\,\widehat \D,$ and two for the error due to the mass lumping. \par
Using the zero row sum property of matrices $\D,\,\widehat \D$, cf. \eqref{D_def},\eqref{D_widehat_def}, we can show, see e.g., \cite{kuzmin2002}, that the residual admits a conservative decomposition into internodal fluxes,
\begin{equation}\label{internodal_fluxes_f}
\bff_i = \sum_{j\neq i}\mathsf{f}_{ij},\quad \mathsf{f}_{ji} = -\mathsf{f}_{ij}\;\;\text{and}\;\;\widehat\bff_i = \sum_{j\neq i}\widehat{\mathsf{f}}_{ij},\quad \widehat{\mathsf{f}}_{ji} = -\widehat{\mathsf{f}}_{ij},\quad i = 1,\dots,\N.
\end{equation}
Also, the differences $(\M_L - \M)\al,\,(\M_L - \M)\be,$
\begin{equation}\label{internodal_fluxes_g}
\bfg_i = \sum_{j\neq i}\mathsf{g}_{ij},\quad \mathsf{g}_{ji} = -\mathsf{g}_{ij}\;\;\text{and}\;\;\widehat\bfg_i = \sum_{j\neq i}\widehat{\mathsf{g}}_{ij},\quad \widehat{\mathsf{g}}_{ji} = -\widehat{\mathsf{g}}_{ij},\quad i = 1,\dots,\N,
\end{equation}
where the amount of mass transported by the raw \textit{antidiffusive flux} is given by
\begin{alignat}{3}
\mathsf{f}_{ij} &:= \mathsf{f}_{ij}(\alpha(t)) &&= d_{ij}\!\left(\alpha_{j}(t)-\alpha_i(t)\right), 
&\qquad \forall j\neq i, \\
\widehat{\mathsf{f}}_{ij} &:= \widehat{\mathsf{f}}_{ij}(\beta(t)) &&= \widehat d_{ij}\!\left(\beta_{i}(t)-\beta_j(t)\right), 
&\qquad \forall j\neq i, \\
\mathsf{g}_{ij} &:= \mathsf{g}_{ij}(\alpha(t)) &&= m_{ij}\!\left(\alpha_{i}(t)-\alpha_j(t)\right), 
&\qquad \forall j\neq i, \\
\widehat{\mathsf{g}}_{ij} &:= \widehat{\mathsf{g}}_{ij}(\beta(t)) &&= m_{ij}\!\left(\beta_{i}(t)-\beta_j(t)\right), 
&\qquad \forall j\neq i.
\end{alignat}
The correction terms are defined as
\begin{equation}\label{correction_term}
\overline{\mathsf{f}}_{i}=\sum_{j\neq i}\mathfrak{a}_{ij}\mathsf{f}_{ij},\;\;\overline{\widehat{\mathsf{f}}}_{i}=\sum_{j\neq i}\widehat{\mathfrak{a}}_{ij}\widehat{\mathsf{f}}_{ij},\;\;\overline{\mathsf{g}}_{i}=\sum_{j\neq i}\mathrm{a}_{ij}\mathsf{g}_{ij}\;\;\text{and}\;\;\overline{\widehat{\mathsf{g}}}_{i}=\sum_{j\neq i}\widehat{\mathrm{a}}_{ij}\widehat{\mathsf{g}}_{ij},\;\;\;\;i=1,\ldots,\N,
\end{equation}
where the correction factors $\mathfrak{a}_{ij} = \mathfrak{a}_{ij}(\al),\,\widehat{\mathfrak{a}}_{ij} = \widehat{\mathfrak{a}}_{ij}(\be)$ and $\mathrm{a}_{ij} = \mathrm{a}_{ij}(\al^\prime,\al),\,\widehat{\mathrm{a}}_{ij} = \widehat{\mathrm{a}}_{ij}(\be^\prime,\be).$ In addition, $\mathfrak{a}_{ij}=\mathfrak{a}_{ji},\,\widehat{\mathfrak{a}}_{ij}=\widehat{\mathfrak{a}}_{ji},\,\mathrm{a}_{ij}=\mathrm{a}_{ji},\,\widehat{\mathrm{a}}_{ij}=\widehat{\mathrm{a}}_{ji}\in[0,1],\;i,j=1,\ldots,\N,$ are appropriately defined.\par
For the rest of this paper we will call the internodal fluxes as anti-diffusive fluxes. Every anti-diffusive flux $\mathsf{f}_{ij},\,\widehat{\mathsf{f}}_{ij},\,\mathsf{g}_{ij},\,\widehat{\mathsf{g}}_{ij}$ is multiplied by a solution-depended correction factor $\mathfrak{a}_{ij},\,\widehat{\mathfrak{a}}_{ij},\,\mathrm{a}_{ij},\,\widehat{\mathrm{a}}_{ij}\in[0,1]$, respectively, before it is inserted into the equation. Hence, the AFC scheme is the following: We seek $\al(t),\,\be(t)\in\Rn$ such that, for $t\in[0,T]$,
\begin{align}
\M_L\al^\prime(t) + (\mu\, \S + \Q + \D + \sigma\,\M_L )  \al(t) & = \bfr(Bu_h) + \bfr(G) + \overline{\mathsf{f}}(\al(t)) + \overline{\mathsf{g}}(\al^\prime(t),\al(t)),\label{AFC_matrix_y}\\
-\M_L \be^\prime(t) + (\mu\,\S - \Q  - \widehat\D + \sigma\,\M_L) \be(t) & = \M\,\al(t) - \bfr(y_d) + \overline{\widehat{\mathsf{f}}}(\be(t)) - \overline{\widehat{\mathsf{g}}}(\be^\prime(t),\be(t)),\label{AFC_matrix_p}
\end{align}
with $\al(0)=\al^0,\,\be(T)=\mathbf{0}.$ \par
To ensure that the AFC scheme satisfies the discrete maximum principle, we choose the correction factors 
$\mathsf{a}_{ij}$ such that the sum of anti-diffusive fluxes is constrained by, see, cf. e.g., \cite{kuzmin2010book},
\begin{equation}\label{led_2D_weakened}
q_i\left(\delta_{i}^{\min}(t)-\delta_{i}(t)\right)\leq \sum_{j\neq i}\mathsf{a}_{ij}\mathsf{p}_{ij}\leq q_i\left(\delta_{i}^{\max}(t)-\delta_{i}(t)\right),\;\;q_i\geq 0,
\end{equation}
where either $\mathsf{a}_{ij} = \mathfrak{a}_{ij}$ and $\mathsf{p}_{ij} = \mathsf{f}_{ij}$ with $\delta(t) = \alpha(t),$ or $\mathsf{a}_{ij} = \widehat{\mathfrak{a}}_{ij}$ and $\mathsf{p}_{ij} = \widehat{\mathsf{f}}_{ij}$ with $\delta(t) = \beta(t),$ for the corrections factors associated with the matrices $\D,\,\widehat\D,$ respectively. In addition, either $\mathsf{a}_{ij} = \mathrm{a}_{ij}$ and $\mathsf{p}_{ij} = \mathsf{g}_{ij}$ with $\delta(t) = \alpha(t)$ or $\mathsf{a}_{ij} = \widehat{\mathrm{a}}_{ij}$ and $\mathsf{p}_{ij} = \widehat{\mathsf{g}}_{ij}$ with $\delta(t) = \beta(t)$ for the correction factors associated with the mass lumping. Here, $\delta^{\max}_i,\,\delta^{\min}_i,$ denotes the local maximum and local minimum at $\omega_i,$ respectively. \par

To formulate the variational form, we define for a function $s\in \Sh,$ its nodal values as $s_i = s(Z_i),\;i=1,\ldots,\N,$ the bilinear form $\widehat d_{\D,h}(w;\cdot,\cdot):\CC \times \CC \to \R,$ with $w\in\CC,$ is defined by
\begin{equation}\label{stab_term_D_widehat}
\begin{aligned}
\widehat d_{\D,h}(w;v,z) & := \sum_{i,j=1}^{\N}\,d_{ij}(1-\mathfrak{a}_{ij}(w))(v_i - v_j)z_i\\
& = \sum_{i<j}\,d_{ij}(1-\mathfrak{a}_{ij}(w))(v_i - v_j)(z_i - z_j),\quad\forall v,\,z\in{\CC},
\end{aligned}
\end{equation}
and the in a similar way, the bilinear form $\widehat d_{\widehat \D,h}(w;\cdot,\cdot):\CC \times \CC \to \R,$ with $w\in\CC,$ is defined by
\begin{equation}\label{stab_term_D_widehat2}
\begin{aligned}
\widehat d_{\widehat \D,h}(w;v,z) & := \sum_{i,j=1}^{\N}\,\widehat d_{ij}(1-\widehat{\mathfrak{a}}_{ij}(w))(v_i - v_j)z_i\\
& = \sum_{i<j}\,\widehat d_{ij}(1 - \widehat{\mathfrak{a}}_{ij}(w))(v_i - v_j)(z_i - z_j),\quad\forall v,\,z\in{\CC}.
\end{aligned}
\end{equation}
In both of the stabilization terms, the last equality is due to the symmetry of matrix $\D$ and $\widehat \D,$ respectively, see, e.g., \cite{barrenechea2018}.

\begin{remark}\label{remark:estimate_optimal}
Since the bounds \eqref{est_dij_2D_L_infty} and \eqref{est_dij_2D_L_infty_hat} hold in both the convection-dominated and diffusion-dominated cases, we can, by Lemma \ref{lemma:stability_widehat_D_main} and Remark \ref{remark:optimal}, derive an optimal estimate for the stabilization terms \eqref{stab_term_D_widehat} and \eqref{stab_term_D_widehat2}. 
\end{remark}

Further, the stabilization term due to the mass lumping, is the bilinear form $\widehat d_{\M_L,h}(w;\cdot,\cdot):\CC \times \CC \to \R,$ with $w,\,s\in\CC,$ which defined by
\begin{equation}\label{stab_term_M_widehat}
\begin{aligned}
\widehat d_{\M_L,h}(w,s;v,z) & := \sum_{i,j=1}^{\N}\,m_{ij}(1-\mathit{a}_{ij}(w,s))(v_i - v_j)z_i \\
& = \sum_{i<j}\,m_{ij}(1-\mathit{a}_{ij}(w,s))(v_i - v_j)(z_i - z_j),\quad\forall v,\,z\in{\CC},
\end{aligned}
\end{equation}
where $\mathit{a}_{ij} = \mathrm{a}_{ij}$ or $\mathit{a}_{ij} = \widehat{\mathrm{a}}_{ij}.$

The variational formulation of the algebraic flux correction scheme \eqref{AFC_matrix_y}--\eqref{AFC_matrix_p} is as follows. We seek $(y_h,p_h,u_h)\in \Sh \times \Sh \times \Uad$, such that
\begin{align}
(y_{h,t}, \chi) + \alpha_h(y_h, \chi)  - \widehat d_{\D,h}(y_h;y_h,\chi) + \widehat d_{\M_L,h}(y_{h,t},y_h;y_{h,t},\chi) & = (Bu_h + G, \chi),\label{optimal_AFC_semi_weak_y}
\end{align}
for all $\chi \in\Sh$ with $y_h(0)=y_h^0,$ and
\begin{align}
- (p_{h,t}, \chi) + \alpha_h(\chi, p_h) + \widehat d_{\widehat\D,h}(p_h;p_h,\chi) - \widehat d_{\M_L,h}(p_{h,t},p_h;p_{h,t},\chi) & = (y_h - y_d,\chi),\label{optimal_AFC_semi_weak_p}
\end{align}
for all $\chi \in\Sh$ with $p_h(T)=0.$ The control variable is defined as
\begin{align}
(\lambda\, u_h +  B^*p_h, \psi - u_h) & \geq 0, \;\;\;\psi\in\Uad.\label{optimal_AFC_semi_weak_u}
\end{align}
Similarly to \eqref{optimal_weak_u2}, one can express \eqref{optimal_AFC_semi_weak_u} by using \eqref{operator_Pi}, as
\begin{align}\label{optimal_AFC_semi_weak_u2}
u_h(t) := \Pi_{[\mathsf{u_a},\mathsf{u_b}]}(-\lambda^{-1} B^*p_h(t)),\;\;\;t\in[0,T].
\end{align}

The semi--discrete scheme \eqref{optimal_AFC_semi_weak_y}--\eqref{optimal_AFC_semi_weak_u2} is a coupled non-linear scheme. To prove the existence and the uniqueness of the solution, we will use a fixed point argument. To be able to form this argument, we first need to specify the correction factors of the stabilization terms.

\begin{remark}\label{remark:equivalence_of_schemes}
Using the identity {\normalfont{\cite[Equation 15.7]{thomee2006}}} as well as the symmetry of $\M$ (cf. {\normalfont{\cite[Section 2.1]{barrenechea2018}}}), one can show that
\begin{align}\label{mass_equival}
(\psi, \chi) = (\psi, \chi)_h - \sum_{i,j=1}^{\N}m_{ij}(\psi_i - \psi_j)\chi_i  = (\psi, \chi)_h - \sum_{i<j}m_{ij}(\psi_i - \psi_j)(\chi_i - \chi_j),
\end{align}
and since the last term on the right-hand side is included in the definition of $\widehat d_{\M_L,h},$ 
in \eqref{optimal_AFC_semi_weak_y} and \eqref{optimal_AFC_semi_weak_p}, 
the first term on the left-hand side is $(\cdot,\cdot)$ rather than $(\cdot,\cdot)_h$.
\end{remark}

\subsection{Correction factors}\label{subsection:corr_factors}

To compute the coefficients $\mathsf{a}_{ij}$, we need to fix a set of nonnegative $q_{i}$, $i=1,\dots,\N$ in \eqref{led_2D_weakened}. Their choice can be arbitrary, but efficiency and accuracy can dictate a strategy, see, e.g,  \cite{kuzmin2010book}. For a detailed presentation we refer to \cite[Chapter 4]{kuzmin2010book}, \cite[Section 6]{barrenechea2024}, \cite[Chapter 10]{barrenechea2025}. Examples of correction factors can be found also in \cite{barrenechea2016,barrenechea2018,barrenechea2024,kuzmin2010book,kuzmin2002,kuzmin2004} and the references therein.
We shall assume that the correction factors $\mathsf{a}_{ij} = \mathfrak{a}_{ij}$ or $\mathsf{a}_{ij} = \widehat{\mathfrak{a}}_{ij}$ satisfies the following two properties.

\begin{assumption}(Linearity Preservation)\label{assumption:linearity_preservation}
Let $Z_i$ an inner node of $\Th.$ The limiters $\mathsf{a}_{ij} = \mathfrak{a}_{ij}$ or $\mathsf{a}_{ij} = \widehat{\mathfrak{a}}_{ij}$ are linearity preserving, i.e., for every edge $e\in\mathcal{E}_h^0$ with endpoints $Z_i,\,Z_j,$
\begin{align*}
\mathsf{a}_{ij}(v) = 1\;\;\text{ if }\;\;v\in \mathbb{P}_1(\mathbb{R}^2).
\end{align*}
\end{assumption}

\begin{assumption}(Local Estimate)\label{assumption:local_estimate_factors}
Let $Z_i$ an inner node of $\Th.$ We assume that for every edge $e\in\mathcal{E}_h^0,$ with endpoints $Z_i,\,Z_j$, the following estimate holds for the limiters $\mathsf{a}_{ij} = \mathfrak{a}_{ij}$ or $\mathsf{a}_{ij} = \widehat{\mathfrak{a}}_{ij},$
\begin{align}\label{local_estimate}
\vert d_{ij} \mathsf{a}_{ij}(\chi)(\chi_i - \chi_j) & - d_{ij} \mathsf{a}_{ij}(\widetilde\chi)(\widetilde\chi_i - \widetilde\chi_j) \vert \leq C\,h\sum_{\ell\in \Zh(\omega_i)}|\chi_\ell - \widetilde\chi_\ell|,
\end{align}
for all $\chi,\,\widetilde\chi\in \Sh.$ The constant $C,$ may depend on the norms of the unknown solution or the data of the problem, but is independent of $h.$
\end{assumption}

We compute the correction factors using Algorithm  \ref{algorithm-1}, which has been proposed  by \cite{barrenechea2018}.  

\begin{algorithm}[H]
Given data: 
\begin{enumerate}
\item The positive  coefficients $q_{i}$, $i,j=1,\dots,\N,$
\item The fluxes $\mathsf{p}_{ij}$, $i\neq j$, $i,j=1,\dots,\N$, i.e., $\mathsf{p}_{ij} = \mathsf{f}_{ij},\,\widehat{ \mathsf{f}}_{ij}$ or $\mathsf{g}_{ij} = \mathsf{g}_{ij},\,\widehat{\mathsf{g}}_{ij}$
\item The coefficients $\delta_j,$ for  $j=1,\dots,\N$.
\end{enumerate}
\noindent
The computation of factors $\mathsf{a}_{ij},$ for $i,\,j\in\Nh,$ that satisfies \eqref{led_2D_weakened} are defined as follows.
\begin{enumerate}
\item
Compute for $i\in\Nh^0,\,j\in\Nh,$ the limited sums $P_{i}^{\pm}:= P_{i}^{\pm}(\del)$ for $\mathsf{p}_{ij} = \mathsf{f}_{ij},\,\widehat{ \mathsf{f}}_{ij}$ or $P_{i}^{\pm}:= P_{i}^{\pm}(\del^\prime)$ for $\mathsf{g}_{ij} = \mathsf{g}_{ij},\,\widehat{\mathsf{g}}_{ij}$, of positive and negative anti-diffusive fluxes
\begin{align*}
P_{i}^{+} = \sum_{j\in \Zh^i}\max\{0, \mathsf{p}_{ij}\},
\quad \text{ and }\quad P_{i}^{-} = \sum_{j\in \Zh^i}\min\{0, \mathsf{p}_{ij}\}.
\end{align*}
\item
Retrieve for $i\in\Nh^0,\,j\in\Nh,$ the local extremum diminishing upper and lower bounds  ${Q}_{i}^{\pm} : = {Q}_{i}^{\pm}(\del),$
\begin{align*}
{Q}_{i}^{+} = q_i(\delta_i^{\max} - \delta_i),
\quad\text{ and }\quad{Q}_{i}^{-} = q_i(\delta_i^{\min} - \delta_i),
\end{align*}
where $\delta_i^{\max},\;\delta_i^{\min}$ are the local maximum and local minimum at $\omega_i.$
\item
Compute for $i\in\Nh^0,\,j\in\Nh,$ also the coefficients $\overline{\mathsf{a}}_{ij},$ for $j\neq i$ are given by
\begin{equation}\label{correction_factors_definition}
\begin{aligned}
R_{i}^+=\min\left\lbrace 1, \frac{Q_{i}^+}{P_{i}^+} \right\rbrace,\quad R_{i}^-= \min\left\lbrace 1, \frac{Q_{i}^-}{P_{i}^-}\right\rbrace\quad \text{and}\quad \overline{\mathsf{a}}_{ij} = 
\begin{cases}
R_{i}^+ , &\text{if} \quad\mathsf{p}_{ij} > 0,\\
1, &\text{if} \quad\mathsf{p}_{ij}= 0,\\
R_{i}^-, &\text{if} \quad\mathsf{p}_{ij} < 0.
\end{cases}
\end{aligned}
\end{equation}
For $P_i^{\pm} = 0,$ we set $R_i^{\pm} = 1.$
\end{enumerate}  
Then, the coefficients $\mathsf{a}_{ij},$ for $j\neq i$ with $i\in\Nh^0,\,j\in\Nh,$ are given by $\mathsf{a}_{ij} = \min\{\overline{\mathsf{a}}_{ij}, \overline{\mathsf{a}}_{ji}\}$ and $\mathsf{a}_{ji} = \mathsf{a}_{ij}.$ For the Dirichlet nodes, i.e., $i\in\Nh^b,\,j\in\Nh^b,$ we set $\mathsf{a}_{ij} = 1.$
\caption{{\normalfont\cite[Lemma 6]{barrenechea2018}} Computation of the correction factors}\label{algorithm-1}
\end{algorithm}

\begin{remark}\label{remark:linearity_preservation1}
In view of {\normalfont{\cite[Section 4.2]{barrenechea2018}}}, the limiter \eqref{correction_factors_definition} is linearity preserving if $q_i,\,i\in\Nh^0,$ are computed by
\begin{align}\label{def:q_i-1}
q_i : = \gamma_i \sum_{j\in \Zh^i} c_{ij},\quad i\in\Nh^0,\;\;\gamma_i=\dfrac{\max_{Z_j\in \partial\omega_i}|Z_i-Z_j|}{\mathrm{dist}(Z_i,\partial \omega_i^{\mathrm{conv}})},
\end{align}
with $\omega_i^{\mathrm{conv}}$  the convex hull of $\omega_i$, for $i\in \Nh^0.$ Note that if $\omega_i$ is symmetric with respect to $Z_i$ then $\gamma_i = 1$, see \cite{barrenechea2017b}.\par
In addition, for $c_{ij} : = d_{ij},$ the limiter \eqref{correction_factors_definition} satisfies the local estimate \eqref{local_estimate}, see, e.g., {\normalfont{\cite[Lemma 2.15]{chatzipantelidis2022}}} for $q_i = \mathcal{O}(h),$ which is true, see, e.g., \eqref{def:q_i-1}. A similar estimate can be also found in {\normalfont{\cite[Lemma 10.37]{barrenechea2025}}}. Note that $d_{ij} = \mathcal{O}(h).$ These local estimates are the basis of the derivation of an optimal estimate for the stabilization term, see, e.g., {\normalfont{\cite[Lemma 2.17]{chatzipantelidis2022}}}, also Lemmas \ref{lemma:stability_widehat_D_main}, \ref{lemma:stability_widehat_D_difference}, \ref{lemma:stability_widehat_M_difference}. 
\end{remark}

\begin{remark}\label{remark:gamma}
In {\normalfont{\cite[Lemma 6.1]{barrenechea2018}}}, the authors proved that if $\omega_i$ is symmetric with respect to $Z_i,$ then $\gamma_i=1.$ An example of such a partition depicted in Fig. \ref{fig:triangulation}.
\end{remark}

Similar to \cite[Lemma 7]{barrenechea2016}, the correction factors defined above can be written in the following generalization formula.

\begin{definition}\label{definition:corr_factors_gen}
Let $\bpsi,\,\zet\in\R^{\N},$ the coefficient vectors that correspond to the finite element functions $\psi,\,\zeta,$ respectively. Then, 
\begin{equation}\label{correction_factors_formulation_old}
\mathtt{a}_{ij}(\bpsi,\zet) = \frac{A_{ij}(\bpsi,\zet)}{\vert\psi_j - \psi_i\vert + B_{ij}(\bpsi)}, \quad i,j\in\Nh^0
\end{equation}
with $\psi_j\neq \psi_i$, and $A_{ij}$ and $B_{ij}$ non-negative functions which are continuous functions in $\bpsi,$ with $\psi_i\neq\psi_j,$ defined by for $i,j\in\Nh^0$ and $c_{ij} > 0,$ as
\begin{align}\label{def_A_ij_first}
A_{ij}(\bpsi,\zet) := \frac{1}{c_{ij}}\begin{cases}
\min\{ -P_{i}^-(\bpsi), -Q_{i}^-(\zet)\}, & \quad\text{ if }\ \psi_i < \psi_j,\\
\min\{ P_{i}^+(\bpsi), Q_{i}^+(\zet)\}, & \quad\text{ if }\ \psi_i > \psi_j,
\end{cases}
\end{align}
and 
\begin{align}\label{def_B_ij_second}
B_{ij}(\bpsi) := \frac{1}{c_{ij}}\begin{cases}
-\widetilde {P}_{i}^-(\bpsi), & \quad\text{ if }\ \psi_i < \psi_j,\\
\widetilde {P}_{i}^+(\bpsi), & \quad\text{ if }\ \psi_i > \psi_j,
\end{cases}
\end{align}
where $P_{i}^{\pm}(\cdot,\cdot),\,Q_{i}^{\pm}(\cdot,\cdot)$ as defined as
\begin{align*}
P_{i}^{+}(\bpsi) = \sum_{j\in\Zh^i}\max\{0, c_{ij}(\psi_i - \psi_j)\}\quad\text{and}\quad P_{i}^{-}(\bpsi) = \sum_{j\in\Zh^i}\min\{0, c_{ij}(\psi_i - \psi_j)\},
\end{align*}
\begin{align*}
Q_{i}^{+}(\zet) = (\zeta_i^{\max} - \zeta_i)\gamma_i\sum_{j\in\Zh^i}c_{ij}\quad\text{and}\quad Q_{i}^{+}(\zet) = (\zeta_i^{\min} - \zeta_i)\gamma_i\sum_{j\in\Zh^i}c_{ij},
\end{align*}
where $\zeta_i^{\max},\,\zeta_i^{\min},$ are the local extrema at $\omega_i.$ In addition,
\begin{align*}
\widetilde {P}_{i}^+(\bpsi) & = \sum_{\substack{k=1, k\neq j}}^{\N}\max\{0, c_{ik}(\psi_i - \psi_k)\}\quad\text{and}\quad\widetilde {P}_{i}^-(\bpsi) = \sum_{\substack{k=1, k\neq j}}^{\N}\min\{0, c_{ik}(\psi_i - \psi_k)\}.
\end{align*}
For the case where $c_{ij} = 0,$ we set $\mathtt{a}_{ij}=0.$
\end{definition}

\begin{remark}\label{remark:def_remark}
In view of Definition \ref{definition:corr_factors_gen}, the correction factors defined in \eqref{correction_term}--\eqref{AFC_matrix_p}, can be written as $\mathfrak{a}_{ij}(\al) = \mathtt{a}_{ij}(\al,\al),$ with $c_{ij} = d_{ij}$, the $\widehat{\mathfrak{a}}_{ij}(\be) = \mathtt{a}_{ij}(\be,\be),$  with $c_{ij} = \widehat d_{ij}$, and $\mathrm{a}_{ij}(\al^\prime,\al) = \mathtt{a}_{ij}(\al^\prime,\al),\,\widehat{\mathrm{a}}_{ij}(\be^\prime,\be) = \mathtt{a}_{ij}(\be^\prime,\be),$  with $c_{ij} = m_{ij}$.
\end{remark}

\begin{remark} 
Note that for a finite element function $\psi\in\Sh$ with associated coefficient vector $\bpsi\in\R^{\N},$ we do not
distinguish between them when writing the correction factors or the non-negative functions of Definition \ref{definition:corr_factors_gen}. In particular, we use the notation $\mathfrak{a}_{ij}(\psi) = \mathfrak{a}_{ij}(\bpsi)$ and suppress the dependence of the correction factors from the same finite element function.
\end{remark}

\subsection{Auxiliary results}\label{subsec:auxiliary}

To estimate the error we will split the error by intersecting a suitable projection. We recall the Ritz projection as defined e.g.,  \cite[Chapter 1]{thomee2006} and the references therein. The projection $R_h : H^{1}\to \Sh,$ for a function $v\in H^{1},$ is defined as
\begin{align}
(\nabla (R_hv -  v), \nabla \chi) &  = 0,\quad\forall \chi\in \Sh.\label{ritz_projection_2D}
\end{align} 
In view of the mesh Assumption \ref{mesh-assumption}, $R_h$ satisfy the following bounds, cf. e.g., \cite[Chapter 8]{brenner2008} and \cite{thomee2006}. 
\begin{align}
\|v - R_hv\|_{L^{2}} + h\|v - R_hv\|_{1} &\le Ch^2\|v\|_{2},\,\quad\forall v\in H^{2},\label{ritz_projection_est2_2D}\\
\|\nabla R_h v\|_{L^{2}} &\le \|\nabla v\|_{L^{2}},\quad\;\;\forall\,v\in H^{1}.\label{ritz_projection_stab}
\end{align}

In what follows, we recall various useful results for the stabilization terms of the AFC scheme, that are crucial to our analysis. Recall the bilinear form $d_{\widehat\D,h}$, introduced in \eqref{stab_term}, and hence also $\widehat{d}_{\widehat\D,h}$, defined in \eqref{stab_term_D_widehat2}, induces a seminorm on $\CC$, thanks to \cite[Lemma 1]{barrenechea2016}. This is also true for $-d_{\D,h},\,-\widehat d_{\D,h},$ since the off-diagonal elements of $\D,$ are non-positive, see their definition in \eqref{D_def}.\par

Therefore, $ d_{\widehat\D,h}(\cdot,\cdot),\,\widehat d_{\widehat\D,h}(w;\cdot,\cdot):{\CC}\times {\CC}\to{\R},$ with $w\in\Sh$, are a non-negative symmetric bilinear forms 
which satisfies the Cauchy-Schwartz's inequality,
\begin{equation} \label{Schwartz_ineq_afc_2D}
|\widehat d_{\widehat \D,h}(w;v,z)|^2\leq \widehat d_{\widehat \D,h}(w;v,v)\,\widehat d_{\widehat \D,h}(w;z,z),\quad\forall v,\,z\in{\CC},
\end{equation}  
and thus induces a seminorm on ${\CC}$. For the stabilization term \eqref{stab_term}, we have the following two results, see \cite{barrenechea2016} since $d_{ij} = \mathcal{O}(h),\,\widehat d_{ij} = \mathcal{O}(h),$ see \eqref{est_dij_2D_L_infty}.

\begin{lemma}\normalfont{\cite[Lemma 16]{barrenechea2016}}\label{lemma:estimate_stab_low}
There are exists a positive constant $C,$ independent of $h,$ such that,
\begin{align*}
\vert\widehat d_{\D,h}(\psi;\psi,\chi)\vert + \vert\widehat d_{\widehat\D,h}(\psi;\psi,\chi)\vert + |d_{\D,h}(\psi, \chi)| + |d_{\widehat \D,h}(\psi, \chi)| \leq Ch\|\bfb\|_{L^{\infty}}\|\nabla \psi\|_{L^{2}}\|\nabla \chi\|_{L^{2}},
\end{align*}
where the constant $C,$ depend on the shape regularity constant $\varpi,$ see \eqref{shape_regularity}.
\end{lemma}

For the bilinear forms  \eqref{stab_term_D_widehat}, \eqref{stab_term_D_widehat2} a more sharp estimate holds, see \cite[Lemma 2.17]{chatzipantelidis2022}. While in \cite{chatzipantelidis2022} the proof is based on a slightly modified version of Algorithm \ref{algorithm-1} suitable for Neumann boundary conditions, the proof can be easily extended in the case of Dirichlet boundary conditions.

\begin{lemma}{\normalfont{\cite[Lemma 2.17]{chatzipantelidis2022}}}\label{lemma:stability_widehat_D_main}
For $v\in H^{2}$ and the correction factors satisfies Assumptions \ref{assumption:linearity_preservation}, \ref{assumption:local_estimate_factors} and in view of Remarks \ref{remark:linearity_preservation1}, \ref{remark:gamma}, there are exists a positive constant $C,$  independent of $h,$ such that for all $\psi,\,\chi\in \Sh,$
\begin{align}
\vert\widehat d_{\D,h}(\psi;\psi,\chi)\vert + \vert\widehat d_{\widehat\D,h}(\psi;\psi,\chi)\vert & \leq Ch(\|\nabla (\psi - v)\|^2_{L^{2}} + h^2\|v\|^2_{2})^{1/2}\|\nabla \chi\|_{L^{2}}.\label{stability_stab_term_widehat_D}
\end{align}
\end{lemma}
\begin{proof}
The proof can be found in \cite[Lemma 2.17]{chatzipantelidis2022} for a slightly modified version of Algorithm \ref{algorithm-1} that applies when the boundary of $\Omega$ has a part with Neumann boundary conditions. For completeness reasons, we describe the briefly the idea of the proof.

Using a similar argument to \cite[Lemma 3]{barrenechea2018}, \cite[Lemma 2.17]{chatzipantelidis2022}, we can define a local problem on $\omega_i,$ with polynomial solution, i.e., $\zeta_i\in\mathbb{P}_1(\omega_i)$ and $\zeta_i = 0,$ outside $\omega_i.$ In particular, recall the auxiliary problem defined in \cite[Lemma 2.17]{chatzipantelidis2022}. For a $v\in H^{2},$ let $\zeta^i\in \mathbb{P}_1(\omega_{i})$ be the unique solution of 
\begin{align*}
(\nabla \zeta^{i}, \nabla \chi)_{L^{2}(\omega_{i})} & = (\nabla v, \nabla \chi)_{L^{2}(\omega_{i})},\;\;\;\forall\,\chi\in \mathbb{P}_1(\omega_{i}),\\
(\zeta^{i},1)_{L^{2}(\omega_{i})} & = (v,1)_{L^{2}(\omega_{i})}.
\end{align*}
We can extend $\zeta^{i}$ arbitrarily as a function on $\Sh$. Also, using standard estimates, there exists a constant $C$ independent of $h$, such that
\begin{align}\label{est_1}
\|\nabla (v - \zeta^{i})\|_{L^{2}(\omega_{i})} \leq Ch\|v\|_{H^{2}(\omega_{i})}.
\end{align}

Then, since the limiter defined in \eqref{correction_factors_definition} satisfies the linearity preservation property, see, e.g., Assumption \ref{assumption:linearity_preservation}, we immediately obtain that $\widetilde\rho_{{ij}}(\zeta_{i}) = 0,$ where $\widetilde\rho_{ij}(\psi) : = (1 - \mathfrak{a}(\psi))(\psi_i - \psi_j),\,\psi\in\Sh.$ Hence, one get
\begin{align*}
\widehat{d}_{\D,h}(\psi;\psi,\chi) & = \sum_{i<j}d_{ij}( \widetilde{\rho}_{{ij}}(\psi) - \widetilde\rho_{{ij}}(\zeta_{i}))(\chi_i - \chi_j),
\end{align*}
and then the desired result obtained by \eqref{local_estimate}, the Cauchy-Schwartz inequality, the above estimate and the shape regularity of the triangulation. In particular, using the local estimate \eqref{local_estimate}, we obtain,
\begin{equation}\label{est_2}
\begin{aligned}
\sum_{i<j}d_{ij}( \widetilde{\rho}_{{ij}}(\psi) - \widetilde\rho_{{ij}}(\zeta_{i}))(\chi_i - \chi_j) & \leq \left(\sum_{i<j}d_{ij}^2( \widetilde{\rho}_{{ij}}(\psi) - \widetilde\rho_{{ij}}(\zeta_{i}))^2\right)^{1/2}\left(\sum_{i<j}(\chi_i - \chi_j)^2\right)^{1/2}\\
& \leq \left(\sum_{i<j}d_{ij}^2( \widetilde{\rho}_{{ij}}(\psi) - \widetilde\rho_{{ij}}(\zeta_{i}))^2\right)^{1/2}\|\nabla \chi\|_{L^{2}},
\end{aligned}
\end{equation}
where, using the local estimate \eqref{local_estimate} and \eqref{def:q_i-1}, we obtain,
\begin{align*}
d_{ij}( \widetilde{\rho}_{{ij}}(\psi) - \widetilde\rho_{{ij}}(\zeta_{i})) & \leq C\,h\sum_{\ell\in \Zh(\omega_i)}|\psi_\ell - (\zeta_i)_\ell|\\
& \leq Ch\|\nabla (\psi - \zeta_i)\|_{L^{2}(\omega_i)} \\
& \leq Ch\left(\|\nabla (\psi - v)\|_{L^{2}(\omega_i)} + \|\nabla (v - \zeta_i)\|_{L^{2}(\omega_i)}\right) \leq Ch\left(\|\nabla (\psi - v)\|_{L^{2}(\omega_i)} + h\|v\|_{H^{2}(\omega_i)}\right).
\end{align*}
Hence, the result follows by combining the latter results with \eqref{est_2}, together with the shape regularity of the triangulation.
\end{proof}

\begin{remark}\label{remark:optimal}
To obtain an optimal estimate for the stabilization term in Lemma \ref{lemma:stability_widehat_D_main}, we need to ensure that the term $\|\nabla (\psi - v)\|_{L^{2}}$ is of order $\mathcal{O}(h),$ for some $\psi\in\Sh$ and $v\in H^{2}.$ This can be done if, e.g., if $\psi = Y^n$ and $v := \overline y^n,$  see also how the term $I_2$ on the proof of Theorem \ref{theorem:error_estimates_AFC} is estimated.
\end{remark}

\begin{remark}\label{remark:optimal2}
Note that the third part of Assumption~\ref{mesh-assumption} is crucial. If it were not true, the matrices $\D$ and $\widehat{\D}$ would have to be chosen so that $\mu \S + \Q + \D$ and $\mu \S - \Q - \widehat{\D}$ have non-positive off-diagonal entries and positive diagonal entries. This would imply that $|d_{ij}| \leq 2\mu |s_{ij}| + |\tau_{ij}| + |\tau_{ij}|$, hence $|d_{ij}| \leq C\mu + C(\varpi)\|\bfb\|_{L^{\infty}}h$. Similarly, $|\widehat d{ij}| \leq C\mu + C(\varpi)\|\bfb\|_{L^{\infty}}h$. In the diffusion-dominant regime, i.e., when $C\mu \gg C(\varpi)\|\bfb\|_{L^{\infty}}h$, the estimates derived in Lemma \ref{lemma:stability_widehat_D_main} would then behave (cf. Remark~\ref{remark:optimal}) as $\mathcal{O}(h)$.
\end{remark}

In view of the stabilization terms defined in \eqref{stab_term} and \eqref{stab_term_D_widehat}, we define the bilinear form $\widetilde d_{\D,h}(w;\cdot,\cdot):\CC \times \CC \to \R,$ with $w\in\Sh,$ is defined by
\begin{equation}\label{stab_term_D_widetilde}
\widetilde d_{\D,h}(w;v,z) := d_{\D,h}(v,z) - \widehat d_{\D,h}(w;v,z) = \sum_{i<j}\,d_{ij}\mathfrak{a}_{ij}(w)(v_i - v_j)(z_i - z_j),\quad\forall v,\,z\in{\CC},
\end{equation}
and the in a similar way, the bilinear form $\widetilde d_{\widehat \D,h}(w;\cdot,\cdot):\CC \times \CC \to \R,$ with $w\in\Sh,$ is defined by
\begin{equation}\label{stab_term_D_widetilde2}
\widetilde d_{\widehat \D,h}(w;v,z) := d_{\widehat \D,h}(v,z) - \widehat d_{\widehat\D,h}(w;v,z) = \sum_{i<j}\,\widehat d_{ij}\widehat{\mathfrak{a}}_{ij}(w)(v_i - v_j)(z_i - z_j),\quad\forall v,\,z\in{\CC}.
\end{equation}

For the stabilization terms defined \eqref{stab_term_D_widehat}, \eqref{stab_term_D_widehat2} as well as those defined in \eqref{stab_term_D_widetilde}, \eqref{stab_term_D_widetilde2}, the following estimates are hold.

\begin{lemma}{\normalfont{\cite[Lemma 2.17]{chatzipantelidis2022}}}\label{lemma:stability_widehat_D_difference}
Let the bilinear forms $\widehat d_{\D,h},\,\widehat d_{\widehat\D,h},\,\widetilde d_{\D,h},\, \widetilde d_{\widehat\D,h}$ defined in \eqref{stab_term_D_widehat}, \eqref{stab_term_D_widehat2}, \eqref{stab_term_D_widetilde}, \eqref{stab_term_D_widetilde2}, respectively. If their correction factors are satisfying Assumption \ref{assumption:local_estimate_factors} and in view of Remarks \ref{remark:linearity_preservation1}, \ref{remark:gamma}, there exists for $w,\,s,\,v,\,\chi\in \Sh,$ a positive constant $C$, independent of $h,$ such that,
\begin{align}
\vert\widehat{d}_{\D,h}(v;v, \chi) - \widehat{d}_{\D,h}(w;w,\chi)\vert + \vert\widetilde{d}_{\D,h}(v;v, \chi) - \widetilde{d}_{\D,h}(w;w,\chi)\vert & \leq Ch\|\nabla (v - w)\|_{L^{2}}\|\nabla \chi\|_{L^{2}},\label{stability_stab_term_widetilde_D}\\
\vert\widehat{d}_{\widehat\D,h}(v;v, \chi) - \widehat{d}_{\widehat\D,h}(w;w,\chi)\vert  + \vert\widetilde{d}_{\widehat\D,h}(v;v, \chi) - \widetilde{d}_{\widehat\D,h}(w;w,\chi)\vert & \leq Ch\|\nabla (v - w)\|_{L^{2}}\|\nabla \chi\|_{L^{2}}.\label{stability_stab_term_widetilde_D2}
\end{align}
\end{lemma}
\begin{proof}
The proof can be found in \cite[Lemma 2.18]{chatzipantelidis2022}. For completeness reasons, we provide a short proof.\par
Using the definition $\widehat{d}_{\D,h},$ see, e.g., \eqref{stab_term_D_widehat} and setting $\rho_{ij}(v) := d_{ij}(1-\mathfrak{a}_{ij}(v))(v_i - v_j),\,v\in\Sh,$ we get 
\begin{align*}
\widehat{d}_{\D,h}(v;v, \chi) - \widehat{d}_{\D,h}(w;w, \chi) & = \sum_{i<j} (\rho_{ij}(v) - \rho_{ij}(w)) (\chi_i - \chi_j)\\
& \leq \left(\sum_{i<j} (\rho_{ij}(v) - \rho_{ij}(w))^2\right)^{1/2}\left(\sum_{i<j} (\chi_i - \chi_j)^2\right)^{1/2}\\
& \leq Ch\left(\sum_{i<j}\sum_{\ell\in \Zh(\omega_i)}|v_\ell - w_\ell|^2\right)^{1/2} \|\nabla \chi\|_{L^{2}} \leq Ch\|\nabla (v - w)\|_{L^{2}}\|\nabla \chi\|_{L^{2}},
\end{align*}
where we have used the Assumption \ref{assumption:local_estimate_factors}, the fact that $d_{ij} = \mathcal{O}(h)$ and that the triangulation is shape regular. To prove the result for \eqref{stab_term_D_widetilde}, we set $\rho_{ij}(v) := d_{ij}\mathfrak{a}_{ij}(v)(v_i - v_j),\,v\in\Sh,$ and performs similar calculations. The estimate \eqref{stability_stab_term_widetilde_D2}, can be concluded analogously.
\end{proof}

Similarly to the stabilization terms due to matrices $\D,\,\widehat\D,$ we also define the corresponding one for the matrix $\M_L.$ We define the bilinear forms $\widetilde{d}_{\M_L,h}(s;\cdot,\cdot),\, d_{\M_L,h}(\cdot,\cdot)\,:\,{\CC}\times {\CC}\to{\R},$ with $w,\,s\in\CC$, for $v,\,z\in{\CC}$,
\begin{align}
d_{\M_L,h}(v,z) & := \sum_{i<j}m_{ij}(v_i - v_j)(z_i - z_j),\label{stab_term_M_low}\\
\widetilde d_{\M_L,h}(w,s;v,z) & := d_{\M_L,h}(v,z) - \widehat d_{\M_L,h}(w,s;v,z)  = \sum_{i<j}m_{ij}\mathit{a}_{ij}(w,s)(v_i - v_j)(z_i - z_j),\label{stab_term_M_widetilde}
\end{align} 
where $\mathit{a}_{ij} = \mathrm{a}_{ij}$ or $\mathit{a}_{ij} = \widehat{\mathrm{a}}_{ij}.$\par

Following the proof of \cite[Lemma 2.10]{chatzipantelidis2022}, the following estimate for the stabilization terms $d_{\M_L,h},\,\widetilde{d}_{\M_L,h},\,\widehat{d}_{\M_L,h}$ can be proved.
\begin{lemma}\label{lemma:stability_widehat_M_main}
There are exists a positive constant $C,$  independent of $k,\,h,$ such that, for $\zeta,\eta,\psi,\chi\in \Sh$,
\begin{align}
\vert d_{\M_L,h}(\psi,\chi)\vert + \vert \widetilde d_{\M_L,h}(\zeta,\eta;\psi,\chi)\vert + \vert \widehat{d}_{\M_L,h}(\zeta,\eta;\psi,\chi)\vert & \leq Ch^2\|\nabla \psi\|_{L^{2}}\|\nabla \chi\|_{L^{2}}.\label{stability_stab_term_widehat_M}
\end{align}
\end{lemma}
\begin{proof}
Given that for the elements of $\M,$ the $m_{ij}$ we have that $m_{ij} = \mathcal{O}(h^2),$ we have for all $\chi\in\Sh,$
\begin{align*}
 \widehat{d}_{\M_L,h}(\zeta,\eta;\chi,\chi) &\leq  Ch^2\sum_{K\in\Th}\sum_{Z_i,Z_j\in \Zh(K)} ( \chi_i - \chi_j)^2 \leq Ch^2\|\nabla \chi\|_{L^{2}}^2. 
\end{align*}
In view \eqref{Schwartz_ineq_afc_2D}, we get the desired result. Similarly, can be also obtained the remaining estimates.
\end{proof}

\section{Fully--discrete scheme}\label{section:fully_discrete}

We now turn our attention to the fully--discrete scheme. For the temporal discretization, we will use the backward Euler method in a uniform partition of the temporal domain, i.e., let $\NT\in\mathbb{N}$, $\NT\ge1$, $k=T/\NT$ and $t^n=nk$, $n=0,\dots, \NT$. First, we define a time--discrete admissible set $\Uadk\subset[\Uad]^{\NT},$  be such that
\begin{align}\label{Uad_subset}
\Uadk :=   \left\lbrace \bv = (v^1,\ldots,v^{\NT})^T\in [L^{2}(\Omega)]^{\NT}\,:\, \mathsf{u_a} \leq v^{n} (:= v(\x,t^n)) \leq \mathsf{u_b},\;\;\text{a.e in }\Omega,\;n=1,\ldots,\NT \right\rbrace.
\end{align}
We seek $Y^n\in \Sh$, approximation of $\overline y^n=\overline y(\cdot,t^n)$ for $n=1,\dots,\NT$, such that, for all $\chi\in\Sh,$
\begin{equation}\label{optimal_AFC_fully_weak_y}
(\overline{\partial}Y^n, \chi)  + \alpha_h(Y^n , \chi) - \widehat d_{\D,h}(Y^n;Y^n,\chi) + \widehat d_{\M_L,h}(Y^n-Y^{n-1},Y^n;\overline{\partial}Y^n,\chi) = (BU^n + G^n, \chi),
\end{equation}
with $Y^0 = y_h^0\in\Sh,\,\overline{\partial} Y^n = (Y^n - Y^{n-1})/k$ and $G^n := G(\cdot,t^n).$ Moreover, we seek $P^{n-1}\in \Sh,$ approximation of $\overline p^{n-1} = \overline p(\cdot,t^{n-1}),$ for $n=\NT,\dots,1$, such that, for all $\chi\in\Sh,$
\begin{equation}\label{optimal_AFC_fully_weak_p}
- (\overline{\partial}P^n, \chi)  + \alpha_h(\chi, P^{n-1}) + \widehat d_{\widehat \D,h}(P^{n-1};P^{n-1},\chi) - \widehat d_{\M_L,h}(P^{n-1} - P^n,P^{n-1};\overline{\partial}P^{n},\chi) = (Y^n - y^n_d,\chi),
\end{equation}
with $P^{\NT} = 0.$  The control variable $\bU\in\Uadk,$ is given by
\begin{align}
(\lambda\, U^n +  B^*P^{n-1}, \psi - U^{n}) & \geq 0, \;\;\;\;\forall\,\psi\in L^{2},\;\text{with}\;\mathsf{u_a}\leq \psi \leq \mathsf{u_b},\;\text{a.e. in}\;\Omega.\label{optimal_AFC_fully_weak_u}
\end{align}
Similarly to \eqref{optimal_weak_u2}, one can express \eqref{optimal_AFC_fully_weak_u} by using \eqref{operator_Pi}, as
\begin{align}\label{optimal_AFC_fully_weak_u2}
U^n := \Pi_{[\mathsf{u_a},\mathsf{u_b}]}(-\lambda^{-1} B^*P^{n-1}),\;n = \NT,\ldots,1.
\end{align}
 
The resulting fully--discrete scheme, \eqref{optimal_AFC_fully_weak_y}--\eqref{optimal_AFC_fully_weak_u2} is nonlinear, therefore we need to prove its existence and uniqueness. Before we prove the existence and uniqueness, we will prove that the stabilized scheme satisfies the local discrete maximum principle. To do this, let us write its matrix formulation. We introduce the following notation. Let $\al^n=(\alpha^n_1,\dots,\alpha^n_{\N})^T,\,\be^n=(\beta^n_1,\dots,\beta^n_{\N})^T$ the coefficients, with respect to the basis of $\Sh$ of $Y^n,\,P^n\in\Sh,$ respectively. Then \eqref{optimal_AFC_fully_weak_y}--\eqref{optimal_AFC_fully_weak_p} can be written as
\begin{align}
(\M_L +  k\,(\mu\S + \Q + \D + \sigma\,\M_L)) \al^n & = \M_L\al^{n-1} + k\,\bfr(BU^n + G^n) +  k\,\overline{\mathsf{f}}(\al^n) + \overline{\mathsf{g}}(\al^n-\al^{n-1},\al^{n}),\label{fully_AFC_matrix_y}\\
(\M_L + k\,(\mu\S - \Q - \widehat\D + \sigma\,\M_L)) \be^{n-1} & = \M_L\be^{n} + k\,\M \al^n - k\,\bfr( y^n_d) +  k\,\overline{\widehat{\mathsf{f}}}(\be^{n-1}) + \overline{\widehat{\mathsf{g}}}(\be^{n-1}-\be^{n},\be^{n-1}),\label{fully_AFC_matrix_p}
\end{align}
with $\al^0 = \widehat v,\,\be^{\NT} = \mathbf{0},$ where $\widehat v,\,\mathbf{0}$ are the coefficient vectors of $y_h^0,\,p_h^{\NT}\in\Sh.$ \par

The next definition describes how the correction factors are computed.

\begin{definition}\label{definition:correction_factors_AFC}
The correction factors in the stabilization terms in \eqref{fully_AFC_matrix_y},\eqref{fully_AFC_matrix_p}, are computed using Algorithm \ref{algorithm-1}. More specifically, for a finite element function $\psi^n\in\Sh,$ at time step $t=t^n,$ with coefficient vector $\varthet^n\in\mathbb{R}^{\N},$ i.e., $\psi^n = \sum_{j=1}^{\N}\vartheta_j^n\phi_j,$ the correction factors $\mathfrak{a}_{ij}(\psi^n),\,\widehat{\mathfrak{a}}_{ij}(\psi^n)$, $i,j=1,\ldots,\N,$ are computed as follows.
The $\mathfrak{a}_{ij}(\psi^n)$ are computed from Algorithm \ref{algorithm-1} with $Q^{\pm}(\varthet^n),\,P^{\pm}(\varthet^n)$ and $q_i = \gamma_i\sum_{j\in\Zh^i}d_{ij}.$ In addition, the $\widehat{\mathfrak{a}}_{ij}(\psi^n)$ are computed from Algorithm \ref{algorithm-1} with $Q^{\pm}(\varthet^n),\,P^{\pm}(\varthet^n)$ and $q_i = \gamma_i\sum_{j\in\Zh^i}\widehat d_{ij}.$ Further, in view of Definition \ref{definition:corr_factors_gen}, we define the correction factors $\mathrm{a}_{ij}(\psi^n-\psi^{n-1},\psi^n):=\mathtt{a}_{ij}(\psi^n-\psi^{n-1},\psi^{n-1})$ which are computed from Algorithm \ref{algorithm-1} with $Q^{\pm}(\varthet^n),\,P^{\pm}(\varthet^n-\varthet^{n-1})$ and $q_i = \gamma_i\sum_{j\in\Zh^i}m_{ij}.$ Similarly, $\widehat{\mathrm{a}}_{ij}(\psi^{n-1}-\psi^{n},\psi^{n-1}) :=\mathtt{a}_{ij}(\psi^{n-1}-\psi^{n},\psi^{n-1}),\,i,j=1,\ldots,\N,$ are computed from Algorithm \ref{algorithm-1} with $Q^{\pm}(\varthet^{n-1}),\,P^{\pm}(\varthet^{n-1}-\varthet^{n})$ and $q_i = \gamma_i\sum_{j\in\Zh^i}m_{ij}.$
\end{definition}

Every interior node $i\in\Nh^0,$ the $i-$th equation of the fully--discrete scheme \eqref{fully_AFC_matrix_y}--\eqref{fully_AFC_matrix_p} at $t=t^n$ can be written as follows,
\begin{align}
( m_i + k\,\sigma\,m_i) \alpha_i^n + \sum_{j\neq i} (\mu s_{ij} + \tau_{ij} + d_{ij}) (\alpha_j^n - \alpha_i^n) & = b_i^n,\label{fully_AFC_matrix_y_i}\\
( m_i + k\,\sigma\,m_i) \beta_i^{n-1} + \sum_{j\neq i} (\mu s_{ij} - \tau_{ij} - \widehat d_{ij}) (\beta_j^{n-1} - \beta_i^{n-1}) & = \widehat b_i^n,\label{fully_AFC_matrix_p_i}
\end{align}
where
\begin{align*}
b_i^n & := m_i\alpha_i^{n-1} + k\,r_i(G^n + BU^n) + k\,\overline{\mathsf{f}}_i^n + \overline{\mathsf{g}}_i^n,\\
\widehat b_i^n & := m_i\beta_i^{n} + k\,\sum_{j=1}^{\N}\,m_{ij}\alpha_j^n - k\,r_i(y_d^n) + k\,\overline{\widehat{\mathsf{f}}}_i^n + \overline{\widehat{\mathsf{g}}}_i^n.
\end{align*}
To obtain the representation \eqref{fully_AFC_matrix_y_i}--\eqref{fully_AFC_matrix_p_i}, we have used the fact that the matrices $\S,\,\Q,\,\D$ have zero--row sum.

Next, we define the local discrete maximum principle (DMP) for the system \eqref{fully_AFC_matrix_y_i}--\eqref{fully_AFC_matrix_p_i}, in analogy with the local DMP for stationary problems, see, e.g., \cite[Theorem 3.1]{barrenechea2017b}, \cite[Theorem 3.3]{barrenechea2025}, \cite[Theorem 2]{barrenechea2018}, \cite[Definition 1 and Theorem 1]{knobloch2017}. For nonlinear discretizations \cite[Theorem 3.18]{barrenechea2024} and in the context of optimal control problem governed by an elliptic convection--diffusion--reaction equation, \cite[Theorem 5.16]{baumgartner2022}. A similar to the following definition can also be found in \cite[Lemmas 8.2, 8.12]{barrenechea2024},  \cite[Definition 4.21]{lohmann2019}.

\begin{definition}(Local DMP)\label{definition:DMP_fully}
We say that the solution $\al^n \in \mathbb{R}^{\N}$ of \eqref{fully_AFC_matrix_y_i} with $\sigma>0$ in $\omega_i$ at $t=t^n,$ satisfies the local DMP if for $i\in\Nh^0,$ we have
\begin{align}
G^n + BU^n & \leq 0 \;\;\text{in}\;\;\omega_i  \Rightarrow \;\alpha_i^n \leq \alpha_i^{\max,+},\label{local_dmp1}\\
G^n + BU^n & \geq 0  \;\;\text{in}\;\;\omega_i  \Rightarrow \;\alpha_i^n \geq \alpha_i^{\min,-},\label{local_dmp2}
\end{align}
where 
\begin{align*}
\alpha_i^{\max,+} &:= \max\left\{0, \alpha_i^{\max}\right\},\;\;\text{with}\;\;\alpha_i^{\max} := \max\left\{\max_{j\in\Zh^i\cup\{i\}}\alpha_j^{n-1}, \max_{j\in\Zh^i} \alpha_j^n\right\}\\
\alpha_i^{\min,-} &:= \min\left\{0, \alpha_i^{\min}\right\},\;\;\;\text{with}\;\;\alpha_i^{\min} := \min\left\{\min_{j\in\Zh^i\cup\{i\}}\alpha_j^{n-1}, \min_{j\in\Zh^i} \alpha_j^n\right\}.
\end{align*}
If $\sigma = 0$ in $\omega_i,$ then,
\begin{align}
G^n + BU^n & \leq 0 \;\;\text{in}\;\;\omega_i  \Rightarrow \;\alpha_i^n \leq \alpha_i^{\max},\label{local_dmp1_2}\\
G^n + BU^n & \geq 0  \;\;\text{in}\;\;\omega_i  \Rightarrow \;\alpha_i^n \geq \alpha_i^{\min}.\label{local_dmp2_2}
\end{align}
Furthermore, we say that the solution $\be^{n-1} \in \mathbb{R}^{\N}$ of \eqref{fully_AFC_matrix_p_i} $\sigma>0$ in $\omega_i$ at $t=t^{n-1}$, satisfies the local DMP if for $i\in\Nh^0,$ we have
\begin{align}
Y^n - y_d^n & \leq 0  \;\;\text{in}\;\;\omega_i  \Rightarrow \;\beta_i^{n-1} \leq \beta_i^{\max,+},\label{local_dmp3}\\
Y^n - y_d^n & \geq 0  \;\;\text{in}\;\;\omega_i \Rightarrow \;\beta_i^{n-1} \geq \beta_i^{\min,-},\label{local_dmp4}
\end{align}
where 
\begin{align*}
\beta_i^{\max,+} &:= \max\left\{0, \beta_i^{\max}\right\},\;\;\text{with}\;\;\beta_i^{\max} := \max\left\{\max_{j\in\Zh^i\cup\{i\}}\beta_j^{n}, \max_{j\in\Zh^i} \beta_j^{n-1}\right\}\\
\beta_i^{\min,-} &:= \min\left\{0, \beta_i^{\min}\right\},\;\;\;\text{with}\;\;\beta_i^{\min} := \min\left\{\min_{j\in\Zh^i\cup\{i\}}\beta_j^{n}, \min_{j\in\Zh^i} \beta_j^{n-1}\right\}.
\end{align*}
If $\sigma = 0$ in $\omega_i,$ then,
\begin{align}
Y^n - y_d^n & \leq 0  \;\;\text{in}\;\;\omega_i  \Rightarrow \;\beta_i^{n-1} \leq \beta_i^{\max},\label{local_dmp3_2}\\
Y^n - y_d^n & \geq 0  \;\;\text{in}\;\;\omega_i \Rightarrow \;\beta_i^{n-1} \geq \beta_i^{\min}.\label{local_dmp4_2}
\end{align}
\end{definition}

In the next theorem, we will show that the fully--discrete scheme \eqref{fully_AFC_matrix_y}--\eqref{fully_AFC_matrix_p} satisfies the local DMP.  

\begin{theorem}\label{theorem:dmp_fully}
Assume the correction factors are computed as in Definition \ref{definition:correction_factors_AFC}. Then, in view of Definition \ref{definition:DMP_fully}, the solution let $(Y^n,P^{n-1},U^n)$ of the fully--discrete scheme \eqref{optimal_AFC_fully_weak_y}--\eqref{optimal_AFC_fully_weak_u2} satisfies the local DMP for all $k,\,h.$
\end{theorem}
\begin{proof}
First, we will prove the case \eqref{local_dmp1}. 
Let $i\in\Nh^0,$ we assume that $G^n + BU^n \leq 0$ in $\omega_i.$ Then, $r_i(G^n + BU^n) \leq 0,$ since $\phi_i>0,\,i\in\Nh^0.$ Let that \eqref{local_dmp1} does not hold, i.e., $\alpha_i^{n} > \alpha_i^{\max} \geq 0$ and $\alpha_i^n = \max_{j\in\Zh^i\cup \{i\}}\alpha_j^n,$ i.e., $\alpha_i^n$ is a local maximum on $\overline\omega_i.$ 
Then, in view of \eqref{fully_AFC_matrix_y_i}, the fact that the off-diagonal elements of $\mu\S + \Q + \D$ are non-positive and the latter matrix has zero row-sum, we have
\begin{align*}
( m_i + k\,\sigma\,m_i) \alpha_i^n & \leq  m_i\alpha_i^{n-1} + \overline{\mathsf{f}}_i^n + \overline{\mathsf{g}}_i^n \leq  m_i\alpha_i^{\max} + \overline{\mathsf{f}}_i^n + \overline{\mathsf{g}}_i^n, 
\end{align*}
since $\alpha_i^{\max} \geq \alpha_i^{n-1},$ due to its definition, see \eqref{local_dmp1}. In addition, because $\alpha_i^n$ is a local maximum on $\overline\omega_i,$ we have that $\overline{\mathsf{f}}_i^n = \overline{\mathsf{g}}_i^n = 0,$ by the application of the Algorithm \ref{algorithm-1}. Then, we directly get that $( m_i + k\,\sigma\,m_i) \alpha_i^n  \leq  m_i\alpha_i^{\max}.$ Due to the fact that $\sigma > 0,$ we obtain that $\alpha_i^n \leq \alpha_i^{\max},$ which cannot be hold due to the hypothesis, therefore, we arrive in a contradiction. The cases \eqref{local_dmp2}--\eqref{local_dmp4_2} can be treated analogously.
\end{proof}

Next, we prove that \eqref{fully_AFC_matrix_y}--\eqref{fully_AFC_matrix_p} satisfies the global DMP. 

\begin{definition}[Global DMP]\label{definition:global_DMP}
We say that the solution $\al^n\in\mathbb{R}^{\N}$ of \eqref{fully_AFC_matrix_y_i} with $\sigma>0$ in $\Omega,$ satisfies the global DMP if for $i\in\Nh^0,$ we have
\begin{align}
G^n + BU^n & \leq 0 \;\;\text{in}\;\;\Omega  \Rightarrow \;\alpha_i^n \leq \alpha^{\max,n-1,+},\label{global_dmp1}\\
G^n + BU^n & \geq 0  \;\;\text{in}\;\;\Omega  \Rightarrow \;\alpha_i^n \geq \alpha^{\min,n-1,-},\label{global_dmp2}
\end{align}
where 
\begin{align*}
\alpha^{\max,n-1,+} &:= \max\left\{0, \alpha^{\max,n-1}\right\},\;\;\text{with}\;\;\alpha^{\max,n-1} := \max_{j\in\Nh}\alpha_j^{n-1}\\
\alpha^{\min,n-1,-} &:= \min\left\{0, \alpha^{\min,n-1}\right\},\;\;\;\text{with}\;\;\alpha^{\min,n-1} := \min_{j\in\Nh}\alpha_j^{n-1}.
\end{align*}
If $\sigma = 0$ in $\Omega,$ then,
\begin{align}
G^n + BU^n & \leq 0 \;\;\text{in}\;\;\Omega  \Rightarrow \;\alpha_i^n \leq \alpha^{\max,n-1},\label{global_dmp1_2}\\
G^n + BU^n & \geq 0  \;\;\text{in}\;\;\Omega  \Rightarrow \;\alpha_i^n \geq \alpha^{\min,n-1}.\label{global_dmp2_2}
\end{align}
Furthermore, we say that the solution $\be^{n-1} \in \mathbb{R}^{\N}$ of \eqref{fully_AFC_matrix_p_i} $\sigma>0$ in $\Omega,$ satisfies the global DMP if for $i\in\Nh^0,$ we have
\begin{align}
Y^n - y_d^n & \leq 0  \;\;\text{in}\;\;\Omega  \Rightarrow \;\beta_i^{n-1} \leq \beta^{\max,n,+},\label{global_dmp3}\\
Y^n - y_d^n & \geq 0  \;\;\text{in}\;\;\Omega \Rightarrow \;\beta_i^{n-1} \geq \beta^{\min,n,-},\label{global_dmp4}
\end{align}
where 
\begin{align*}
\beta^{\max,n,+} &:= \max\left\{0, \beta^{\max,n}\right\},\;\;\text{with}\;\;\beta^{\max,n} := \max_{j\in\Nh}\beta_j^{n}\\
\beta^{\min,n,-} &:= \min\left\{0, \beta^{\min,n}\right\},\;\;\;\text{with}\;\;\beta^{\min,n} := \min_{j\in\Nh}\beta_j^{n}.
\end{align*}
If $\sigma = 0$ in $\Omega,$ then,
\begin{align}
Y^n - y_d^n & \leq 0  \;\;\text{in}\;\;\Omega  \Rightarrow \;\beta_i^{n-1} \leq \beta^{\max,n},\label{global_dmp3_2}\\
Y^n - y_d^n & \geq 0  \;\;\text{in}\;\;\Omega \Rightarrow \;\beta_i^{n-1} \geq \beta^{\min,n}.\label{global_dmp4_2}
\end{align}
\end{definition}

In the next theorem, we will show that the fully--discrete scheme \eqref{fully_AFC_matrix_y}--\eqref{fully_AFC_matrix_p} satisfies the global DMP.  

\begin{theorem}\label{theorem:global_dmp_fully}
Assume the correction factors are computed as in Definition \ref{definition:correction_factors_AFC}. Then, in view of Definition \ref{definition:global_DMP}, the solution let $(Y^n,P^{n-1},U^n)$ of the fully--discrete scheme \eqref{optimal_AFC_fully_weak_y}--\eqref{optimal_AFC_fully_weak_u2} satisfies the global DMP for all $k,\,h.$
\end{theorem}
\begin{proof}
First, we will prove the case \eqref{global_dmp1}. 
Let $i\in\Nh^0,$ we assume that $G^n + BU^n \leq 0$ in $\Omega.$ Then, $r_i(G^n + BU^n) \leq 0,$ since $\phi_i>0,\,i\in\Nh^0.$ We assume that $\alpha_i^n = \max_{j\in\Nh}\alpha_j^n > \alpha^{\max,n-1,+}\geq 0,$ otherwise \eqref{global_dmp1} holds trivially.

Similar to the proof of Theorem \ref{theorem:dmp_fully}, the off-diagonal elements of $\mu\S + \Q + \D$ are non-positive and the latter matrix has zero row-sum, hence,
\begin{align*}
( m_i + k\,\sigma\,m_i) \alpha_i^n & \leq  m_i\alpha_i^{n-1} + \overline{\mathsf{f}}_i^n + \overline{\mathsf{g}}_i^n \leq  m_i\alpha_i^n + \overline{\mathsf{f}}_i^n + \overline{\mathsf{g}}_i^n, 
\end{align*}
since $\alpha_i^{\max,n-1,+} < \alpha_i^{n}.$
In addition, because $\alpha_i^n$ is a global maximum on $\Omega$  and thus local maximum on $\overline\omega_i,$ we have that $\overline{\mathsf{f}}_i^n = \overline{\mathsf{g}}_i^n = 0,$ by the application of the Algorithm \ref{algorithm-1}. Then, we directly get that $( m_i + k\,\sigma\,m_i) \alpha_i^n  \leq  m_i\alpha_i^n.$ Due to the fact that $\sigma > 0,$ we arrive in a contradiction. The cases \eqref{global_dmp2}--\eqref{global_dmp4_2} can be treated analogously.
\end{proof}

\subsection{Error estimates for the state and co-state variable}

For our error analysis purposes, we require that the solution $(\overline{y}, \overline{p}, \overline{u})$ of \eqref{optimal_weak_y}--\eqref{optimal_weak_u2} is sufficiently smooth. This regularity is presented on the Lemma \ref{lemma:regularity} and it can be ensured by the following assumption.

\begin{assumption}\label{assumption:u_regularity}
We assume that $\Omega \subset \mathbb{R}^2$ is a bounded convex polygonal domain with Lipschitz boundary $\partial \Omega$ and $G,\,y_d \in H^{1}(0,T; L^{2})$, $y_0 \in H^{2}$, $\mathbf{b}(\mathbf{x},t),\,\mathbf{b}_t(\mathbf{x},t) \in [L^{\infty}]^2$ for all $t \in [0,T].$ For the control operator we assume that $B\overline u \in H^{1}(0,T;L^{2}).$
 In addition, the following compatibility conditions are true, $G(0),\,B\overline u(0),\, y_d(T),\,\Delta\,\overline y_0,\,\bfb(\cdot,0)\cdot\nabla \overline y_0\in H^{1}_0.$
\end{assumption}

We refer to Remark \ref{remark:reg_example} for some examples of control operators $B$ that satisfy the above assumption. In view of the latter assumption, the following lemma can be proved concerning the regularity of \eqref{optimal_weak_y}--\eqref{optimal_weak_u2}.

\begin{lemma}\label{lemma:regularity}
In view of the Assumption \ref{assumption:u_regularity}, for the unique solution $(\overline y, \overline p, \overline u)$ of \eqref{optimal_weak_y}--\eqref{optimal_weak_u2}, the following estimates are hold,
\begin{equation} 
\begin{aligned}
\overline u &\in L^{2}(0,T;W^{1}_p) \cap H^{1}(0,T;L^{2}),\;\; p<\infty,\\
\overline y,\,\overline p &\in L^{\infty}(0,T;H^{2})\cap H^{2}(0,T;L^{2}), \quad
\overline y_t,\,\overline p_t \in L^{2}(0,T;H^{2})\cap L^{\infty}(0,T;H^{1}_0).
\end{aligned}
\end{equation}
In particular, there exist a constant $M>0,$ that may depend on $\mu^\nu,\,\nu<0,$ such that for $0\leq t\leq T,$
\begin{equation}\label{M0-bound}
\begin{aligned}
& \max_{0\leq s\leq t} \left(\|\overline y(s)\|_{2}^2 + \|\overline p(s)\|_{2}^2 + \|\overline y_t(s)\|_{1}^2 + \|\overline p_t(s)\|_{1}^2 \right)\\
&\qquad\quad + \int_0^t\left( \|\overline u(s)\|_{1,p}^2 + \|\overline y_t(s)\|_2^2  + \|\overline p_t(s)\|_2^2 + \|\overline y_{tt}(s)\|_{L^{2}}^2  + \|\overline p_{tt}(s)\|_{L^{2}}^2 \right)\,ds \leq M.
\end{aligned}
\end{equation}
\end{lemma}
\begin{proof}
The proof uses standard techniques, see, e.g., \cite[Chapter 7]{evans2010} in combination with the stability estimate for the pointwise projection \eqref{operator_Pi}, see, e.g., \cite[Equation (2.9)]{meidner2011}, \cite[Equation (3.5)]{jin2020}, 
\begin{align}\label{Pi_stab_estimate}
\|\Pi_{[\mathsf{u_a},\mathsf{u_b}]} v\|_{W^{s}_{p}(0,T;L^{2})} \leq \|v\|_{W^{s}_{p}(0,T;L^{2})},\;\;\forall\,v\in {W^{s}_{p}(0,T;L^{2})},\;\;0\leq s\leq 1,\,\;1 \leq p \leq \infty.
\end{align}
Similar proofs for the optimal control problem governed by the heat equation can be found in \cite[Propositions 2.1 and 2.3]{meidner2011}, \cite[Lemma 2]{springer2014} and \cite[Lemma 2.2]{zhang2023}. For the semilinear heat equation, a related proof is given in \cite[Proposition 2.1]{neitzel2012}. For completeness reasons, we describe the proof briefly.\par

Since $G,\,B\overline{u}\in L^{2}(0,T;L^{2}),$ we get, see, \cite[Section 7.1.3, Theorem 5(i)]{evans2010}, that $\overline y \in L^{2}(0,T;H^{2})\cap  L^{\infty}(0,T;H^{1}_0)$ with $\overline y_t \in L^{2}(0,T;L^{2}).$ Next, since $\overline y - y_d \,\in L^{2}(0,T;L^{2}),$ we also have $\overline p \in L^{2}(0,T;H^{2})\cap  L^{\infty}(0,T;H^{1}_0)$ with $\overline p_t \in L^{2}(0,T;L^{2}).$ \par

From the embedding $H^{2} \hookrightarrow W^{1}_p,\,p<\infty$ and using \eqref{Pi_stab_estimate}, we obtain the estimate for the control $\overline u.$ Further, in view of the Assumption \ref{assumption:u_regularity}, $G,\,B\overline u\in H^{1}(0,T;L^{2}),$ we obtain, see, \cite[Section 7.1.3, Theorem 5(ii)]{evans2010} that $\overline y \in L^{\infty}(0,T;H^{2})$ with $\overline y_t \in L^{\infty}(0,T;L^{2})\cap L^{2}(0,T;H^{1}_0).$ To derive the estimate for the second derivative with respect to time, we take the derivative with respect to time on \eqref{optimal_weak_y}, to get
\begin{align*}
(\overline y_{tt} , v) + \alpha(\overline y_t, v) = ((B\overline u)_t + G_t - \bfb_t \cdot \nabla \overline y,v),\qquad  \forall\, v\in H^{1}_0.
\end{align*}
Testing $v=\overline y_{tt}$ and using standard estimates, see \cite[Section 7.1.3, Theorem 5(i)]{evans2010}, as well as the compatibility condition $\overline y_t(0)=B\overline u(0) + G(0) + \mu\,\Delta\,\overline y_0 - \bfb(\cdot,0)\cdot\nabla \overline y_0 - \sigma \,\overline y_0 \in H^{1}_0,$ we obtain that $\overline y_{t} \in L^{2}(0,T;H^{2})\cap L^{\infty}(0,T;H^{1}_0),\,\overline y_{tt} \in L^{2}(0,T;L^{2}).$  Note that we have used Assumption \ref{assumption:u_regularity}.\par

In addition, $\overline y - y_d \in H^{1}(0,T;L^{2}),$ we obtain that $\overline p \in L^{\infty}(0,T;H^{2})$ with $\overline p_t \in L^{\infty}(0,T;L^{2})\cap L^{2}(0,T;H^{1}_0).$ To obtain the desired regularity for $\overline p_t,\,\overline p_{tt},$ we proceed using similar arguments as before, together with the embedding $L^{2}(0,T;H^{2}\cap H^{1}_0)\cap H^{1}(0,T;L^{2}) \hookrightarrow \CC([0,T];H^{1}_0)$ from which we conclude to the following compatibility condition at $t=T,$ $\overline y(T)-y_d(T)\in H^{1}_0.$
\end{proof}

\begin{remark}\label{remark:reg_example}
A choice for the control operator $B$ that satisfy Assumption \ref{assumption:u_regularity}, is $(B\overline u)(\x,t) = \overline u(\x,t)$ for $(\x,t)\in \Omega \times [0,T]$ and $(B\overline u)(\x,t) = 0,\,(\x,t)\in\partial\Omega \times [0,T].$ Note that, in view of Lemma \ref{lemma:regularity}, we have that $\overline p \in L^{\infty}(0,T,H^{2}),$ which can be ensured if $y_d^n,\,\overline y\in H^{1}(0,T;L^{2}).$ While the former is a hypothesis, the latter can be obtained if $G,\,B\overline u\in L^{2}(0,T;L^{2}),$ which are both true. Thus, using \eqref{Pi_stab_estimate}, we obtain that $\overline u \in L^{2}(0,T;W^{1}_p) \cap H^{1}(0,T;L^{2}),\,p<\infty,$ and hence the $B\overline u \in L^{2}(0,T;W^{1}_p) \cap H^{1}(0,T;L^{2}),\,p<\infty.$  Further,  in view of \eqref{Pi_stab_estimate} and the above estimate for $\overline p,$ we obtain $\overline u(t)\in H^{1},\,t\in [0,T].$ For $\mathsf{u_a}\mathsf{u_b}\leq 0,$ we further get that $\min\{\mathsf{u_b}, - \lambda^{-1}p(\x,t)\} = 0,\,\x \in\partial\Omega,\,t\in [0,T].$ Thus, $\overline u(t)\in H^{1}_0,\,t\in [0,T].$ This ensures the compatibility condition that we have used in the proof of Lemma \ref{lemma:regularity}.

In addition, let $B$ represent a finite--dimensional control operator, where the controls are purely time-dependent,
\begin{align*}
(B\overline u)(\x,t) = \sum_{i=1}^D \overline u_i(t)\, g_i(\x),
\end{align*}
where $D$ is a positive integer and the functions $u_i\in L^{2}(0,T),\,g_i \in H^{1}_0,\,i=1,\ldots,D.$ Such a operator can be found e.g., in \cite{springer2014,vondaniels2015,meidner2011}.
\end{remark}

\begin{theorem}\label{theorem:error_estimates_AFC} 
Let $(\overline y,\overline p, \overline u)$ be the unique solution of \eqref{optimal_semi_weak_y}--\eqref{optimal_semi_weak_u2}, which is sufficiently regular by Lemma \ref{lemma:regularity} and let $(Y^n,P^{n-1},U^n)$ be a solution of the fully--discrete scheme \eqref{optimal_AFC_fully_weak_y}--\eqref{optimal_AFC_fully_weak_u2} with $Y^0\in\Sh,$ such that $\|Y^0 - \overline y^0\|_{L^{2}} + h\tribar Y^0 - \overline y^0 \tribar  = \mathcal{O}(h^2).$ Assume the correction factors are computed as in Definition \ref{definition:correction_factors_AFC}. Then, for $k,\,h$ sufficiently small with $k = \mathcal{O}(h),$ there exists $C_{M,\mu}>0$, independent of $k,\,h,$ depending on $M,\,\mu^\nu,\,\nu<0,$ such that at time level $t=t^n$ for the state and at $t=t^{n-1}$ for the co-state function, we have,
\begin{equation}\label{error_estimate_afc_state}
\begin{aligned}
& \|Y^n - R_h\overline y^n\|_{L^{2}} + \|P^{n-1} - R_h\overline p^{n-1}\|_{L^{2}} \\
&\quad + \left( k \sum_{n=1}^{\NT}\left(\tribar Y^n - R_h\overline y^n \tribar^2 + \tribar P^{n-1} - R_h\overline p^{n-1}\tribar^2 \right)\right)^{1/2} \leq C_{M,\mu}\left(k\sum_{j=1}^{\NT}\|U^j - \overline u^j\|_{L^{2}}^2 + k^2 + h^4 \right)^{1/2},\\
& \tribar Y^n - R_h\overline y^n \tribar + \tribar P^{n-1} - R_h\overline p^{n-1} \tribar \leq C_{M,\mu}\left(k\sum_{j=1}^{\NT}\|U^j - \overline u^j\|_{L^{2}}^2 + k^2 + h^2 + k^{-1}h^4\right)^{1/2}.
\end{aligned}
\end{equation}
\end{theorem}
\begin{proof}
Define the finite element functions $\theta^{n}_y,\,\theta^{n}_p\in \Sh,$ defined as $\theta^{n}_y  = Y^{n} - R_h\overline y^{n},\,\theta^{n}_p  = P^{n} - R_h\overline p^{n}$ and  $\rho^{n}_y  = R_h\overline y^{n}-\overline y^{n},\,\rho^n_p  = R_h\overline p^n-\overline p^n$, for $n\ge0$, where $R_h\,:\,H^{1} \to \Sh,$ the usual elliptic projection defined in \eqref{ritz_projection_2D}. \par
Combining \eqref{optimal_weak_y} and \eqref{optimal_AFC_fully_weak_y}, the error equation for $\theta^{n}_y$, is for all $\chi\in\Sh,$
\begin{equation}\label{error_eq_1}
\begin{aligned}
(\overline{\partial}\theta^n_y, \chi)& + \alpha_h(\theta^n_y , \chi) + \widehat d_{\M_L,h}(Y^n-Y^{n-1},Y^n;\overline{\partial} \theta_y^n,\chi)  =  (\overline y^n_t - \overline{\partial}R_h \overline y^n, \chi) + \widehat d_{\D,h}(Y^n;Y^n,\chi)\\
& - \widehat d_{\M_L,h}(Y^n-Y^{n-1},Y^n;\overline{\partial} R_h \overline y^n,\chi) - (B(\overline u^n - U^n), \chi)
 = I_1(\chi) + \ldots + I_4(\chi).
\end{aligned}
\end{equation}
For the first functional $I_1$ using the Taylor expansion, the stability of $R_h,$ see, e.g., \eqref{ritz_projection_stab} and its estimates \eqref{ritz_projection_est2_2D}, we get
\begin{align*}
|I_1(\chi)| \leq C_{M}(k + h^2)\|\chi\|_{L^{2}}.
\end{align*}
The second functional can be estimated in view of Lemmas \ref{lemma:stability_widehat_D_main}, \ref{lemma:regularity} with $v := \overline y^n,$ to get
\begin{align*}
|I_2(\chi)| & \leq Ch(\|\nabla (\theta^n_y + \rho_y^n)\|_{L^{2}}^2 + h^2\|\overline y^n\|_{2}^2)^{1/2}\|\nabla \chi\|_{L^{2}} \leq C_M\,(\|\theta^n_y\|_{L^{2}} + h^2)\|\nabla \chi\|_{L^{2}},
\end{align*}
where we have used also the inverse inequality \eqref{eq:inverse_estimate} and \eqref{ritz_projection_est2_2D}.
For stabilization term due to mass lumping, we use Lemma \ref{lemma:stability_widehat_M_main}, the Taylor expansion and the Lemma \ref{lemma:regularity}, to get
\begin{align*}
I_3(\chi) = \widehat d_{\M_L,h}(Y^n-Y^{n-1},Y^n;\overline{\partial} R_h \overline y^n,\chi) \leq C\,h^2 \|\overline{\partial} \nabla R_h \overline y^n\|_{L^{2}}\|\nabla \chi\|_{L^{2}} \leq C_{M}\,h^2\|\nabla \chi\|_{L^{2}}.
\end{align*}
Note now that in view of the symmetry of $(\cdot,\cdot)$, we obtain
\begin{equation}\label{identity-form}
(\overline{\partial}\theta^n_y, \theta^n_y) = \dfrac1{2k}(\|\theta^n_y\|_{L^{2}}^2-\|\theta^{n-1}_y\|_{L^{2}}^2)+\dfrac{k}2\|\overline{\partial}\theta^n_y\|_{L^{2}}^2.
\end{equation}
A similar identity holds also for $\widehat{d}_{\M_L,h}(w,s;\cdot,\cdot),\,w,s\in\Sh,$ i.e.,
\begin{equation}\label{identity-form_D}
\begin{aligned}
\widehat{d}_{\M_L,h}(w,s;\overline{\partial}\theta_y^n,\theta_y^n) & = \dfrac1{2k}\left(\widehat{d}_{\M_L,h}(w,s;\theta_y^n,\theta_y^n) - \widehat{d}_{\M_L,h}(w,s;\theta_y^{n-1},\theta_y^{n-1}) \right) + \frac{k}{2}\widehat{d}_{\M_L,h}(w,s;\overline\partial \theta_y^n ,\overline\partial \theta_y^n),
\end{aligned}
\end{equation}
and since $\widehat{d}_{\M_L,h}(w,s;\cdot,\cdot),$ forms a seminorm on $\CC,$ the last term on the last identity is non-negative.\par
Testing $\chi = \theta^n_y$ in \eqref{error_eq_1} and using the above estimates as well as the Young inequality, we get
\begin{align*}
(1 - C_{M,\mu}k)\|\theta^n_y\|_{L^{2}}^2  + k\tribar \theta^{n}_y \tribar_h^2 & + \widehat{d}_{\M_L,h}(Y^n - Y^{n-1},Y^n;\theta_y^n,\theta_y^n) - \widehat{d}_{\M_L,h}(Y^n - Y^{n-1},Y^n;\theta_y^{n-1},\theta_y^{n-1})\\
& \leq \|\theta^{n-1}_y\|_{L^{2}}^2  + C_{M,\mu}\,k\,(k^2 + h^4 + \|U^n - \overline u^n\|_{L^{2}}^2) + \frac{k}{2}\tribar \theta^n_y \tribar^2_h,
\end{align*}
where the dependence of constant $C_{M,\mu}$ from $M,\,\mu^\nu,\,\nu<0.$
Therefore, by rearranging the terms, we obtain, for small $k,$
\begin{align*}
\|\theta^n_y\|_{L^{2}}^2 & + k\tribar \theta^{n}_y \tribar_h^2 + \widehat{d}_{\M_L,h}(Y^n - Y^{n-1},Y^n;\theta_y^n,\theta_y^n)  \\
& \leq (1 + C_{M,\mu}k)(\|\theta^{n-1}_y\|_{L^{2}}^2 + \widehat{d}_{\M_L,h}(Y^n - Y^{n-1},Y^n;\theta_y^{n-1},\theta_y^{n-1})) + C_{M,\mu}\,k\,(k^2 + h^4 + \|U^n - \overline u^n\|_{L^{2}}^2).
\end{align*}
Summing from $1$ to $n \leq \NT,$ we get, for small $k,$
\begin{equation}\label{theta_y_est}
\begin{aligned}
\|\theta^n_y\|_{L^{2}}^2  & + k\,\sum_{j=1}^n\tribar \theta^{j}_y\tribar_h^2\\
& \leq \|\theta^{0}_y\|_{L^{2}}^2 + \widehat{d}_{\M_L,h}(Y^n - Y^{n-1},Y^n;\theta_y^{0},\theta_y^{0})  + C_{M,\mu}\,\left(k\sum_{j=1}^n\|U^j - \overline u^j\|_{L^{2}}^2 + k^2 + h^4\right).
\end{aligned}
\end{equation}
Next, we estimate some terms of \eqref{error_eq_1} slightly differently. In particular, all the terms expect the first stabilization term, can be treated similarly to previous calculations. For that, we use Lemmas \ref{lemma:stability_widehat_D_main}, \ref{lemma:regularity} with $v := \overline y^n,$ to get
\begin{align*}
\widehat d_{\D,h}(Y^n;Y^n,\chi) & \leq Ch(\|\nabla (\theta^n_y + \rho_y^n)\|_{L^{2}}^2 + h^2\|\overline y^n\|_{2}^2)^{1/2}\|\nabla \chi\|_{L^{2}} \leq C_M\,(\|\nabla \theta^n_y\|_{L^{2}} + h)\|\chi\|_{L^{2}}.
\end{align*}

Testing $\chi = \overline{\partial}\theta_y^n$ into \eqref{error_eq_1}, and using similar manipulations, we further obtain
\begin{align*}
\| \overline{\partial}\theta_y^n\|_{L^{2}}^2 & + \frac{1}{2k}(\tribar \theta^{n}_y \tribar_h^2 - \tribar \theta^{n-1}_y \tribar_h^2) + \frac{k}{2}\tribar \overline\partial  \theta_y^n \tribar_h^2 + \widehat d_{\M_L,h}(Y^n-Y^{n-1},Y^n;\overline \partial \theta^n_y,\overline \partial \theta^n_y)\\
& \leq C_{M,\mu}\left((k+h + \|\nabla \theta_y^n\|_{L^{2}} + \|U^n - \overline u^n\|_{L^{2}})\|\overline\partial\theta_y^n\|_{L^{2}} + h^2\|\overline\partial \nabla \theta_y^n\|_{L^{2}}\right)\\
& \leq  C_{M,\mu}\left(k^2 + h^2 + \|\nabla \theta_y^n\|_{L^{2}}^2 + \|U^n - \overline u^n\|_{L^{2}}^2\right) + \frac{1}{2}\|\overline\partial\theta_y^n\|_{L^{2}}^2 + C_{M,\mu}\,k^{-1}h^4 + \frac{\mu}{4}k\|\overline\partial \nabla \theta_y^n\|_{L^{2}}^2.
\end{align*}
Summing from $1$ to $n\leq \NT,$ rearranging the terms, using the properties of the stabilization term, i.e., its non-negativity and for small $k,$ we get
\begin{equation}\label{nabla_theta_y_est}
\begin{aligned}
k\sum_{j=1}^n\| \overline{\partial}\theta_y^j\|_{L^{2}}^2 + \tribar \theta^{n}_y\tribar_h^2 \leq \tribar \theta^{0}_y\tribar_h^2 + C_{M,\mu}\left(k\sum_{j=1}^n\|U^j - \overline u^j\|_{L^{2}}^2 + k^2 + h^2 + k^{-1}h^4\right).
\end{aligned}
\end{equation}
Next, the error for the co-state function, is given by combining \eqref{optimal_weak_p} and \eqref{optimal_AFC_fully_weak_p}, 
\begin{equation}\label{error_eq_3}
\begin{aligned}
- (\overline\partial\theta_p^{n}, \chi)  + \alpha_h(\chi,\theta^{n-1}_p) &  - \widehat d_{\M_L,h}(P^{n-1}-P^n,P^{n-1};\overline{\partial}\theta_p^n,\chi)  = - (\overline p^{n-1}_t - \overline{\partial}R_h\overline p^n, \chi) - \widehat d_{\widehat\D,h}(P^{n-1};P^{n-1},\chi)\\
& + \widehat d_{\M_L,h}(P^{n-1}-P^n,P^{n-1};\overline{\partial}R_h\overline p^n,\chi) + (\theta_y^n + \rho_y^n , \chi) = I_1(\chi) + \ldots + I_4(\chi),
\end{aligned}
\end{equation}
for all $\chi\in\Sh.$ Using similar arguments to estimate the functionals $I_1,\ldots,I_4$ as before, we get by setting $\chi = \theta_p^{n-1},$
\begin{align*}
(1 - C_{M,\mu}k) \|\theta_p^{n-1}\|_{L^{2}}^2  & + k\,\tribar \theta^{n-1}_p \tribar^2_h + \widehat d_{\M_L,h}(P^{n-1}-P^n,P^{n-1};\theta_p^{n-1},\theta_p^{n-1}) - \widehat d_{\M_L,h}(P^{n-1}-P^n,P^{n-1};\theta_p^{n},\theta_p^{n})\\
&\quad \leq  \|\theta_p^{n}\|_{L^{2}}^2 + C_{M,\mu}\,k\,(k^2 + h^4  + \|\theta_y^n\|_{L^{2}}^2 ) + \frac{k}{2}\tribar \theta_p^{n-1}\tribar^2_h.
\end{align*}
Using standard estimates and summing over $\NT$ to $n\geq 1,$ we get for small $k,$ 
\begin{align*}
\|\theta_p^{n-1}\|_{L^{2}}^2  + k\,\sum_{j=n}^{\NT}\tribar \theta^{j-1} \tribar^2_h +  \widehat d_{\M_L,h}(P^{n-1}-P^n,P^{n-1};\theta_p^{n-1},\theta_p^{n-1}) & \leq  C_{M,\mu}\left(k\sum_{j=n}^{\NT}\|\theta_y^j\|_{L^{2}}^2 + k^2 + h^4 \right).
\end{align*}
In view of the estimate \eqref{theta_y_est}, we further obtain,
\begin{align*}
k\sum_{j=1}^{\NT}\|\theta_y^j\|_{L^{2}}^2 & \leq C_{M,\mu}\,k\,\sum_{j=1}^{\NT}\left( \|\theta^{0}_y\|_{L^{2}}^2 + \widehat{d}_{\M_L,h}(Y^n - Y^{n-1}, Y^n;\theta_y^{0},\theta_y^{0}) + k\sum_{\ell=1}^j\|U^\ell - \overline u^\ell\|_{L^{2}}^2 + k^2 + h^4 \right)\\
& \leq C_{M,\mu} \left(\|\theta^{0}_y\|_{L^{2}}^2 + \widehat{d}_{\M_L,h}(Y^n - Y^{n-1}, Y^n;\theta_y^{0},\theta_y^{0}) + k\sum_{\ell=1}^{\NT}\|U^\ell - \overline u^\ell\|_{L^{2}}^2 + k^2 + h^4 \right),
\end{align*}
where we have used that $k\sum_{j=1}^{\NT}1 = T.$ Thus, we get
\begin{equation}\label{theta_p_est}
\begin{aligned}
\|\theta_p^{n-1}\|_{L^{2}}^2 + k\,\sum_{j=n}^{\NT}\tribar \theta^{j-1} \tribar^2_h & \leq C_{M,\mu}\left(\|\theta^{0}_y\|_{L^{2}}^2 + \widehat{d}_{\M_L,h}(Y^n - Y^{n-1},Y^n;\theta_y^{0},\theta_y^{0}) + k\sum_{j=1}^{\NT}\|U^j - \overline u^j\|_{L^{2}}^2 + k^2 + h^4 \right).
\end{aligned}
\end{equation}
We note that similar to \eqref{identity-form}, we can prove that
\begin{align*}
\alpha_h(\theta_p^{n-1} - \theta_p^n,\theta^{n-1}_p) & = \frac{1}{2}\left(\tribar \theta_p^{n-1} \tribar^2_h - \tribar \theta_p^{n} \tribar^2_h\right) + \frac{1}{2}\tribar \theta_p^{n-1} - \theta_p^{n} \tribar_h^2 - (\bfb^n \cdot \nabla \theta_p^{n-1}, \theta_p^{n-1} - \theta_p^n).
\end{align*}
Testing $\chi = - \overline\partial \theta_p^n$ into \eqref{error_eq_3}, and using similar calculations as the derivation of the estimate \eqref{nabla_theta_y_est}, we get
\begin{equation}\label{nabla_theta_p_est}
\begin{aligned}
k\sum_{j=n}^{\NT}\| \overline{\partial}\theta_p^j\|_{L^{2}}^2 + \tribar \theta^{n-1}_p \tribar^2_h \leq  C_{M,\mu}\left(k\sum_{j=1}^{\NT}\|U^j -\overline  u^j\|_{L^{2}}^2 + k^2 + h^2 + k^{-1}h^4\right).
\end{aligned}
\end{equation}
\end{proof}

\subsection{Error estimates for the control variable}

It remains to derive an error estimate for the control discretization. To do this, we follow standard arguments previously used for the discretization of the control variables, see, e.g., \cite{hinze2005, yan2009}, along with a treatment of the stabilization terms.

We define the following auxiliary variables. For any $\psi\in\Uad,$ we seek $(\overline y_\psi,\overline p_\psi) := (\overline y(\psi),\overline p(\psi)),$ that satisfies the following optimality conditions for $t\in[0,T],$
\begin{align}
(\overline y_{\psi,t} , v) + \alpha(\overline y_\psi ,  v)  & = ( B\psi + G,v), \;\;\forall\, v\in H^{1}_0,\label{optimal_weak_y_v}\\
- (\overline p_{\psi,t}, v) + \alpha( v,\overline p_\psi) & = (\overline y_\psi  -  y_d ,v), \;\;\;\forall\, v\in H^{1}_0,\label{optimal_weak_p_v}
\end{align}
with $\overline y_\psi(0) = y^0,\,\overline p_\psi(T) = 0.$

\begin{remark}\label{remark:auxiliary_weak_eq}
The weak formulations \eqref{optimal_weak_y}--\eqref{optimal_weak_p} and \eqref{optimal_weak_y_v}--\eqref{optimal_weak_p_v} coincide for $\psi = \overline u,$ where $\overline u\in\Uad$ as defined in \eqref{optimal_weak_cond} or equivalently in \eqref{optimal_weak_u2}.
\end{remark}

Similarly to \eqref{optimal_AFC_fully_weak_y}--\eqref{optimal_AFC_fully_weak_u2}, we can define the fully--discrete scheme for the auxiliary variables \eqref{optimal_weak_y_v}--\eqref{optimal_weak_p_v}. In particular, for given $\bpsi\in\Uadk,$ we seek $Y^n_{\bpsi}:=Y^n(\bpsi)\in \Sh$, approximation of $\overline y^n_{\psi}= \overline y_{\psi}(\cdot,t^n)$ for $n=1,\dots,\NT$, such that, for all $\chi\in\Sh,$
\begin{equation}\label{optimal_AFC_fully_weak_y_psi}
\begin{aligned}
(\overline{\partial}Y^n_{\bpsi}, \chi)  + \alpha_h(Y^n_{\bpsi},\chi)  - \widehat d_{\D,h}(Y^n_{\bpsi};Y^n_{\bpsi},\chi) + \widehat d_{\M_L,h}(Y^n_{\bpsi}-Y^{n-1}_{\bpsi},Y^n_{\bpsi};\overline{\partial}Y^n_{\bpsi},\chi)  = (B\psi^n + G^n,\chi),
\end{aligned}
\end{equation}
with $Y^0 := y_h^0\in\Sh$ and $\overline{\partial} Y^n_{\bpsi} = (Y^n_{\bpsi} - Y^{n-1}_{\bpsi})/k.$ After the computation of $Y^n_{\bpsi},\,n=1,\ldots,\NT,$ we seek at time level $t=t^{n-1},$ for known $P^n_{\bpsi} := P^n(\bpsi)\in \Sh,$ the function $P^{n-1}_{\bpsi} := P^{n-1}(\bpsi)\in \Sh,$ for $n=\NT,\dots,1$, such that, for all $\chi\in\Sh,$
\begin{equation}\label{optimal_AFC_fully_weak_p_psi}
\begin{aligned}
- (\overline{\partial}P^n_{\bpsi}, \chi) & + \alpha_h(\chi,P^{n-1}_{\bpsi})  + \widehat d_{\widehat \D,h}(P^{n-1}_{\bpsi};P^{n-1}_{\bpsi},\chi) - \widehat d_{\M_L,h}(P^{n-1}_{\bpsi} - P^n_{\bpsi}, P^{n-1}_{\bpsi};\overline{\partial}P^{n}_{\bpsi},\chi) = (Y^n_{\bpsi} - y^n_d,\chi),
\end{aligned}
\end{equation}
with $P^{\NT} := 0,$ and $\overline{\partial} P^n_{\bpsi} = (P^n_{\bpsi} - P^{n-1}_{\bpsi})/k$.

\begin{remark}\label{remark:auxiliary_fully_eq}
The weak formulations \eqref{optimal_AFC_fully_weak_y}--\eqref{optimal_AFC_fully_weak_p} and \eqref{optimal_AFC_fully_weak_y_psi}--\eqref{optimal_AFC_fully_weak_p_psi} coincide for $\bU = \bpsi,$ where $\bU\in\Uadk$ as defined in \eqref{optimal_AFC_fully_weak_u} or equivalently in \eqref{optimal_AFC_fully_weak_u2}.
\end{remark}

\begin{lemma}\label{lemma:truncation_error_fully}
Let $G,\,y_d\in H^{1}(0,T;L^{2}).$ For a $\psi\in\Uad$ with $B\psi \in H^{1}(0,T;L^{2}),$ let $(\overline y_{\psi},\overline p_{\psi})$ be the unique solution of \eqref{optimal_weak_y_v}--\eqref{optimal_weak_p_v}, which is sufficient regular, i.e., $\overline y_{\psi},\,\overline p_{\psi} \in L^{\infty}(0,T;H^{2})\cap H^{2}(0,T;L^{2})$ with $\overline y_{\psi,t},\,\overline p_{\psi,t} \in L^{2}(0,T;H^{2})\cap L^{\infty}(0,T;H^{1}_0).$ Also, let $(Y_{\bpsi}^n,P_{\bpsi}^{n-1})$ be the unique solution of \eqref{optimal_AFC_fully_weak_y_psi}--\eqref{optimal_AFC_fully_weak_p_psi}, where the correction factors computed as in Definition \ref{definition:correction_factors_AFC}. Then, for $k,\,h$ sufficiently small with $k = \mathcal{O}(h),$ there exists $C_{M,\mu}>0,$ independent of $k,\,h,$ depending on $M,\,\mu^\nu,\,\nu<0,$ such that,
\begin{align*}
& \|\overline y_{\bpsi}^n - Y_{\bpsi}^n \|_{L^{2}} + \|\overline p_{\bpsi}^{n-1} - P_{\bpsi}^{n-1} \|_{L^{2}} + \left(k\sum_{n=1}^{\NT}\left( \tribar\overline y_{\bpsi}^n - Y_{\bpsi}^n \tribar^2 + \tribar\overline p_{\bpsi}^{n-1} - P_{\bpsi}^{n-1} \tribar^2 \right)\right)^{1/2} \leq C_{M,\mu}\,(k + h^2),\\
& \tribar\overline y_{\bpsi}^n - Y_{\bpsi}^n \tribar + \tribar\overline p_{\bpsi}^{n-1} - P_{\bpsi}^{n-1} \tribar  \leq C_{M,\mu}\,(k + h + k^{-1/2}h^2).
\end{align*}
\end{lemma}
\begin{proof}
The proof follows the same arguments to the proof of Theorem \ref{theorem:error_estimates_AFC}, where now, since both \eqref{optimal_weak_y_v} and \eqref{optimal_AFC_fully_weak_y_psi} depend on the same control variable $\psi^n,$ its contribution vanishes and the estimate follows from \eqref{error_estimate_afc_state} and the elliptic projection estimate \eqref{ritz_projection_est2_2D}.
\end{proof}

\begin{remark}\label{remark:auxilliary_problem}
If $\psi := \overline u,$ the assumption on $B\psi$ on Lemma \ref{lemma:truncation_error_fully} is true, see e.g., Lemma \ref{lemma:regularity}. 
\end{remark}

\begin{theorem}\label{theorem:error_estimates_AFC_u} 
Let $(\overline y,\overline p, \overline u)$ be the unique solution of \eqref{optimal_semi_weak_y}--\eqref{optimal_semi_weak_u2}, which is sufficiently regular by Lemma \ref{lemma:regularity} and let $(Y^n,P^{n-1},U^n)$ be a solution of the fully--discrete scheme \eqref{optimal_AFC_fully_weak_y}--\eqref{optimal_AFC_fully_weak_u2}. Then, there exists $C_{M}>0,$ depending on $M,$ but independent of $k,\,h,$ such that at time level $t=t^n,$ we have for the control function
\begin{equation}\label{error_estimate_afc_control}
\begin{aligned}
\lambda\,\|U^n - \overline u^n\|_{L^{2}}^2 & \leq C_{M}\,( k + h^2 )\|\overline u^n - U^n\|_{L^{2}} + (P^{n-1} - P^{n-1}_{\overline \bu}, B(\overline u^n - U^n)).
\end{aligned}
\end{equation}
\end{theorem}
\begin{proof}
The proof uses ideas of \cite[Theorem 2.4]{hinze2005} and a similar result can be found in various works for various temporal discretizations, cf., \cite[Theorem 4.4]{zhou2010}, \cite[Theorem 3.6]{zhang2023}, \cite[Theorem 3.1]{hinze2009}.
The proof is based in \eqref{optimal_weak_cond}, \eqref{optimal_AFC_fully_weak_u}, the Taylor expansion and Lemma \ref{lemma:regularity}. More specifically, we have,
\begin{align*}
& \lambda\,\|\overline u^n - U^n\|_{L^{2}}^2  = \lambda\,(\overline u^n - U^n,\overline u^n - U^n)\\
& \; = (\lambda\,\overline u^n + B^*\overline p^{n-1}, \overline u^n - U^n) + (B^*(P^{n-1} - \overline p^{n-1}),\overline u^n - U^n) - (\lambda\,U^n + B^*P^{n-1} ,\overline u^n - U^n).
\end{align*}
For the first term, by Taylor expansion, we have
\begin{align*}
 (\lambda\,\overline u^n + B^*\overline p^{n-1}, \overline u^n - U^n) & = (\lambda\,\overline u^n + B^*\overline p^{n},\overline  u^n - U^n) + (B^*(\overline p^{n-1} - \overline p^{n}), \overline u^n - U^n)\\
&  \leq C\|\overline p^{n-1} - \overline p^{n}\|_{L^{2}}\| \overline u^n - U^n\|_{L^{2}}\\
& \leq  C\int_{t^{n-1}}^{t^n}\|\overline p_t(t)\|_{L^{2}}\,dt\,\|\overline u^n - U^n\|_{L^{2}} \leq  C_M\,k\,\|\overline u^n - U^n\|_{L^{2}},
\end{align*}
where the first term is non-positive for $\widehat u = U^n$ into \eqref{optimal_weak_cond}. In addition, to obtain the latter estimate we have used Lemma \ref{lemma:regularity}.\par

The second term can be estimated as follows,
\begin{align*}
(B^*(P^{n-1} - \overline p^{n-1}),\overline u^n - U^n) & = (P^{n-1} - P^{n-1}_{\overline \bu}, B(\overline u^n - U^n)) + (B^*(P^{n-1}_{\overline \bu} - \overline p^{n-1})),\overline u^n - U^n)\\
& \leq (P^{n-1} - P^{n-1}_{\overline \bu}, B(\overline u^n - U^n)) + C\|P^{n-1}_{\overline \bu} - \overline p^{n-1}\|_{L^{2}}\|\overline u^n - U^n\|_{L^{2}}\\
& \leq (P^{n-1} - P^{n-1}_{\overline \bu}, B(\overline u^n - U^n)) + C_{M,\mu}(k+h^2)\|\overline u^n - U^n\|_{L^{2}},
\end{align*}
where in the last step, we have used Lemma \ref{lemma:truncation_error_fully}. 

Further, using the fact that $(\lambda\,U^n + B^*P^{n-1}, \overline u^n - U^n) \geq 0$ when $\psi=\overline u^n$ into \eqref{optimal_AFC_fully_weak_u}, we get,
\begin{align*}
\lambda\,\|U^n - \overline u^n\|_{L^{2}}^2 & \leq C_{M}\,( k + h^2 )\|\overline u^n - U^n\|_{L^{2}} + (P^{n-1} - P^{n-1}_{\overline \bu}, B(\overline u^n - U^n)).
\end{align*}
\end{proof}

\begin{remark}\label{remark:lemma_comp}
Unlike the proofs in \cite[Theorem 4.4]{zhou2010} and \cite[Theorem 3.6]{zhang2023}, where the last term in \eqref{error_estimate_afc_control} was shown to be non-positive and could therefore be absorbed, in our case this cannot be established because the of the presence of the stabilization terms in the fully--discrete scheme \eqref{optimal_AFC_fully_weak_y}–\eqref{optimal_AFC_fully_weak_p} which are nonlinear.
\end{remark}

Before we prove an estimation of the control error, we derive a similar result to Lemma \ref{lemma:stability_widehat_D_difference}, that can be proved also for the stabilization terms due to mass lumping.

\begin{lemma}\label{lemma:stability_widehat_M_difference}
Assume that the correction factors are satisfying Algorithm \ref{algorithm-1} and Definitions \ref{definition:corr_factors_gen}, \ref{definition:correction_factors_AFC}. There are exists a positive constant $C_\varpi,$ that depends on \eqref{shape_regularity}, but independent of $h,\,k,$ such that at $t=t^n,\,n=1,\ldots,\NT,$ for $\zeta^n,\psi^n,\chi\in \Sh$, we have
\begin{equation}\label{stability_stab_term_widehat_M_diff}
\begin{aligned}
& \vert \widetilde d_{\M_L,h}(\psi^n - \psi^{n-1},\psi^n;\psi^n - \psi^{n-1},\chi) - \widetilde d_{\M_L,h}(\zeta^n - \zeta^{n-1},\zeta^n;\zeta^n - \zeta^{n-1},\chi)\vert\\
& + \vert \widehat d_{\M_L,h}(\psi^n - \psi^{n-1},\psi^n;\psi^n - \psi^{n-1},\chi) - \widehat d_{\M_L,h}(\zeta^n - \zeta^{n-1},\zeta^n;\zeta^n - \zeta^{n-1},\chi) \vert\\
&\qquad\qquad\qquad\qquad\qquad\qquad\qquad \leq C_{\varpi}h^2(\|\nabla (\psi^n - \zeta^n)\|_{L^{2}} + \|\nabla (\psi^{n-1} - \zeta^{n-1})\|_{L^{2}})\|\nabla \chi\|_{L^{2}}.
\end{aligned}
\end{equation}
\end{lemma}
\begin{proof}
Using the definition $\widehat{d}_{\M_L,h},$ see, e.g., \eqref{stab_term_M_widehat} and in view of Remark \ref{remark:def_remark} and Definition \ref{definition:correction_factors_AFC}, we set at $t=t^n,\,n=1,\ldots,\NT,$ the $\rho_{ij}(\widehat\psi^n,\psi^{n}) := m_{ij}(1-\mathtt{a}_{ij}(\widehat\psi^n,\psi^n))(\widehat\psi^n_i - \widehat\psi^n_j),\,\psi^n\in\Sh$ with $\widehat \psi^n:= \psi^n-\psi^{n-1},$ we get 
\begin{align*}
 \widehat d_{\M_L,h}(\widehat \psi^n,\psi^n;\widehat \psi^n,\chi) & - \widehat d_{\M_L,h}(\widehat \zeta^n,\zeta^n;\widehat \zeta^n,\chi)   = \sum_{i<j} (\rho_{ij}(\widehat\psi^n,\psi^{n}) - \rho_{ij}(\widehat\zeta^n,\zeta^{n})) (\chi_i - \chi_j),
\end{align*}
where $\widehat \zeta^n: = \zeta^n - \zeta^{n-1}$ and then an inequality similar to \eqref{local_estimate} is needed. We have
\begin{align*}
\vert \rho_{ij}(\widehat\psi^n,\psi^{n}) - \rho_{ij}(\widehat\zeta^n,\zeta^{n}) \vert & \leq \vert \rho_{ij}(\widehat\psi^n,\psi^{n}) - \rho_{ij}(\widehat\zeta^n,\psi^{n}) \vert + \vert \rho_{ij}(\widehat\zeta^n,\psi^{n}) - \rho_{ij}(\widehat\zeta^n,\zeta^{n}) \vert = I_1 + I_2.
\end{align*}
For $I_1,$ we follow the proof of \cite[Lemma 6]{barrenechea2016}, \cite[Lemma 2.13, Lemma 2.15]{chatzipantelidis2022}, where we split cases depending on the sign of $(\psi_i^n - \psi_j^n)(\zeta_i^n - \zeta_j^n).$ 
For positive sign, we have,
\begin{align*}
I_1 & = |\widehat A_{ij}(\widehat\psi^n,\psi^{n}) - \widehat A_{ij}(\widehat\zeta^n,\psi^{n})| + |\widehat B_{ij}(\widehat\psi^n) - \widehat B_{ij}(\widehat\zeta^n)| + m_{ij}|(\widehat\psi_i^n - \widehat\psi_j^n) - (\widehat\zeta_i^n - \widehat\zeta_j^n)|,
\end{align*}
where $\widehat A_{ij} = m_{ij}A_{ij},\,\widehat B_{ij} = m_{ij}B_{ij}.$
Initially, we assume that $\widehat\psi_i^n - \widehat\psi_j^n>0,\,\widehat\zeta_i^n - \widehat\zeta_j^n>0.$ We have in view of Definition \ref{definition:corr_factors_gen}, that
\begin{align*}
\vert \widehat A_{ij}(\widehat\psi^n,\psi^{n}) - \widehat A_{ij}(\widehat\zeta^n,\psi^{n}) \vert &= \vert \min\{ P_{i}^+(\widehat\psi^n), Q_{i}^+(\psi^{n})\} - \min\{ P_{i}^+(\widehat\zeta^n), Q_{i}^+(\psi^n)\}\vert \leq \vert P_{i}^+(\widehat\psi^n) - P_{i}^+(\widehat\zeta^n) \vert.
\end{align*}
Then, by calculations, we obtain
\begin{align*}
\vert P_{i}^+(\widehat\psi^n) - P_{i}^+(\widehat\zeta^n) \vert & = \left\vert \sum_{\ell\in\Zh^i}\max\{0, m_{i\ell}(\widehat\psi^n_i - \widehat\psi^n_\ell)\} - \sum_{\ell\in\Zh^i}\max\{0, m_{i\ell}(\widehat\zeta^n_i - \widehat\zeta^n_\ell)\} \right\vert \\
& \leq \sum_{\ell\in\Zh^i}m_{i\ell} \vert (\widehat\psi^n_i - \widehat\psi^n_\ell) - (\widehat\zeta^n_i - \widehat\zeta^n_\ell)\vert.
\end{align*}
Combining the above estimates, one get
\begin{align*}
I_1 \leq \sum_{\ell\in\Zh^{i}}m_{i\ell}\vert (\widehat\psi^n_i - \widehat\zeta^n_i) - (\widehat\psi^n_\ell - \widehat\zeta^n_\ell)\vert.
\end{align*}
The case where $\widehat\psi_i^n - \widehat\psi_j^n<0,\,\widehat\zeta_i^n - \widehat\zeta_j^n<0$ can be treated analogously. For non-positive sign, we have,
\begin{align*}
I_1 = \vert \rho_{ij}(\widehat\psi^n,\psi^{n}) - \rho_{ij}(\widehat\zeta^n,\zeta^{n}) \vert \leq  m_{ij}\vert (\widehat\psi^n_i - \widehat\psi^n_j) - (\widehat\zeta^n_i - \widehat\zeta^n_j)\vert,
\end{align*}
and we conclude to the same result for $I_1.$ \par

To estimate the term $I_2,$ we have from Definition \ref{definition:corr_factors_gen} and the non-negativity of $B_{ij}(\cdot),$ that
\begin{align*}
I_2 & = m_{ij}\left\vert \frac{A(\widehat\zeta^n,\psi^n) - A(\widehat\zeta^n,\zeta^n)}{|\widehat\zeta_i^n - \widehat\zeta_j^n| + B_{ij}(\widehat\zeta^n)}(\widehat\zeta_i^n - \widehat\zeta_j^n)\right\vert \leq \vert\widehat A(\widehat\zeta^n,\psi^n) -  \widehat A(\widehat\zeta^n,\zeta^n)\vert\\
& \leq \vert Q_i^+(\psi^n) - Q_i^+(\zeta^n) \vert \leq q_i\sum_{\ell\in\Zh^i\cup\{i\}}\vert \psi_\ell^n - \zeta_\ell^n\vert.
\end{align*} 
Therefore, combining all these estimates, we obtain
\begin{align*}
\vert \rho_{ij}(\widehat\psi^n,\psi^{n}) - \rho_{ij}(\widehat\zeta^n,\zeta^{n}) \vert & \leq \sum_{\ell\in\Zh^{i}}m_{i\ell}\vert (\widehat\psi^n_i - \widehat\zeta^n_i) - (\widehat\psi^n_\ell - \widehat\zeta^n_\ell)\vert + q_i\sum_{\ell\in\Zh^i\cup\{i\}}\vert \psi_\ell^n - \zeta_\ell^n\vert\\
& \leq C_{\varpi}h^2\sum_{\ell\in\Zh^i\cup\{i\}}\left( \vert \psi_\ell^n - \zeta_\ell^n\vert + \vert \psi_\ell^{n-1} - \zeta_\ell^{n-1}\vert\right),
\end{align*}
since $m_{ij} = \mathcal{O}(h^2)$ and $q_i = \gamma_i\sum_{j\neq i}m_{ij} = \mathcal{O}(h^2).$ The desired inequality follows using similar arguments to the proof of Lemma \ref{lemma:stability_widehat_D_difference} and the latter estimate.

\end{proof}

\begin{lemma}\label{lemma:truncation_error_fully2_2}
For $\overline u \in\Uad$ as defined in \eqref{optimal_weak_u2}, let $(Y_{\overline \bu}^n,P_{\overline \bu}^{n-1})$ be the unique solution of \eqref{optimal_AFC_fully_weak_y_psi}--\eqref{optimal_AFC_fully_weak_p_psi} with $\bpsi = \overline\bu.$ In addition, let $(Y^n,P^{n-1})$ be the unique solution of \eqref{optimal_AFC_fully_weak_y}--\eqref{optimal_AFC_fully_weak_u2}. We assume that in both fully--discrete schemes the correction factors computed as in Definition \ref{definition:correction_factors_AFC}. Then, for small $k,\,h$ with $k=\mathcal{O}(h),$ there exists a positive constant $C_{M,\mu}>0$, independent of $k,\,h,$ but depending on $M,\,\mu^\nu,\,\nu<0,$ such that,
\begin{align*}
k\sum_{j=1}^{\NT}\|U^j - \overline u^j\|_{L^{2}}^2 \leq C_{M,\mu}(k^2 + h^4).
\end{align*}
\end{lemma}
\begin{proof}
Define the finite element functions $\eta_y^n,\,\eta_p^n\in \Sh,$ defined as $\eta_y^n  = Y^{n}_{\overline \bu} - Y^{n},\,n=1,\ldots,\NT$ and $\eta^{n-1}_p  = P^{n-1} - P^{n-1}_{\overline \bu}$ for $n=\NT, \ldots,1$. From \eqref{optimal_AFC_fully_weak_y_psi} and \eqref{optimal_AFC_fully_weak_y}, we have 
\begin{equation}\label{eq_dummy_1}
\begin{aligned}
(B(\overline u^n - U^n),\chi) & = (\overline{\partial}\eta^n_y, \chi) + \alpha_h(\eta^n_y ,\chi) + \widehat d_{\D,h}(Y^n;Y^n,\chi) -  \widehat d_{\D,h}(Y^n_{\overline \bu};Y^n_{\overline \bu},\chi)  \\
&  + \widehat d_{\M_L,h}(Y^n_{\overline \bu}-Y^{n-1}_{\overline \bu},Y^n_{\overline \bu};\overline{\partial}Y^n_{\overline \bu},\chi) - \widehat d_{\M_L,h}(Y^n-Y^{n-1},Y^n;\overline{\partial}Y^n,\chi).
\end{aligned}
\end{equation}
Also, from \eqref{optimal_AFC_fully_weak_p_psi} and \eqref{optimal_AFC_fully_weak_p}, we have 
\begin{equation}\label{eq_dummy_2}
\begin{aligned}
- (\overline{\partial}\eta^{n}_p, \chi) & + \alpha_h(\chi, \eta^{n-1}_p) + \widehat d_{\widehat \D,h}(P^{n-1};P^{n-1},\chi) - \widehat d_{\widehat \D,h}(P^{n-1}_{\overline \bu};P^{n-1}_{\overline \bu},\chi) \\
&  + \widehat d_{\M_L,h}(P^{n-1}_{\overline \bu}-P^{n}_{\overline \bu}, P^{n-1}_{\overline \bu};\overline{\partial}P^n_{\overline \bu},\chi) - \widehat d_{\M_L,h}(P^{n-1}-P^{n}, P^{n-1};\overline{\partial}P^n,\chi) = - (\eta_y^n, \chi).
\end{aligned}
\end{equation}
We use similar arguments as in Theorem \ref{theorem:error_estimates_AFC}. To estimate the stabilization terms $\widehat d_{\D,h},\,\widehat d_{\widehat \D,h},$ we use Lemma \ref{lemma:stability_widehat_D_difference}. In particular, we have
\begin{align*}
\vert \widehat d_{\D,h}(Y^n_{\overline \bu};Y^n_{\overline \bu},\chi)  - \widehat d_{\D,h}(Y^n;Y^n,\chi) \vert & \leq Ch\|\nabla \eta_y^n\|_{L^{2}}\|\nabla \chi\|_{L^{2}},\\
\vert \widehat d_{\D,h}(P^{n-1}_{\overline \bu};P^{n-1}_{\overline \bu},\chi)  - \widehat d_{\D,h}(P^{n-1};P^{n-1},\chi) \vert & \leq Ch\|\nabla \eta_p^{n-1}\|_{L^{2}}\|\nabla \chi\|_{L^{2}}.
\end{align*}
For the stabilization term due to mass lumping, we use Lemma \ref{lemma:stability_widehat_M_difference}, to obtain, that
\begin{align*}
\vert \widehat d_{\M_L,h}(Y_{\overline \bu}^n-Y_{\overline \bu}^{n-1},Y^n_{\overline \bu};\overline{\partial}Y_{\overline \bu}^n,\chi) - \widehat d_{\M_L,h}(Y^n-Y^{n-1},Y^n;\overline{\partial}Y^n,\chi) \vert & \leq C\,\frac{h^2}{k}\left(\|\nabla \eta_y^n\|_{L^{2}} + \|\nabla \eta_y^{n-1}\|_{L^{2}}\right)\|\nabla \chi\|_{L^{2}},\\
\vert \widehat d_{\M_L,h}(P_{\overline \bu}^{n-1}-P_{\overline \bu}^{n},P_{\overline \bu}^{n-1};\overline{\partial}P_{\overline \bu}^n,\chi) - \widehat d_{\M_L,h}(P^{n-1}-P^{n},P^{n-1};\overline{\partial}P^n,\chi) \vert & \leq C\,\frac{h^2}{k}\left(\|\nabla \eta_p^n\|_{L^{2}} + \|\nabla \eta_p^{n-1}\|_{L^{2}}\right)\|\nabla \chi\|_{L^{2}}.
\end{align*}
We set $\chi = \eta_p^{n-1}$ into \eqref{eq_dummy_1} and $\chi = \eta_y^n$ in \eqref{eq_dummy_2}, to get
\begin{align*}
(B(\overline u^n - U^n),\eta_p^{n-1}) & = (\overline{\partial}\eta^n_y, \eta_p^{n-1}) + (\overline{\partial}\eta^{n}_p, \eta_y^n) + \widehat d_{\D,h}(Y^n;Y^n,\eta_p^{n-1}) -  \widehat d_{\D,h}(Y^n_{\overline \bu};Y^n_{\overline \bu},\eta_p^{n-1})  \\
&  + \widehat d_{\M_L,h}(Y^n_{\overline \bu}-Y^{n-1}_{\overline \bu},Y^n_{\overline \bu};\overline{\partial}Y^n_{\overline \bu},\eta_p^{n-1}) - \widehat d_{\M_L,h}(Y^n-Y^{n-1},Y^n;\overline{\partial}Y^n,\eta_p^{n-1})\\
& + \widehat d_{\widehat \D,h}(P^{n-1}_{\overline \bu};P^{n-1}_{\overline \bu},\eta_y^{n}) - \widehat d_{\widehat \D,h}(P^{n-1};P^{n-1},\eta_y^{n}) \\
& + \widehat d_{\M_L,h}(P^{n-1}-P^{n},P^{n-1};\overline{\partial}P^n,\eta_y^{n}) - \widehat d_{\M_L,h}(P^{n-1}_{\overline \bu}-P^{n}_{\overline \bu},P_{\overline \bu}^{n-1};\overline{\partial}P^n_{\overline \bu},\eta_y^{n}) - \|\eta_y^n\|^2_{L^{2}}.
\end{align*}
Now, using some elementary calculations and the estimates derived for the stabilization differences above, we get
\begin{align*}
& (B(\overline u^n - U^n),\eta_p^{n-1}) \leq \frac{1}{k}(\eta_y^n, \eta_p^n) - \frac{1}{k}(\eta_y^{n-1}, \eta_p^{n-1})\\
& \qquad + C\,(h\|\nabla \eta_y^n\|_{L^{2}}\|\nabla \eta_p^{n-1}\|_{L^{2}} + h^2k^{-1}(\|\nabla \eta_y^n\|_{L^{2}} + \|\nabla \eta_y^{n-1}\|_{L^{2}})(\|\nabla \eta_p^{n-1}\|_{L^{2}} + \|\nabla \eta_p^{n}\|_{L^{2}})).
\end{align*}

Thus, using the latter estimate together with the result of Theorem \ref{theorem:error_estimates_AFC_u}, we get
\begin{align}\label{eq_dummy_3}
k\sum_{j=1}^{\NT}\|U^j - \overline u^j\|_{L^{2}}^2 \leq (\eta_y^{\NT}, \eta_p^{\NT}) - (\eta_y^0, \eta_p^0) + Ch(k + h)\sum_{j=1}^{\NT}(\|\nabla \eta_y^j\|_{L^{2}}^2 + \|\nabla \eta_p^{j-1}\|_{L^{2}}^2) + C_{M,\mu}(k^2 + h^4).
\end{align}
To estimate the right hand side, we need to obtain estimates for the norms inside the sum. Recall the equation \eqref{eq_dummy_1}. We will test $\chi = \eta_y^n,$ but now, we use different estimation to the stabilization terms. In particular, to estimate the stabilization term $\widehat d_{\D,h},$ we use Lemma \ref{lemma:stability_widehat_D_difference} together with the inverse inequality \eqref{eq:inverse_estimate} and we get
\begin{align*}
\vert \widehat d_{\D,h}(Y^n_{\overline \bu};Y^n_{\overline \bu},\chi)  - \widehat d_{\D,h}(Y^n;Y^n,\chi) \vert & \leq Ch\|\nabla \eta_y^n\|_{L^{2}}\|\nabla \chi\|_{L^{2}} \leq C\|\nabla \eta_y^n\|_{L^{2}}\|\chi\|_{L^{2}}.
\end{align*}
For the stabilization term due to mass lumping, we use Lemmas \ref{lemma:stability_widehat_M_difference}, \ref{lemma:regularity}, to obtain, that
\begin{align*}
& \vert \widehat d_{\M_L,h}(Y^n-Y^{n-1},Y^n;\overline{\partial}Y^n,\chi) -  d_{\M_L,h}(Y^n_{\overline \bu}-Y^{n-1}_{\overline \bu},Y^n_{\overline \bu};\overline{\partial}Y^n_{\overline \bu},\chi)\vert  \\
& \qquad \leq C\,\frac{h^2}{k}\left( \| \nabla (Y^n - Y^n_{\overline \bu})\|_{L^{2}} + \| \nabla (Y^{n-1} - Y^{n-1}_{\overline \bu})\|_{L^{2}}\right)\|\nabla \chi\|_{L^{2}}\\
& \qquad \leq C\,\frac{h^2}{k}\left(\|\nabla (Y^n - \overline y^n)\|_{L^{2}} + \|\nabla(Y^n_{\overline \bu} - \overline y^n)\|_{L^{2}} + \|\nabla (Y^{n-1} - \overline y^{n-1})\|_{L^{2}} + \|\nabla(Y^{n-1}_{\overline \bu} - \overline y^{n-1})\|_{L^{2}}\right)\|\nabla \chi\|_{L^{2}}\\
& \qquad \leq C_M\,\frac{h^2}{k}\left( \|\nabla (\theta_y^n + \rho_y^n)\|_{L^{2}} + \|\nabla (\theta_y^{n-1} + \rho_y^{n-1})\|_{L^{2}} +  k + h + k^{-1/2}h^2 \right)\|\nabla \chi\|_{L^{2}}\\
& \qquad \leq C_M\,\frac{h^2}{k}\left( \|\nabla \theta_y^n\|_{L^{2}} + \|\nabla \theta_y^{n-1}\|_{L^{2}} +  k + h + k^{-1/2}h^2 \right)\|\nabla \chi\|_{L^{2}},
\end{align*}
where $\theta_y^n,\,\rho_y^n,\,n=1,\ldots,\NT,$ as defined in Theorem \ref{theorem:error_estimates_AFC}. Further, we have used Lemma \ref{lemma:truncation_error_fully} with $\bpsi = \overline\bu,$ see also Remark \ref{remark:auxilliary_problem}. 

We notice in view of Theorem \ref{theorem:error_estimates_AFC}, the inverse inequality \eqref{eq:inverse_estimate} and the condition $k = \mathcal{O}(h)$, that
\begin{align*}
\frac{h^4}{k^2}\left(\|\nabla \theta_y^n\|_{L^{2}}^2 + \|\nabla \theta_y^{n-1}\|_{L^{2}}^2\right) & \leq C \left(\|\theta_y^n\|_{L^{2}}^2 + \| \theta_y^{n-1}\|_{L^{2}}^2\right) \leq C_{M,\mu}\,\left(k\sum_{j=1}^{\NT}\|U^j - \overline u^j\|_{L^{2}}^2 + k^2 + h^4 \right).
\end{align*}
Combining all these estimates and testing $\chi = \eta^n_y$, we get for a small $k$ and $k^2 + h^2 + k^{-1}h^4 = \mathcal{O}(k^2),$
\begin{align*}
& \|\eta^n_y\|_{L^{2}}^2 + k\tribar \eta^n_y \tribar^2_h \leq \|\eta^{n-1}_y\|_{L^{2}}^2 + C_{M,\mu}\,k\,\left(k\sum_{j=1}^{\NT}\|U^j - \overline u^j\|_{L^{2}}^2 + k^2 + h^4 \right).
\end{align*}
Summing over $1$ to $\NT,$
\begin{align*}
& \|\eta^{\NT}_y\|_{L^{2}}^2 + k\sum_{j=1}^{\NT}\tribar \eta^j_y \tribar^2_h \leq \|\eta^{0}_y\|_{L^{2}}^2 + C_{M,\mu}\left(k\sum_{j=1}^{\NT}\|U^j - \overline u^j\|_{L^{2}}^2 + k^2 + h^4 \right).
\end{align*}
Next, testing $\chi = \eta_p^{n-1}$ into \eqref{eq_dummy_2} and using similar arguments, we obtain that
\begin{align*}
& \|\eta^{0}_p\|_{L^{2}}^2 + k\sum_{j=1}^{\NT}\tribar \eta^{j-1}_p \tribar^2_h \leq \|\eta^{\NT}_p\|_{L^{2}}^2 + C_{M,\mu}\left(k\sum_{j=1}^{\NT}\|\eta_y^j\|_{L^{2}}^2 + k\sum_{j=1}^{\NT}\|U^j - \overline u^j\|_{L^{2}}^2 + k^2 + h^4 \right).
\end{align*}
Combining the last two estimates with \eqref{eq_dummy_3} and using the fact that $\eta^{0}_y,\, \eta^{\NT}_p = 0,$ we obtain
\begin{align*}
k\sum_{j=1}^{\NT}\|U^j - \overline u^j\|_{L^{2}}^2 \leq C_{M,\,\mu}\left(k^{-1}h(k+h)\,k\,\sum_{j=1}^{\NT}\|U^j - \overline u^j\|_{L^{2}}^2 + k^2 + h^4 \right),
\end{align*}
from which for small $k,\,h,$ with $k=\mathcal{O}(h),$ we conclude the desired result.
\end{proof}

\begin{theorem}\label{theorem:estimate_global}
Let $(\overline y,\overline p, \overline u)$ be the unique solution of \eqref{optimal_semi_weak_y}--\eqref{optimal_semi_weak_u2}, which is sufficiently regular by Lemma \ref{lemma:regularity} and let $(Y^n,P^{n-1},U^n)$ be a solution of the fully--discrete scheme \eqref{optimal_AFC_fully_weak_y}--\eqref{optimal_AFC_fully_weak_u2}. Assume the correction factors computed as in Definition \ref{definition:correction_factors_AFC}. Then for $k$, $h$ sufficiently small and $k=\mathcal{O}(h),$ there exists $C_{M,\mu}>0$, independent of $k,\,h,$ depending on $M,\,\mu^\nu,\,\nu<0,$ such that at time level $t=t^n$ we have for the state and co-state function
\begin{equation}\label{error_estimate_afc_state_global}
\begin{aligned}
& \|Y^n - \overline y^n\|_{L^{2}} + \|P^{n-1} - \overline p^{n-1}\|_{L^{2}} + \|U^n - \overline u^n\|_{L^{2}} \\
& \;\;\quad\qquad\qquad + \left(k \sum_{n=1}^{\NT}\left(\tribar Y^n - \overline y^n \tribar^2 + \tribar P^{n-1} - \overline p^{n-1} \tribar^2 \right)\right)^{1/2} \leq C_{M,\mu}(k + h^2),\\
& \tribar Y^n - \overline y^n \tribar + \tribar P^{n-1} - \overline p^{n-1} \tribar  \leq C_{M,\mu}(k + h).
\end{aligned}
\end{equation}
\end{theorem}
\begin{proof}
The proof of the estimates is a direct consequence of the estimates in Theorems \ref{theorem:error_estimates_AFC}, \ref{theorem:error_estimates_AFC_u}, Lemma \ref{lemma:truncation_error_fully2_2}, the estimates \eqref{ritz_projection_est2_2D}  and the error splittings $Y^n-\overline y^n=(Y^n - R_h\overline y^n)+(R_h\overline y^n-\overline y^n)$ and $P^n-\overline p^n=(P^n - R_h\overline p^n)+(R_h\overline p^n-\overline p^n).$
 \end{proof}

\section{An iterative scheme for approximating the nonlinear problem}\label{section:fixed_point}
In this section, we define an iterative scheme using the fixed-point method. By selecting appropriate correction factors based in subsection \ref{subsection:corr_factors}, we demonstrate the existence and uniqueness of the nonlinear AFC scheme described in \eqref{optimal_AFC_fully_weak_y}--\eqref{optimal_AFC_fully_weak_u2}.

\subsection{An auxiliary problem}
We define the following auxiliary fully--discrete scheme. For given $\bpsi\in\Uadk,$ we seek at time level $t=t^n,$ for known $\widetilde Y^{n-1} := Y^{n-1}(\bpsi)\in\Sh,$ the function $\widetilde Y^{n} :=\widetilde Y^n(\bpsi)\in \Sh$, for $n=1,\dots,\NT$, solution of \eqref{optimal_AFC_fully_weak_y_psi}.  After the computation of $\widetilde Y^n,\,n=1,\ldots,\NT,$ we seek at time level $t=t^{n-1},$ for known $\widetilde P^n := \widetilde P^n(\bpsi)\in \Sh,$ the function $\widetilde P^{n-1} := \widetilde P^{n-1}(\bpsi)\in \Sh,$ solution of \eqref{optimal_AFC_fully_weak_p_psi}. Then, we seek $\widetilde U^n:=\widetilde U^n(\bpsi)$ by
\begin{align}
(\lambda\, \widetilde U^n +  B^*\widetilde P^{n-1}, v - \widetilde U^n) & \geq 0,\;\;\;\;\forall\,v\in L^{2},\;\text{with}\;\mathsf{u_a}\leq v \leq \mathsf{u_b},\;\text{a.e. in}\;\Omega.\label{optimal_AFC_fully_weak_u_psi}
\end{align}
Similarly to \eqref{optimal_weak_u2}, one can express \eqref{optimal_AFC_fully_weak_u_psi} by using \eqref{operator_Pi}, as
\begin{align}\label{optimal_AFC_fully_weak_u2_psi}
\widetilde U^n  := \Pi_{[\mathsf{u_a},\mathsf{u_b}]}(-\lambda^{-1} B^*\widetilde P^{n-1}),\;n = \NT,\ldots,1.
\end{align}

\begin{remark}\label{remark:existence}
For given $\bpsi\in\Uadk,$ the auxiliary fully--discrete scheme \eqref{optimal_AFC_fully_weak_y_psi},\eqref{optimal_AFC_fully_weak_p_psi},\eqref{optimal_AFC_fully_weak_u2_psi} is decoupled with respect to the control variable, but both \eqref{optimal_AFC_fully_weak_y_psi}, \eqref{optimal_AFC_fully_weak_p_psi} are nonlinear due to the presence of the stabilization terms. Moreover, the solution of the AFC scheme \eqref{optimal_AFC_fully_weak_y}--\eqref{optimal_AFC_fully_weak_u2} for which we are interested to
prove its existence and uniqueness, is the unique fixed point of the mapping $\Uadk\ni\bpsi \mapsto \widetilde \bU(\bpsi)\in\Uadk.$
\end{remark}

\subsection{A fixed point scheme}

To define the fixed point scheme, we introduce the following stabilization terms for every $w,\,s\in\CC,$ in accordance to \eqref{stab_term_D_widetilde}, \eqref{stab_term_D_widetilde2},
\begin{align}
\doublewidetilde d_{\D,h}(w,s;v,z) := \sum_{i<j}\,d_{ij}\,\mathfrak{a}_{ij}(w,s)(v_i - v_j)(z_i - z_j),\quad\forall v,\,z\in{\CC},\label{stab_term_D_widetilde3}\\
\doublewidetilde d_{\widehat\D,h}(w,s;v,z) := \sum_{i<j}\,d_{ij}\,\widehat{\mathfrak{a}}_{ij}(w,s)(v_i - v_j)(z_i - z_j),\quad\forall v,\,z\in{\CC}.\label{stab_term_D_widetilde4}
\end{align}
The difference between \eqref{stab_term_D_widetilde} and \eqref{stab_term_D_widetilde3}, (respectively \eqref{stab_term_D_widetilde2} and \eqref{stab_term_D_widetilde4}) are the fact that now the correction terms depends in two variables. This modification is necessary so that the linearized scheme satisfies the local and global DMP. The proof follows the lines of Theorems \ref{theorem:dmp_fully}, \ref{theorem:global_dmp_fully}.

For $\Xh : = \Sh \times \Sh,$ we consider the iteration operator  $\Gn : \Xh \to \Xh$, $(v^n(\bpsi),w^{n-1}(\bpsi)) $ $\mapsto$ $ (\Gn_1 v^n(\bpsi),\Gno_2 w^{n-1}(\bpsi)),$ for $\bpsi\in\Uadk,$ defined as follows. For a given $(v^n,\bpsi)\in \Sh\times \Uadk,$ starting with $v^n = \widetilde Y^{n-1},$ for every time step $t=t^n,\,n=1,\ldots,\NT,$ we define $\Gn_1v^n := \Gn_1v^n(\bpsi)$ as
\begin{equation}\label{existence_linear_y}
\begin{aligned}
& (\Gn_1 v^n - \widetilde Y^{n-1}, \chi)  +  k\,\alpha_h(\Gn_1 v^n,\chi)  - k\,d_{\D,h}(\Gn_1 v^n, \chi) +  d_{\M_L,h}(\Gn_1 v^n - \widetilde Y^{n-1}, \chi)\\
& \qquad\quad - \widetilde d_{\M_L,h}(v^n-\widetilde Y^{n-1},\widetilde Y^{n-1}; v^n - \widetilde Y^{n-1} , \chi)  = - k\,\doublewidetilde d_{\D,h}(v^n,Y^{n-1};v^n, \chi) + k\,( B\psi^n + G^n, \chi),
\end{aligned}
\end{equation}
for the state variable for all $\chi\in\Sh.$ After having computed $\Gn_1 v^n$ that depend on $\psi^n,$ we define $\Gno_2w^{n-1} = \Gno_2w^{n-1}(\bpsi),$ solution of the linearized scheme for the co-state variable, starting with $w^n = \widetilde P^n,$  for every time step $t=t^{n-1},\,n=\NT,\ldots,1,$ as
\begin{equation}\label{existence_linear_p}
\begin{aligned}
  (\Gno_2 w^{n-1} - \widetilde P^n, \chi) & +  k\,\alpha_h(\Gno_2 w^{n-1},\chi) + k\,d_{\widehat\D,h}(\Gno_2 w^{n-1}, \chi) \\
&  +  d_{\M_L,h}(\Gno_2 w^{n-1} - \widetilde P^n, \chi) - \widetilde d_{\M_L,h}(w^{n-1} - \widetilde P^n,\widetilde P^{n}; w^{n-1} - \widetilde P^n , \chi)\\
&  = k\,\doublewidetilde d_{\widehat \D,h}(w^{n-1},P^{n};w^{n-1}, \chi) + k\,(\Gn_1 v^{n} -  y_d^n, \chi),
\end{aligned}
\end{equation}
for all $\chi\in\Sh.$ Furthermore, we define the iteration operator $\F = (\F^1,\ldots,\F^{\NT})^T\,:\,\Uadk\to\Uadk,$ as
\begin{align}
(\lambda\, \F^n \psi^n +  B^*\Gno_2 w^{n-1}, \vartheta - \F^n \psi^n) & \geq 0\;\;\;\;\forall\,\vartheta\in L^{2},\;\text{with}\;\mathsf{u_a}\leq \vartheta \leq \mathsf{u_b},\;\text{a.e. in}\;\Omega.\label{existence_linear_u}
\end{align}
Similarly to \eqref{optimal_weak_u2}, one can express \eqref{existence_linear_u} by using \eqref{operator_Pi}, as
\begin{align}\label{existence_linear_u2}
\F^n \psi^n := \Pi_{[\mathsf{u_a},\mathsf{u_b}]}(-\lambda^{-1} B^* \Gno_2 w^{n-1}),\;\;n = \NT,\ldots,1.
\end{align}
Thus, through $\Gn$, we obtain a sequence of elements $v^n_{j+1} = \Gn_1v^n_j$, for $n=1,\ldots,\NT,\,j\geq 0$ and $w^{n-1}_{j+1} = \Gno_2w^{n-1}_j$, for $n=\NT,\ldots,1,\,j\geq 0.$ In addition, $\bpsi_{j+1} = \F\bpsi_j,\,j\geq 0.$

Next, we write the matrix formulation of the fully--discrete scheme \eqref{existence_linear_y}--\eqref{existence_linear_u2}. To do this, we introduce the following notation. Let $\widetilde \al^n=(\widetilde \alpha_1^n,\dots,\widetilde \alpha_{\N}^n)^T,\,\widetilde \be^{n-1} =(\widetilde \beta^{n-1}_1,\dots,\widetilde \beta_{\N}^{n-1})^T,$ the coefficients, with respect to the basis of $\Sh$ of $\Gn_1 v^n,\,\Gno_2 w^{n-1}\in\Sh,$ respectively. Similarly, let $\widehat \al^n=(\widehat \alpha_1^n,\dots,\widehat \alpha_{\N}^n)^T,\,\widehat \be^{n-1}=(\widehat \beta^{n-1}_1,\dots,\widehat \beta_{\N}^{n-1})^T$ the coefficients, with respect to the basis of $\Sh$ of $v^n,\,w^{n-1}\in\Sh,$ respectively. Then \eqref{existence_linear_y}--\eqref{existence_linear_u2} can be written as
\begin{align}
(\M_L +  k\,(\mu\S + \Q + \D + \sigma\,\M)) \widetilde \al^n & = \M_L\al^{n-1} + k\,\bfr(B\psi^n + G^n) +  k\,\overline{\mathsf{f}}^1(\widehat\al^n,\al^{n-1}) + \overline{\mathsf{g}}^1(\widehat\al^n,\al^{n-1}),\label{linear_fully_AFC_matrix_y}\\
(\M_L + k\,(\mu\S - \Q - \widehat\D + \sigma\,\M)) \widetilde \be^{n-1} & = \M_L\be^{n} + k\,\M\widetilde \al^{n} - k\,\bfr( y^n_d) +  k\,\overline{\mathsf{f}}^2(\widehat\be^{n-1},\be^{n}) + \overline{\mathsf{g}}^2(\widehat\be^{n-1},\be^{n}),\label{linear_fully_AFC_matrix_p}
\end{align}
with $\al^0 = \widehat v_1,\,\be^{\NT} = \mathbf{0},$ where $\widehat v_1,\,\mathbf{0}$ are the coefficient vectors of $y_h^0,\,p_h^{\NT}\in\Sh.$ \par

\begin{definition}\label{definition:correction_factors_fully_AFC_linear}
The correction terms $\overline{f}^\iota,\,\iota=1,2,$ with correction factors $\mathfrak{a}_{ij} = \mathfrak{a}_{ij}(\widehat\al^n,\al^{n-1})$ and $\widehat{\mathfrak{a}}_{ij} = \widehat{\mathfrak{a}}_{ij}(\widehat\be^{n-1},\be^{n}),\,i,j=1,\ldots,\N,$ for $\iota=1,2,$ respectively, are computed in view of Definition \ref{definition:corr_factors_gen} and Algorithm \ref{algorithm-1}, by
\begin{equation}\label{correction_term_linear_f}
\begin{aligned}
(\overline{\mathsf{f}}^{1})_i := (\overline{\mathsf{f}}^{1}(\widehat{\al}^n,\al^{n-1}))_i = \sum_{j\neq i}\mathfrak{a}_{ij}\mathsf{f}_{ij},\;\text{with}\;P^{\pm} & :=  P^{\pm}(\widehat{\al}^n),\,Q^{\pm} := Q^{\pm}(\al^{n-1}),\,q_{i} = \gamma_i\sum_{j\in\Zh^i}d_{ij},\\
(\overline{\mathsf{f}}^{2})_i := (\overline{\mathsf{f}}^{2}(\widehat{\be}^{n-1},\be^n))_i = \sum_{j\neq i}\widehat{\mathfrak{a}}_{ij}\widehat{\mathsf{f}}_{ij},\;\text{with}\;P^{\pm} & :=  P^{\pm}(\widehat{\be}^{n-1}),\,Q^{\pm} := Q^{\pm}(\be^{n}),\,q_{i} = \gamma_i\sum_{j\in\Zh^i}\widehat d_{ij},
\end{aligned}
\end{equation}
for $i=1,\ldots,\N.$  Similarly, the correction terms $\overline{g}^\iota,\,\iota=1,2,$ with correction factors $\mathrm{a}_{ij} = \mathrm{a}_{ij}(\widehat\al^{n}-\al^{n-1},\al^{n-1})$, $\widehat{\mathrm{a}}_{ij} = \widehat{\mathrm{a}}_{ij}(\widehat\be^{n-1}-\be^{n},\be^{n})$, $i,j=1,\ldots,\N,$ respectively, are computed in view of Definition \ref{definition:corr_factors_gen} and Algorithm \ref{algorithm-1}, by
\begin{equation}\label{correction_term_linear_g1}
\begin{aligned}
(\overline{\mathsf{g}}^{1})_i & := (\overline{\mathsf{g}}^{1}(\widehat{\al}^n,\al^{n-1}))_i = \sum_{j\neq i}\mathrm{a}_{ij}\mathsf{g}_{ij},\;\text{with}\;P^{\pm} :=  P^{\pm}(\widehat{\al}^{n}-\al^{n-1}),\,Q^{\pm} := Q^{\pm}(\al^{n-1}),\\
(\overline{\mathsf{g}}^{2})_i & := (\overline{\mathsf{g}}^{2}(\widehat{\be}^{n-1},\be^{n}))_i = \sum_{j\neq i}\widehat{\mathrm{a}}_{ij}\widehat{\mathsf{g}}_{ij},\;\text{with}\;P^{\pm} :=  P^{\pm}(\widehat{\be}^{n-1}-\be^{n}),\,Q^{\pm} := Q^{\pm}(\be^{n}),
\end{aligned}
\end{equation}
with $q_{i} = \gamma_i\sum_{j\in\Zh^i}m_{ij}$ for $i=1,\ldots,\N.$ The amount of mass transported by the raw antidiffusive flux $\widehat{\mathsf{f}}_{ij}$ is given by
\begin{align}
\mathsf{f}_{ij} & := \mathsf{f}_{ij}(\widehat{\al}^{n}) = \left(\widehat{\alpha}_j^{n} - \widehat{\alpha}_i^{n}\right)d_{ij},\;\;\;\qquad\;\;\;\forall \,j\neq i,\label{flux_D_linear}\\
\widehat{\mathsf{f}}_{ij} & := \widehat{\mathsf{f}}_{ij}(\widehat{\be}^{n-1}) = \left(\widehat{\beta}_i^{n-1} - \widehat{\beta}_j^{n-1}\right)\widehat d_{ij},\;\;\forall \,j\neq i,\label{flux_D_widehat_linear}
\end{align}
and the amount of mass due to the replacement of the mass matrix $\M$ by its mass lumped $\M_L,$ by its mass lumped $\M_L,$ the $\widehat{\mathsf{g}}_{ij},$ are given by
\begin{align}
 \mathsf{g}_{ij} & := \mathsf{g}_{ij}(\widehat{\al}^{n},\al^{n-1}) = m_{ij} ( (\widehat{\alpha}_i^{n} - \widehat{\alpha}_j^{n}) - (\alpha^{n-1}_i - \alpha^{n-1}_j)),\;\;\;\;\forall \,j\neq i,\label{flux_M_linear1}\\
\widehat{\mathsf{g}}_{ij} & := \widehat{\mathsf{g}}_{ij}(\widehat{\be}^{n-1},\be^{n}) = m_{ij} ( (\widehat{\beta}_i^{n-1} - \widehat{\beta}_j^{n-1}) - (\beta^{n}_i - \beta^{n}_j)),\,\quad\;\forall \,j\neq i,\label{flux_M_linear2}
\end{align}
\end{definition}

\noindent
The above linearized scheme can be summarized in the following algorithm.
\begin{algorithm}
\begin{algorithmic}[1]
\State Given a tolerance for the while loop, defined as $\mathrm{TOL}.$
\State Set the initial values $\al^0,\,\be^{\NT}$ as the coefficient vectors of $y_h^0,\,p_h^{\NT},$ respectively.
\State Initialize $\widehat \ga^n,\,\widetilde \ga^n = \mathbf{0} \in\R^{\N},\,n=1,\ldots,\NT,$ the zero vector, also $\widehat\al^0 = \widehat\al^1 := \al^0$ and $\widehat\be^{\NT} = \widehat\be^{\NT-1} := \be^{\NT}.$
\While{$1$}
\For{$n=1,\ldots,\NT$}
\While{$1$}\\ 
\hspace{1.3cm} Given $\widehat \al^{n-1}(\widetilde \ga^{n-1}),\,\widehat\al^n(\widehat \ga^n)\in\R^{\N},$ compute the correction factors $\mathfrak{a}_{ij},\,\mathrm{a}_{ij}$ as in Definition \ref{definition:correction_factors_fully_AFC_linear}.\\
\hspace{1.3cm} Given $\widehat \ga^n\in\R^{\N},$ solve \eqref{linear_fully_AFC_matrix_y} to obtain $\widetilde\al^n(\widehat \ga^n)\in\R^{\N}$ for which $\widetilde\al^n(\widehat \ga^n) \approx \al^n(\ga^n).$\\
\hspace{1.3cm} If $\frac{\|\widetilde\al^n(\widehat \ga^n) - \widehat\al^n(\widehat \ga^n)\|_{\max}}{\|\widetilde\al(\widehat\ga^n)\|_{\max}}  \leq \mathrm{TOL}$ set $\widehat\al^n(\widehat \ga^n) = \widetilde\al^n(\widehat \ga^n)$ and break the inner while.
\EndWhile
\EndFor
\For{$n=\NT,\ldots,1$}
\While{$1$}\\ 
\hspace{1.3cm} Given $\widehat \be^{n}(\widehat \ga^n),\,\widehat\be^{n-1}(\widehat \ga^{n-1})\in\R^{\N},$ compute the correction factors $\widehat{\mathfrak{a}}_{ij},\,\widehat{\mathrm{a}}_{ij}$ as in Definition \ref{definition:correction_factors_fully_AFC_linear}.\\
\hspace{1.3cm} Given $\widetilde \al^{n}(\widehat \ga^{n})\in\R^{\N},$ solve \eqref{linear_fully_AFC_matrix_p} to obtain $\widetilde\be^{n-1}(\widehat \ga^{n})\in\R^{\N}$ for which $\widetilde\be^{n-1}(\widehat \ga^{n}) \approx \be^{n-1}( \ga^{n}).$\\
\hspace{1.3cm} Compute the control by from \eqref{existence_linear_u2} and $\widetilde\ga^{n}$ be the coefficients of its $L^{2}-$projection \eqref{L2_modf_projection}.\\
\hspace{1.3cm} If $\frac{\|\widetilde\be^n(\widehat \ga^n) - \widehat\be^n(\widehat \ga^n)\|_{\max}}{\|\widetilde\be(\widehat\ga^n)\|_{\max}}  \leq \mathrm{TOL}$ set $\widehat\be^{n-1}(\widehat \ga^n) = \widetilde\be^{n-1}(\widehat \ga^n)$ and break the inner while.
\EndWhile
\EndFor
\State Compute $A := \max_{0\leq n\leq \NT}\left\lbrace \frac{\|\widetilde\al^n(\widehat \ga^n) - \widetilde\al^n(\widetilde \ga^n)\|_{\max}}{\|\widetilde\al(\widetilde\ga^n)\|_{\max}}\right\rbrace$
\State Compute $B := \max_{0\leq n\leq \NT}\left\lbrace \frac{\|\widetilde\be^n(\widehat \ga^n) - \widetilde\be^n(\widetilde \ga^n)\|_{\max}}{\|\widetilde\be(\widetilde\ga^n)\|_{\max}}\right\rbrace$
\State Compute $\Gamma := \max_{0\leq n\leq \NT}\left\lbrace \frac{\|\widehat \ga^n - \widetilde \ga^n\|_{\max}}{\|\widetilde\ga^n\|_{\max}}\right\rbrace$
\If{$\max\{A,B,\Gamma\} \leq \mathrm{TOL}$}
        \State Set $\widehat \ga^n := \widetilde \ga^n,\,\widehat \al^n(\widetilde \ga^n) := \widetilde \al^n(\widetilde \ga^n)$ for $n=1,\ldots,\NT$ and $\widehat \be^{n-1}(\widetilde \ga^n) := \widetilde \be^{n-1}(\widetilde \ga^n)$ for $n = \NT, \ldots,1.$       
     \State  Break the while loop
\Else
        \State Continue.
\EndIf

\EndWhile
\end{algorithmic}
\caption{Computation of the fixed point scheme \eqref{existence_linear_y}--\eqref{existence_linear_u2}}\label{algorithm-2}
\end{algorithm}

\subsection{Existence and uniqueness of the AFC scheme}\label{subsection:existence}

The proof of existence and uniqueness of \eqref{optimal_AFC_fully_weak_y_psi}--\eqref{optimal_AFC_fully_weak_u2_psi} will be proceed in two steps. First, we prove that for given $\bpsi\in\Uadk,$ the (decoupled) auxiliary problems \eqref{optimal_AFC_fully_weak_y_psi}, \eqref{optimal_AFC_fully_weak_p_psi} have a unique solution. Then, we note that \eqref{optimal_AFC_fully_weak_y}--\eqref{optimal_AFC_fully_weak_u2}  can be viewed as $\bpsi = \widetilde \bU(\bpsi),$  defined through \eqref{optimal_AFC_fully_weak_y_psi},\eqref{optimal_AFC_fully_weak_p_psi},\eqref{optimal_AFC_fully_weak_u2_psi}, thus, our aim is to prove that the latter mapping is contraction in $\Uadk.$ 

Before we prove the above result, we derive a similar estimate to Lemma \ref{lemma:stability_widehat_D_difference} for the stabilization terms defined in \eqref{stab_term_D_widetilde3}, \eqref{stab_term_D_widetilde4}.

\begin{lemma}\label{lemma:stability_widehat_D_difference2}
Let the bilinear forms $\doublewidetilde d_{\D,h},\,\doublewidetilde d_{\widehat\D,h}$ defined in \eqref{stab_term_D_widetilde3}, \eqref{stab_term_D_widetilde4}, respectively. If their correction factors are satisfying the Definition \ref{definition:correction_factors_fully_AFC_linear}, there exists for $w,\,s,\,v,\,\chi\in \Sh,$ a positive constant $C$, independent of $h,$ such that,
\begin{align}
\vert\doublewidetilde{d}_{\D,h}(v,s;v, \chi) - \doublewidetilde{d}_{\D,h}(w,s;w,\chi)\vert + \vert\doublewidetilde{d}_{\widehat\D,h}(v,s;v, \chi) - \doublewidetilde{d}_{\widehat\D,h}(w,s;w,\chi)\vert & \leq Ch\|\nabla (v - w)\|_{L^{2}}\|\nabla \chi\|_{L^{2}}.\label{stability_stab_term_widetilde_D34}
\end{align}
\end{lemma}
\begin{proof}
The proof follows the lines of the proof of Lemma \ref{lemma:stability_widehat_D_difference}. In particular, in that lemma the main idea was the local estimate \eqref{local_estimate} and the shape regularity of the triangulation. Therefore, we need to ensure that for the correction factors defined in \eqref{stab_term_D_widetilde3} and Definition \ref{definition:correction_factors_fully_AFC_linear}, we have for fixed $\psi\in\Sh,$ the following local estimate
\begin{align*}
\vert d_{ij} \mathfrak{a}_{ij}(\chi,\psi)(\chi_i - \chi_j) & - d_{ij} \mathfrak{a}_{ij}(\widetilde\chi,\psi)(\widetilde\chi_i - \widetilde\chi_j) \vert \leq C\,h\sum_{\ell\in \Zh(\omega_i)}|\chi_\ell - \widetilde\chi_\ell|,\;\;\forall\,\chi,\,\widetilde\chi\in\Sh.
\end{align*}
Such an estimate can be verified by same arguments as in proof of Lemma \ref{lemma:stability_widehat_M_difference}. A similar inequality holds also for the correction factors defined \eqref{stab_term_D_widetilde4} and Definition \ref{definition:correction_factors_fully_AFC_linear}.
\end{proof}

\begin{theorem}\label{theorem:existence_uniqueness}
Assume the correction factors in \eqref{optimal_AFC_fully_weak_y_psi}--\eqref{optimal_AFC_fully_weak_p_psi} are computed as in Definition \ref{definition:correction_factors_AFC}. For small $k,\,h$ with $k = \mathcal{O}(h)$, the following statements hold:
\begin{enumerate}
    \item For a given $\bpsi \in \Uadk$, the nonlinear fully discrete scheme \eqref{optimal_AFC_fully_weak_y_psi}--\eqref{optimal_AFC_fully_weak_p_psi} has a unique solution.
    \item The mapping $\widetilde \bU(\bpsi)$, defined through \eqref{optimal_AFC_fully_weak_y_psi},\eqref{optimal_AFC_fully_weak_p_psi},\eqref{optimal_AFC_fully_weak_u2_psi}, is contractive in $\Uadk$.
\end{enumerate}
As a result, the nonlinear fully discrete scheme \eqref{optimal_AFC_fully_weak_y}--\eqref{optimal_AFC_fully_weak_u2} has a unique solution for small $k,\,h$ with $k = \mathcal{O}(h)$.
\end{theorem}
\begin{proof}
Let a given $\bpsi\in\Uadk.$ We note that if the map $\Gn$ defined in \eqref{existence_linear_y}--\eqref{existence_linear_p} has a fixed point $(v^*, w^*)$, then $(\widetilde Y^n, \widetilde P^{n-1} ) := (v^*, w^*)$ is the solution of the
discrete scheme \eqref{optimal_AFC_fully_weak_y_psi}--\eqref{optimal_AFC_fully_weak_p_psi}. Thus, our aim is to prove that for the mapping $\Gn$, is contractive on $\Xh$. \par
Assume $v_1^n,\,v_2^n\in\Sh$ with $v_1^n\neq v_2^n$ and $\widehat v^n := v_1^n - v_2^n.$ Similarly, $w_1^{n-1},\,w_2^{n-1}\in\Sh,$ with $w_1^{n-1}\neq w_2^{n-1}$ and $\widehat w^{n-1}= w_1^{n-1}-w_2^{n-1}.$ Clearly, by \eqref{existence_linear_y}--\eqref{existence_linear_p}, we obtain $\Gn_1 v_1^n,\,\Gn_1v_2^n\in\Sh$ and $\Gno_2w_1^{n-1},\,\Gno_2w_2^{n-1}\in\Sh.$ Using the previous notation for the mapping, i.e., setting $\Gn_1\widehat v^n := \Gn_1v_1^n - \Gn_1v_2^n$ and $\Gno_2\widehat w^n= \Gno_2w_1^{n-1}-\Gno_2w_2^{n-1},$ we get, by \eqref{existence_linear_y},
\begin{align*}
(\Gn_1\widehat v^n, \chi) + k\,\alpha_h(\Gn_1\widehat v^n, \chi) - k\,d_{\D,h}(\Gn_1\widehat v^n, \chi) + d_{\M_L,h}(\Gn_1\widehat v^n, \chi) & = I_1(\chi) - I_2(\chi),
\end{align*}
where the functionals $I_1,\,I_2,$ are defined and can be estimated by the Lemmas \ref{lemma:stability_widehat_D_difference2}, \ref{lemma:stability_widehat_M_difference}, the inverse inequality \eqref{eq:inverse_estimate} and the $\|\nabla \zeta\|_{L^2} \leq \mu^{-1/2}\tribar \zeta\tribar_h,\,\zeta\in\Sh,$
\begin{align*}
\vert I_1(\chi) \vert & = k\,\vert \doublewidetilde d_{\D,h}(v_1^n,Y^{n-1};v_1^n,\chi) - \doublewidetilde d_{\D,h}( v_2^n,Y^{n-1}; v_2^n,\chi)\vert \leq C\,k\,\|\nabla \widehat v^n\|_{L^{2}}\|\chi\|_{L^{2}} \leq C\,k\,\mu^{-1/2}\tribar \widehat v^n\tribar_h\|\chi\|_{L^{2}},
\end{align*}
and
\begin{align*}
\vert I_2(\chi) \vert & = \vert\widetilde d_{\M_L,h}(v_1^n - Y^{n-1},Y^{n-1};v_1^n - Y^{n-1},\chi) - \widetilde d_{\M_L,h}(v_2^n - Y^{n-1},Y^{n-1};v_2^n - Y^{n-1},\chi)\vert\\
& \leq C\,h\,\mu^{-1/2}\tribar \widehat v^n\tribar_h\|\chi\|_{L^{2}}.
\end{align*}
Setting $\chi = \Gn_1\widehat v^n$ and performing standard arguments, we have for small $k,\,h,$
\begin{align*}
\|\Gn_1\widehat v^n\|_{L^{2}}^2  & +  k\,\tribar \Gn_1\widehat v^n \tribar^2_h \leq C_{M,\mu}^2(k+h)^2\,\tribar \widehat v^n\tribar^2_h,
\end{align*}
where the constant $C_{M,\mu}$ depends on $M,\,\mu^\nu,\,\nu<0,$ while is independent of $k,\,h.$
Note that, to obtain the last estimate, we have used the inverse inequality \eqref{eq:inverse_estimate} as well as the fact that $ - d_{\D,h}(\chi, \chi) \geq 0$ and $d_{\M_L,h}(\chi, \chi) \geq 0$ for all $\chi\in\Sh,$ see, e.g., their definition. 
Similarly, by  \eqref{existence_linear_p}, we also get
\begin{align*}
\|\Gno_2\widehat w^{n-1}\|_{L^{2}}^2   +  k\,\tribar  \Gno_2\widehat w^{n-1} \tribar^2_h & \leq C_{M,\mu}^2(k+h)^2\tribar\widehat w^{n-1}\tribar^2_h + C\,k^2\,\|\Gn_1\widehat v^n\|_{L^{2}}^2\\
& \leq C_{M,\mu}^2(k+h)^2\tribar\widehat w^{n-1}\tribar^2_h + C_{M,\mu}^2k^2(k+h)^2\tribar \widehat v^n\tribar^2_h,
\end{align*}
where the constant $C_{M,\mu}$ depends on $M,\,\mu^\nu,\,\nu<0,$ while is independent of $k,\,h.$
Combining the above  estimates, yields,
\begin{align*}
\tribar\Gn_1\widehat v^n\tribar^2_h + \tribar\Gno_2\widehat w^{n-1}\tribar^2_h \leq C_{M,\mu}^2\,k^{-1}(k+h)^2(\tribar\widehat v^n\tribar_h^2 + \tribar\widehat w^{n-1}\tribar_h^2),
\end{align*}
and then, for small $k,\,h$ with $k = \mathcal{O}(h),$ such that $C_{M,\mu}^2\,k^{-1}(k+h)^2<1,$ the sequence $v^n_{j+1}(\bpsi) \to \Gn_1v^n_j(\bpsi)\in\Sh,\,j\to \infty$ as well as $w^{n-1}_{j+1}(\bpsi) \to \Gno_2w^{n-1}_j(\bpsi),\,j\to\infty.$ This means, that there exists a unique solution $(\widetilde Y^n(\bpsi), \widetilde P^{n-1}(\bpsi))\in\Sh\times\Sh$ to \eqref{optimal_AFC_fully_weak_y_psi}--\eqref{optimal_AFC_fully_weak_p_psi}. \par

After ensuring the existence and uniqueness of \eqref{optimal_AFC_fully_weak_y_psi}--\eqref{optimal_AFC_fully_weak_p_psi} for a fixed $\bpsi\in\Uadk,$ we now aim to establish the existence and uniqueness of \eqref{optimal_AFC_fully_weak_y}--\eqref{optimal_AFC_fully_weak_u2}. To achieve this, and in view of Remark \ref{remark:existence}, we consider $\bpsi_1,\,\bpsi_2\in \Uadk,$ and our aim is to prove that the mapping defined through \eqref{optimal_AFC_fully_weak_y_psi},\eqref{optimal_AFC_fully_weak_p_psi},\eqref{optimal_AFC_fully_weak_u2_psi} satisfies that $\Uadk \ni \bpsi \mapsto \widetilde \bU(\bpsi)\in\Uadk$ is contractive in $\Uadk.$ The former can be easily obtained by \eqref{optimal_AFC_fully_weak_u2_psi}. \par
In view of \eqref{optimal_AFC_fully_weak_u2_psi}, we have the following estimate,
\begin{align*}
\|\widetilde U^n(\bpsi_1) - \widetilde U^n(\bpsi_2)\|_{L^{2}} \leq C_\lambda\|\widetilde P^{n-1}(\bpsi_1) - \widetilde P^{n-1}(\bpsi_2)\|_{L^{2}}.
\end{align*}
Further, by \eqref{optimal_AFC_fully_weak_p_psi}, we get for small $k,\,h,$ with $k=\mathcal{O}(h),$
\begin{align*}
\|\widetilde P^{n-1}(\bpsi_1) - \widetilde P^{n-1}(\bpsi_2)\|_{L^{2}}^2 & + k\,\sum_{j=n}^{\NT}\tribar \widetilde P^{j-1}(\bpsi_1) - \widetilde P^{j-1}(\bpsi_2) \tribar^2_h \leq C_{M,\mu}\,k\,\sum_{j=n}^{\NT}\|\widetilde Y^{j}(\bpsi_1) - \widetilde Y^{j}(\bpsi_2)\|_{L^{2}}^2.
\end{align*}
More specifically, to obtain the latter estimate we have used the following estimates for the stabilization term, using the Lemma \ref{lemma:stability_widehat_D_difference},
\begin{align*}
\vert \widehat d_{\widehat \D,h}(P^{n-1}(\bpsi_1);P^{n-1}(\bpsi_1),\chi) & - \widehat d_{\widehat \D,h}(P^{n-1}(\bpsi_2);P^{n-1}(\bpsi_2),\chi)\vert\\
 & \leq C\,h\,\mu^{-1}\tribar P^{n-1}(\bpsi_1) - P^{n-1}(\bpsi_2)\tribar_h \tribar \chi\tribar_h,
\end{align*}
and Lemma \ref{lemma:stability_widehat_M_difference}, the inverse inequality \eqref{eq:inverse_estimate} and the condition $k=\mathcal{O}(h),$
\begin{align*}
& \vert \widehat d_{\M_L,h}(P^{n-1}(\bpsi_1) - P^{n}(\bpsi_1),P^{n-1}(\bpsi_1);\overline\partial P^{n}(\bpsi_1),\chi)  - \widehat d_{\M_L,h}(P^{n-1}(\bpsi_2) - P^{n}(\bpsi_2),P^{n-1}(\bpsi_2);\overline\partial P^{n}(\bpsi_2),\chi)\vert \\
&\qquad \leq C\,k^{-1}\mu^{-1/2}\left( h^2\mu^{-1/2}\tribar P^{n-1}(\bpsi_1) - P^{n-1}(\bpsi_2) \tribar_h + k\| P^{n}(\bpsi_1) - P^{n}(\bpsi_2)\|_{L^{2}}\right) \tribar \chi\tribar_h.
\end{align*}
To derive an estimate for the right hand side, we use \eqref{optimal_AFC_fully_weak_y_psi} and similar estimates as before, to get
\begin{align*}
\|\widetilde Y^{n}(\bpsi_1) - \widetilde Y^{n}(\bpsi_2)\|_{L^{2}}^2  & + k\,\sum_{j=1}^{n}\tribar \widetilde Y^{j}(\bpsi_1) - \widetilde Y^{j}(\bpsi_2) \tribar^2_h \leq C_{M,\mu}\,k\,\sum_{j=1}^{n}\|\psi_1^j - \psi_2^j\|_{L^{2}}^2.
\end{align*}
Combining the above estimates and using the fact that $k\sum_{j=1}^{\NT}1 = T$, yields, for all $n = \NT,\ldots,1,$
\begin{align*}
\|\widetilde U^n(\bpsi_1) - \widetilde U^n(\bpsi_2)\|_{L^{2}}^2 \leq C_{M,\mu,\lambda}\,k\,\sum_{j=1}^{\NT}\|\psi_1^j - \psi_2^j\|_{L^{2}}^2.
\end{align*}
We define a discrete time-dependent norm as $\|\bv\|_{\ell^{2}(L^{2})} := \left(k\sum_{j=1}^{\NT}\|v^j\|_{L^{2}}^2 \right)^{1/2},$ where $v^n := v(\x,t^n),\,t^n = n\,k,\,n=0,\ldots,\NT.$ Hence, we obtain
\begin{align*}
\|\widetilde \bU(\bpsi_1) - \widetilde \bU(\bpsi_2)\|_{\ell^{2}(L^{2})} \leq C_{M,\mu,\lambda}\,k^{1/2}\,\|\bpsi_1 - \bpsi_2\|_{\ell^{2}(L^{2})}.
\end{align*}
For small $k$ such that $C_{M,\mu,\lambda}k^{1/2}<1,$ the sequence $\bpsi_{j+1} \to \widetilde U^n(\bpsi_j),\,j\to \infty.$ Therefore, the fully--discrete scheme \eqref{optimal_AFC_fully_weak_y}--\eqref{optimal_AFC_fully_weak_u2}, has a unique solution.
\end{proof}

\section{Numerical Experiments}\label{section:numerical_results}

\begin{figure}
\centering
\includegraphics[scale=1]{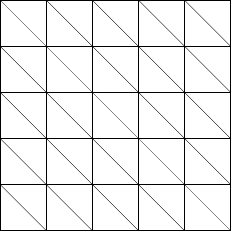}
\caption{A triangulation of a square domain.}\label{fig:triangulation}
\end{figure}

In this section we present several numerical experiments, validating our theoretical results. For this reason, we consider a uniform mesh $\Th$ of the unit square $\Omega = [0,1]^2.$ Each side of $\Omega$ is divided into $M$ intervals of length $h_0=1/M$ for $M\in\mathbb{N}$ and we define the triangulation $\Th$ by dividing each small square by its diagonal as in Fig. \ref{fig:triangulation}. This means that $\Th$ consists of $2M^2$ right-angle triangles with diameter $h = \sqrt{2}h_0.$ Note that $\Th$ satisfies Assumption \ref{mesh-assumption} and therefore, the corresponding stiffness matrix $\S$ has non-positive off-diagonal elements and positive diagonal elements. Also, $\gamma_i=1,\,i=1,\ldots,\N,$ in the Definition \ref{definition:correction_factors_fully_AFC_linear}, see, e.g., Remark \ref{remark:gamma}. In all numerical experiments, we set $\mathrm{TOL} = 10^{-6},$ see Algorithm \ref{algorithm-2}.

In the following subsections, we test our scheme on three optimal control problems governed by a parabolic convection--diffusion equation. Specifically, first we consider one test case in which both the state and co-state variables are smooth functions. Next, we consider that the state variable exhibits steep gradients while the co-state variable remains smooth, and another test case where the state variable develops an interior layer while the co-state variable exhibits steep gradients. Our aim is to compute the corresponding errors over a sequence of uniformly refined triangulations. In particular, we define the following errors,
\begin{equation}\label{computed_errors}
\begin{aligned}
\|(e_y^n, e_p^n)\|_{\max} & := \left( \max_{0\leq n\leq \NT}\|e_y^n\|_{L^{2}}^2 + \max_{0\leq n\leq \NT}\|e_p^n\|_{L^{2}}^2\right)^{1/2},\\
\tribar(e_y^n, e_p^n)\tribar_{\max} & := \left( \max_{0\leq n\leq \NT}\tribar e_y^n\tribar^2 + \max_{0\leq n\leq \NT}\tribar e_p^n\tribar^2\right)^{1/2},\\
\tribar(e_y^n, e_p^n)\tribar_{\ell^{2}} & := \left(k\sum_{n=1}^{\NT}\tribar e_y^n \tribar^2\right)^{1/2} + \left(k\sum_{n=1}^{\NT}\tribar e_p^n\tribar^2\right)^{1/2},
\end{aligned}
\end{equation}
where $e_y^n : = Y^n - \overline y(t^n)$ and $e_p^n : = P^n - \overline p(t^n),\,n=0,\ldots,\NT.$ 

In all the numerical examples, we consider $T=0.5$ as the final time.  To compute the errors, initially, we assume a triangulation $\Th^{(0)}$ with $h_0^{(0)} = 1/M_0,\,M_0=4$ and time step $k^{(0)} = \frac{1}{10}h_0^{(0)}.$ Then, we perform uniform refinement, i.e., on each triangle we connect the midpoints to obtain four similar triangles. Thus, we obtain $\Th^{(\ell)},\,\ell=1,\ldots,6.$ Also, the each time step is computed as $k^{(\ell)} = \frac{1}{10}h_0^{(\ell)},\,\ell=1,\ldots,6.$

\subsection{Test example with a smooth solution}

We consider the following optimal control problem where the source term $G(\x,t)$ is chosen so as the control problem \eqref{control_problem}--\eqref{conv_diff} have the following solution for $(\x,t)\in [0,1]^2 \times [0,0.5],$ with $\lambda = 1,\,\bfb = (2,3)^T$ and $\sigma = 0,$
\begin{equation}\label{non_linear_system_f_sol}
\begin{aligned}
\overline y(\x,t) & = 100 t  (x^2 - x^4) (y - 3y^2 + 2y^3),\\
y_d(\x,t) & = y(\x,t) + p_t(\x,t) + \mu\,\Delta\,p(\x,t) + \bfb\cdot \nabla p(\x,t) - \sigma\,p(\x,t),\\
\overline p(\x,t) & = - 100 (T-t)  (x^2 - x^4) (y - 3y^2 + 2y^3),\\
\overline u(\x,t) & = \max\{\mathsf{u_a}, \min\{\mathsf{u_b}, - \lambda^{-1}p(\x,t)\}\},\\
(\mathsf{u_a},\mathsf{u_b}) & = (-3,3),\;\;\lambda=1,\;\;\bfb = (2,3)^T.
\end{aligned}
\end{equation}

\begin{figure}
\centering

\begin{subfigure}[b]{0.32\textwidth}
\centering
\includegraphics[width=\textwidth]{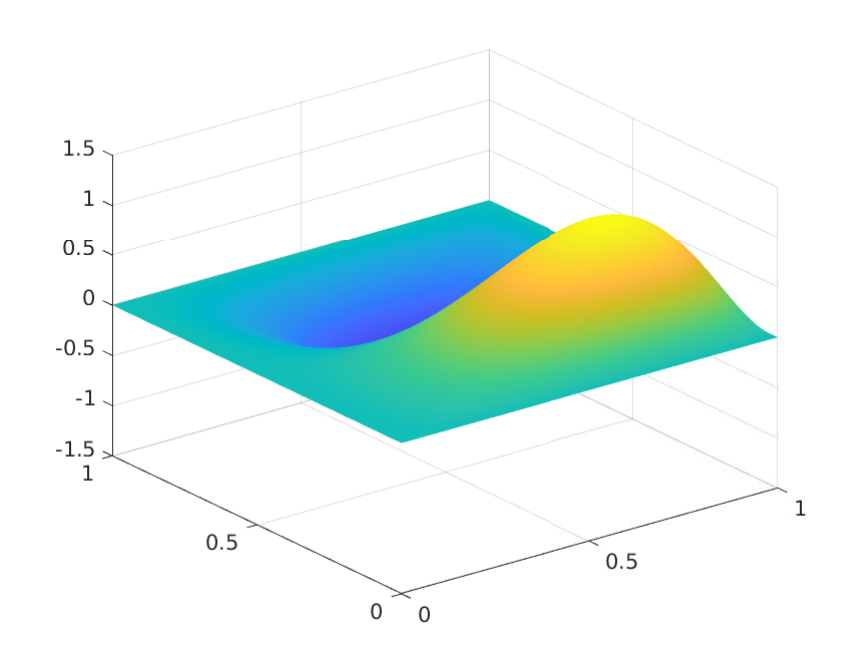}
\caption{State equation with $\mu = 10^{-8}$ using AFC.}
\end{subfigure}
\hspace{.1cm}
\begin{subfigure}[b]{0.32\textwidth}
\centering
\includegraphics[width=\textwidth]{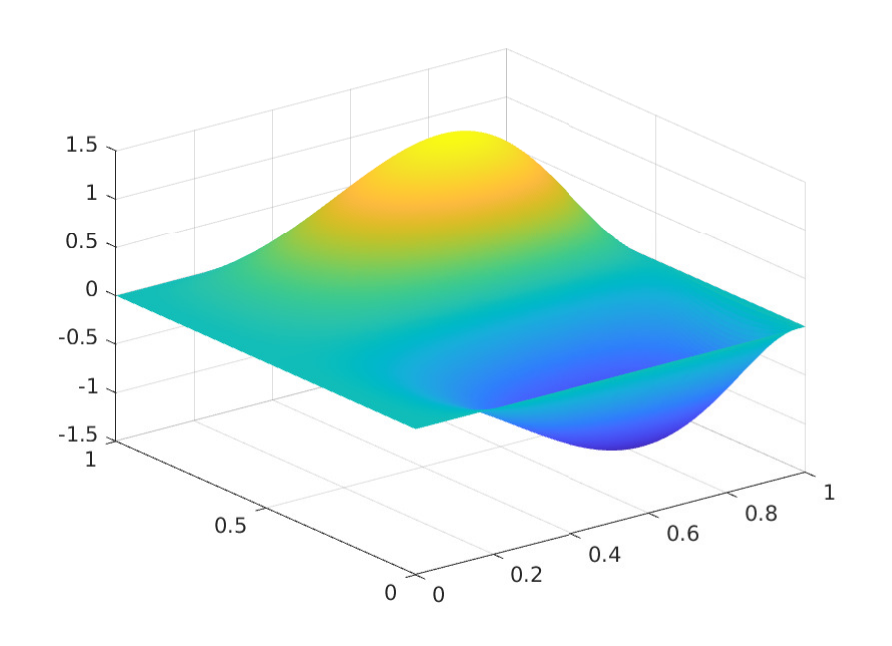}
\caption{Co-state equation with $\mu = 10^{-8}$ using AFC.}
\end{subfigure}

\begin{subfigure}[b]{0.32\textwidth}
\centering
\begin{tikzpicture}[scale=0.6]
\begin{axis}[
    xlabel={$h_0$},
    ylabel={Error},
    xmode=log,
    ymode=log,
    legend entries={
        {$\|(e_y^n, e_p^n)\|_{\max}$}, 
        {$\tribar (e_y^n, e_p^n) \tribar_{\max}$}, 
        {$\tribar (e_y^n, e_p^n) \tribar_{\ell^{2}}$}
    },
    legend style={font=\footnotesize, at={(1,1)}, anchor=north west},
    title={$\mu = 10^{-4}$},
    tick label style={font=\scriptsize},
    xlabel style={font=\small},
    ylabel style={font=\small},
]
\addplot[blue, mark=o, thick] table[x=h0, y=L2_mu1e4] {data_smooth.dat};
\addplot[blue, mark=triangle*, thick] table[x=h0, y=H1_mu1e4_max] {data_smooth.dat};
\addplot[blue, mark=square*, thick] table[x=h0, y=H1_mu1e4_mean] {data_smooth.dat};

\invLogSlopeTriangle{0.15}{0.2}{0.8}{-1}{black};
\logLogSlopeTriangle{0.35}{0.2}{0.1}{2}{black};
\end{axis}
\end{tikzpicture}
\caption{Experimental order of convergence for $\mu = 10^{-4}$.}
\end{subfigure}
\hspace{1.5cm}
\begin{subfigure}[b]{0.32\textwidth}
\centering
\begin{tikzpicture}[scale=0.6]
\begin{axis}[
    xlabel={$h_0$},
    ylabel={Error},
    xmode=log,
    ymode=log,
    legend entries={
        {$\|(e_y^n, e_p^n)\|_{\max}$}, 
        {$\tribar (e_y^n, e_p^n) \tribar_{\max}$}, 
        {$\tribar (e_y^n, e_p^n) \tribar_{\ell^{2}}$}
    },
    legend style={font=\footnotesize, at={(1,1)}, anchor=north west},
    title={$\mu = 10^{-5}$},
    tick label style={font=\scriptsize},
    xlabel style={font=\small},
    ylabel style={font=\small},
]
\addplot[red, mark=o, thick] table[x=h0, y=L2_mu1e5] {data_smooth.dat};
\addplot[red, mark=triangle*, thick] table[x=h0, y=H1_mu1e5_max] {data_smooth.dat};
\addplot[red, mark=square*, thick] table[x=h0, y=H1_mu1e5_mean] {data_smooth.dat};

\invLogSlopeTriangle{0.15}{0.2}{0.8}{-1}{black};
\logLogSlopeTriangle{0.4}{0.2}{0.1}{2}{black};
\end{axis}
\end{tikzpicture}
\caption{Experimental order of convergence for $\mu = 10^{-5}$.}
\end{subfigure}
\begin{subfigure}[b]{0.32\textwidth}
\centering
\begin{tikzpicture}[scale=0.6]
\begin{axis}[
    xlabel={$h_0$},
    ylabel={Error},
    xmode=log,
    ymode=log,
    legend entries={
        {$\|(e_y^n, e_p^n)\|_{\max}$}, 
        {$\tribar (e_y^n, e_p^n) \tribar_{\max}$}, 
        {$\tribar (e_y^n, e_p^n) \tribar_{\ell^{2}}$}
    },
    legend style={font=\footnotesize, at={(1,1)}, anchor=north west},
    title={$\mu = 10^{-8}$},
    tick label style={font=\scriptsize},
    xlabel style={font=\small},
    ylabel style={font=\small},
]
\addplot[magenta, mark=o, thick] table[x=h0, y=L2_mu1e8] {data_smooth.dat};
\addplot[magenta, mark=triangle*, thick] table[x=h0, y=H1_mu1e8_max] {data_smooth.dat};
\addplot[magenta, mark=square*, thick] table[x=h0, y=H1_mu1e8_mean] {data_smooth.dat};

\invLogSlopeTriangle{0.15}{0.2}{0.8}{-1}{black};
\logLogSlopeTriangle{0.4}{0.2}{0.1}{2}{black};
\end{axis}
\end{tikzpicture}
\caption{Experimental order of convergence for $\mu = 10^{-8}$.}
\end{subfigure}

\caption{Results using AFC: (top) state and co-state for $\mu = 10^{-8}$, (bottom) convergence studies for $\mu = 10^{-4},\, 10^{-5},\,10^{-8}$.}\label{Fig:smooth_solution}
\end{figure}

Both the state and co-state variable are smooth functions and their profiles are depicted in Fig. \ref{Fig:smooth_solution} for $\mu = 10^{-8},$ at $\Th^{(6)}$, using the mesh step $h^{(6)}$ and time step $k^{(6)}$. We set $\mu=10^{-4},\,10^{-5},\,10^{-8},$ and we compute the errors in Fig. \ref{Fig:smooth_solution}  as described in \eqref{computed_errors}. The latter figure confirms the results derived in Theorem \ref{theorem:estimate_global}.

\subsection{Test example with a hump changing its height}

We consider the optimal control problem \eqref{control_problem}--\eqref{conv_diff} for $(\x,t)\in [0,1]^2 \times [0,0.5]$ with $\lambda = \sigma = 1,\,\bfb = (2,3)^T.$ The source term $G(\x,t)$ is chosen, so that, the exact solution is 
\begin{equation}\label{non_linear_system_f_sol_hump}
\begin{aligned}
\overline y(\x,t) & = 16 \sin(\pi t) x (1 - x)  y  (1 - y)\left( \frac{1}{2} + \frac{1}{\pi}  \arctan\left( \frac{2}{\sqrt{\mu}} \left(\frac{1}{16} - \left(x-\frac{1}{2}\right)^2 - \left(y-\frac{1}{2}\right)^2 \right)   \right)\right),\\
y_d(\x,t) & = y(\x,t) + p_t(\x,t) + \mu\,\Delta\,p(\x,t) + \bfb\cdot \nabla p(\x,t) - p(\x,t),\\
\overline p(\x,t) & = -(T-t)x (1 - x)  y  (1 - y),\\
\overline u(\x,t) & = \max\{\mathsf{u_a}, \min\{\mathsf{u_b}, - \lambda^{-1}p(\x,t)\}\},\\
(\mathsf{u_a},\mathsf{u_b}) & = (-2,2),\;\;\lambda=1,\;\;\bfb = (2,3)^T.
\end{aligned}
\end{equation}
The state variable $\overline y$ develops steep gradients along the circle $0.0625 - (x-0.5)^2 - (y - 0.5)^2 =0$ of width $\sqrt{\mu}.$ In particular, the state variable is a hump that changes its height in the course of time. The maximal height is attained at $t = 0.5.$ This function has also studied in the context of the time--depended convection--diffusion equations, see, e.g., \cite[Section 7.2]{john2008}. The co-state variable $\overline p,$ is taken to be a smooth function.\par
\begin{figure}
\centering

\begin{subfigure}[b]{0.32\textwidth}
\centering
\includegraphics[width=\textwidth]{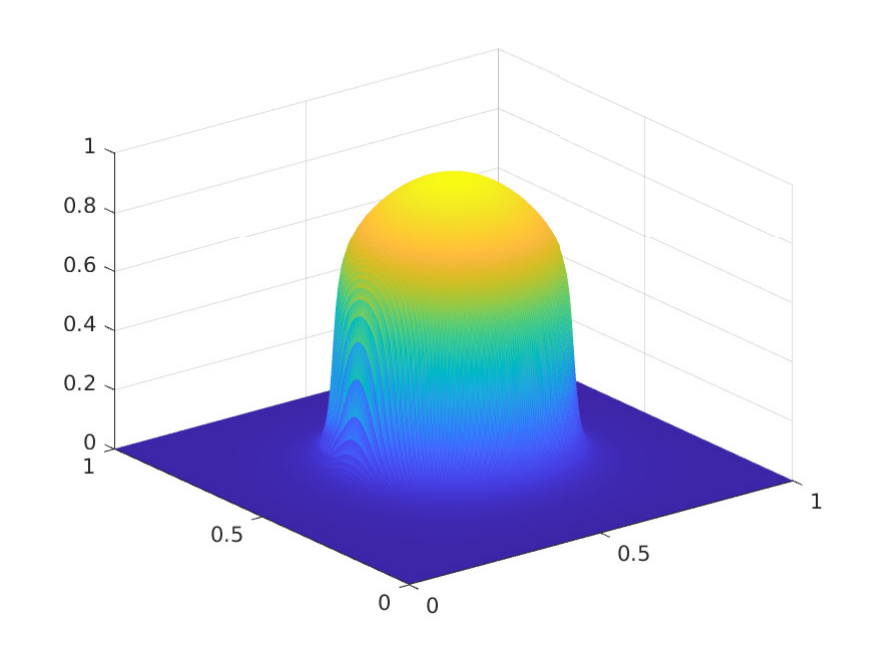}
\caption{State equation with $\mu = 10^{-4}$ using AFC.}
\end{subfigure}
\hspace{.1cm}
\begin{subfigure}[b]{0.32\textwidth}
\centering
\includegraphics[width=\textwidth]{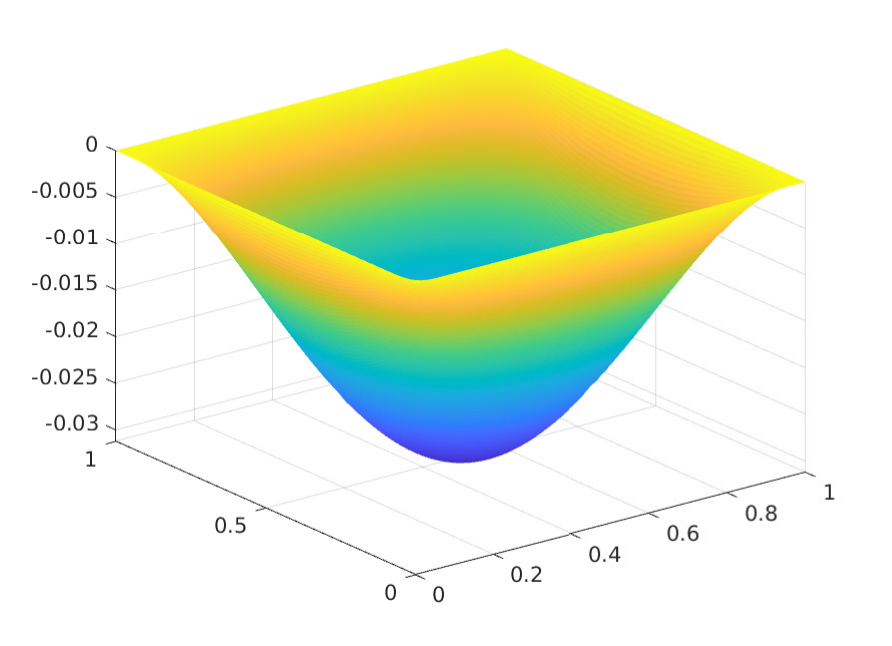}
\caption{Co-state equation with $\mu = 10^{-4}$ using AFC.}
\end{subfigure}

\begin{subfigure}[b]{0.32\textwidth}
\centering
\begin{tikzpicture}[scale=0.6]
\begin{axis}[
    xlabel={$h_0$},
    ylabel={Error},
    xmode=log,
    ymode=log,
    legend entries={
        {$\|(e_y^n, e_p^n)\|_{\max}$}, 
        {$\tribar (e_y^n, e_p^n) \tribar_{\max}$}, 
        {$\tribar (e_y^n, e_p^n) \tribar_{\ell^{2}}$}
    },
    legend style={font=\footnotesize, at={(1,1)}, anchor=north west},
    title={$\mu = 10^{-4}$},
    tick label style={font=\scriptsize},
    xlabel style={font=\small},
    ylabel style={font=\small},
]
\addplot[blue, mark=o, thick] table[x=h0, y=L2_mu1e4] {data_hump.dat};
\addplot[blue, mark=triangle*, thick] table[x=h0, y=H1_mu1e4_max] {data_hump.dat};
\addplot[blue, mark=square*, thick] table[x=h0, y=H1_mu1e4_mean] {data_hump.dat};

\invLogSlopeTriangle{0.15}{0.2}{0.7}{-1}{black};
\logLogSlopeTriangle{0.35}{0.2}{0.1}{2}{black};
\end{axis}
\end{tikzpicture}
\caption{Experimental order of convergence for $\mu = 10^{-4}$.}
\end{subfigure}
\hspace{1.5cm}
\begin{subfigure}[b]{0.32\textwidth}
\centering
\begin{tikzpicture}[scale=0.6]
\begin{axis}[
    xlabel={$h_0$},
    ylabel={Error},
    xmode=log,
    ymode=log,
    legend entries={
        {$\|(e_y^n, e_p^n)\|_{\max}$}, 
        {$\tribar (e_y^n, e_p^n) \tribar_{\max}$}, 
        {$\tribar (e_y^n, e_p^n) \tribar_{\ell^{2}}$}
    },
    legend style={font=\footnotesize, at={(1,1)}, anchor=north west},
    title={$\mu = 10^{-5}$},
    tick label style={font=\scriptsize},
    xlabel style={font=\small},
    ylabel style={font=\small},
]
\addplot[red, mark=o, thick] table[x=h0, y=L2_mu1e5] {data_hump.dat};
\addplot[red, mark=triangle*, thick] table[x=h0, y=H1_mu1e5_max] {data_hump.dat};
\addplot[red, mark=square*, thick] table[x=h0, y=H1_mu1e5_mean] {data_hump.dat};

\invLogSlopeTriangle{0.15}{0.2}{0.7}{-1}{black};
\logLogSlopeTriangle{0.4}{0.12}{0.1}{2}{black};
\end{axis}
\end{tikzpicture}
\caption{Experimental order of convergence for $\mu = 10^{-5}$.}
\end{subfigure}

\caption{Results using AFC: (top) state and co-state for $\mu = 10^{-4}$, (bottom) convergence studies for $\mu = 10^{-4}$ and $\mu = 10^{-5}$.}\label{Fig:hump_solution}
\end{figure}

In Fig. \ref{Fig:hump_solution}, we plot the profiles of the state and co-state solutions for $\mu = 10^{-4}$ at $\Th^{(6)}$, using the mesh step $h^{(6)}$ and time step $k^{(6)}$. In addition, we compute the errors \eqref{computed_errors} for $\mu = 10^{-4}$ and $\mu = 10^{-5}$. It is readily verified that, asymptotically, the estimates established in Theorem \ref{theorem:estimate_global} are confirmed.

\subsection{Test example with interior layer}

We consider the optimal control problem \eqref{control_problem}--\eqref{conv_diff} for $(\x,t)\in [0,1]^2 \times [0,0.5]$ with $\lambda = \sigma = 1,\,\bfb = (2,3)^T.$ The source term $G(\x,t)$ is chosen, so that, the exact solution is 
\begin{equation}\label{non_linear_system_f_sol_interior_layer}
\begin{aligned}
\overline y(\x,t) & = (1 - e^{-t})\arctan\left(\frac{1}{\sqrt{\mu}}\left(-\frac{1}{2}x + y - \frac{1}{4}\right) \right),\\
y_d(\x,t) & = y(\x,t) + p_t(\x,t) + \mu\,\Delta\,p(\x,t) + \bfb\cdot \nabla p(\x,t) - p(\x,t),\\
\overline p(\x,t) & = 32 (T - t) x (1 - x)  y  (1 - y)\left( \frac{1}{2} + \frac{1}{\pi}  \arctan\left( \frac{2}{\sqrt{\mu}} \left(\frac{1}{16} - \left(x-\frac{1}{2}\right)^2 - \left(y-\frac{1}{2}\right)^2 \right)   \right)\right),\\
\overline u(\x,t) & = \max\{\mathsf{u_a}, \min\{\mathsf{u_b}, - \lambda^{-1}p(\x,t)\}\},\\
(\mathsf{u_a},\mathsf{u_b}) & = (-2,2),\;\;\lambda=1,\;\;\bfb = (2,3)^T.
\end{aligned}
\end{equation}
The above example can be found in the context of the elliptic convection--diffusion in \cite[Example 5.2]{hinze2009}. The state variable $\overline y$ develops a interior layer along the line $-0.5x + y - 0.25 = 0,$ while the co-state variable $\overline p,$ develops steep gradients along the circle $0.0625 - (x-0.5)^2 - (y - 0.5)^2 =0.$\par
\begin{figure}
\centering

\begin{subfigure}[b]{0.32\textwidth}
\centering
\includegraphics[width=\textwidth]{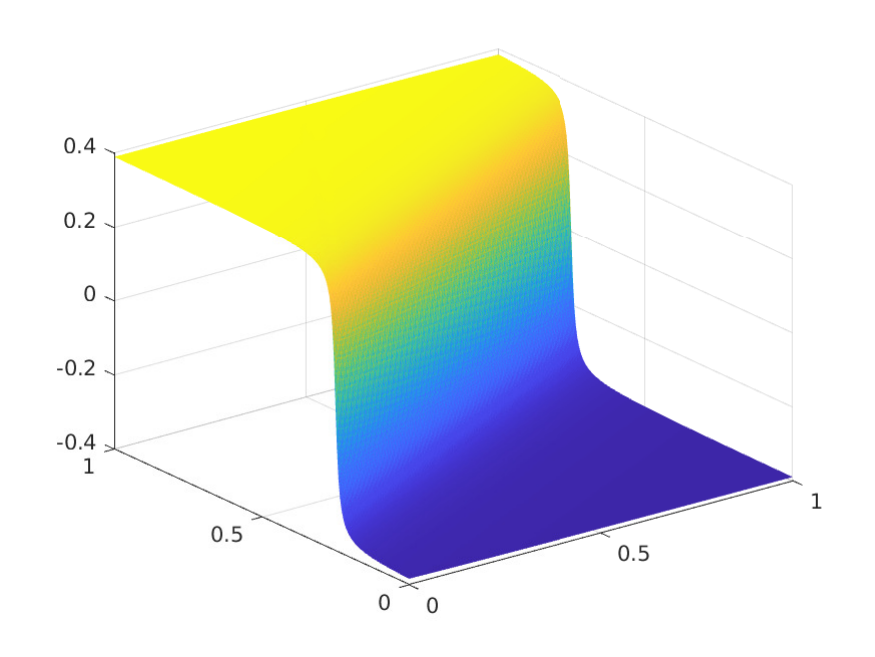}
\caption{State equation with $\mu = 10^{-4}$ using AFC.}
\end{subfigure}
\hspace{.1cm}
\begin{subfigure}[b]{0.32\textwidth}
\centering
\includegraphics[width=\textwidth]{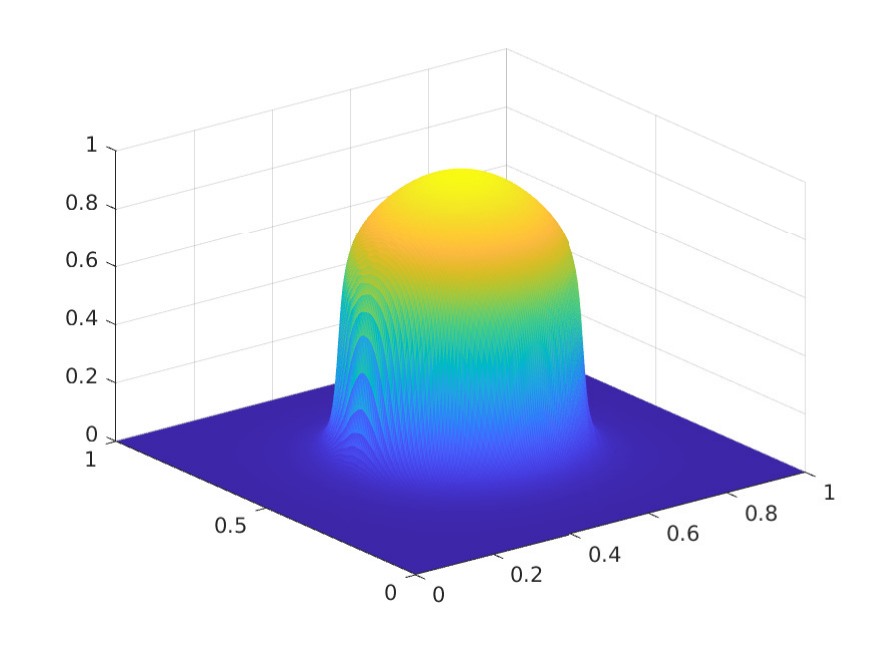}
\caption{Co-state equation with $\mu = 10^{-4}$ using AFC.}
\end{subfigure}

\begin{subfigure}[b]{0.32\textwidth}
\centering
\begin{tikzpicture}[scale=0.6]
\begin{axis}[
    xlabel={$h_0$},
    ylabel={Error},
    xmode=log,
    ymode=log,
    legend entries={
        {$\|(e_y^n, e_p^n)\|_{\max}$}, 
        {$\tribar (e_y^n, e_p^n) \tribar_{\max}$}, 
        {$\tribar (e_y^n, e_p^n) \tribar_{\ell^{2}}$}
    },
    legend style={font=\footnotesize, at={(1,1)}, anchor=north west},
    title={$\mu = 10^{-4}$},
    tick label style={font=\scriptsize},
    xlabel style={font=\small},
    ylabel style={font=\small},
]
\addplot[blue, mark=o, thick] table[x=h0, y=L2_mu1e4] {data_layer.dat};
\addplot[blue, mark=triangle*, thick] table[x=h0, y=H1_mu1e4_max] {data_layer.dat};
\addplot[blue, mark=square*, thick] table[x=h0, y=H1_mu1e4_mean] {data_layer.dat};

\invLogSlopeTriangle{0.15}{0.2}{0.7}{-1}{black};
\logLogSlopeTriangle{0.35}{0.2}{0.1}{2}{black};
\end{axis}
\end{tikzpicture}
\caption{Experimental order of convergence for $\mu = 10^{-4}$.}
\end{subfigure}
\hspace{1.5cm}
\begin{subfigure}[b]{0.32\textwidth}
\centering
\begin{tikzpicture}[scale=0.6]
\begin{axis}[
    xlabel={$h_0$},
    ylabel={Error},
    xmode=log,
    ymode=log,
    legend entries={
        {$\|(e_y^n, e_p^n)\|_{\max}$}, 
        {$\tribar (e_y^n, e_p^n) \tribar_{\max}$}, 
        {$\tribar (e_y^n, e_p^n) \tribar_{\ell^{2}}$}
    },
    legend style={font=\footnotesize, at={(1,1)}, anchor=north west},
    title={$\mu = 10^{-5}$},
    tick label style={font=\scriptsize},
    xlabel style={font=\small},
    ylabel style={font=\small},
]
\addplot[red, mark=o, thick] table[x=h0, y=L2_mu1e5] {data_layer.dat};
\addplot[red, mark=triangle*, thick] table[x=h0, y=H1_mu1e5_max] {data_layer.dat};
\addplot[red, mark=square*, thick] table[x=h0, y=H1_mu1e5_mean] {data_layer.dat};

\invLogSlopeTriangle{0.15}{0.2}{0.7}{-1}{black};
\logLogSlopeTriangle{0.4}{0.12}{0.1}{2}{black};
\end{axis}
\end{tikzpicture}
\caption{Experimental order of convergence for $\mu = 10^{-5}$.}
\end{subfigure}

\caption{Results using AFC: (top) state and co-state for $\mu = 10^{-4}$, (bottom) convergence studies for $\mu = 10^{-4}$ and $\mu = 10^{-5}$.}
\label{Fig:layer_solution}
\end{figure}

Similar to previous example, we plot the state and co-state profile of the solutions, see Fig. \ref{Fig:layer_solution}. Furthermore, for $\mu = 10^{-4}$ and $10^{-5}$. As can be verified, the numerical results asymptotically validate the estimates established in Theorem \ref{theorem:estimate_global}.

\section{Conclusions}

In this paper, we considered an optimal control problem that governed by a parabolic convection--diffusion--reaction equation on a bounded domain $\Omega\subset \mathbb{R}^2.$ We focused on the optimize--then--discretize approach, considering the piecewise finite element method alongside the algebraic flux correction method for stabilization. For temporal discretization, we used the backward Euler method. The discrete control variable is obtained by projecting the discretized adjoint state onto the set of admissible controls. We proved results concerning the existence and uniqueness of the fully--discrete scheme as well as we derived optimal error estimates under a sufficient choice of time and mesh step in $L^{2}$ and in energy norm with respect to the spatial variable and in $\ell^\infty$ in time. We also presented numerical experiments that validate the order of convergence of the stabilized fully--discrete scheme using the algebraic flux correction method.
 
\bigskip
\bibliographystyle{plain} 
\bibliography{ref}

\begin{thebibliography}{10}

\bibitem{arada2002}
N.~Arada, E.~Casas, and F.~Tr{\"o}ltzsch.
\newblock {Error estimates for a semilinear elliptic control problem}.
\newblock {\em Comput. Optim. Appl.}, 23:201--229, 2002.
\newblock \url{https://doi.org/10.1023/A:1020576801966}.

\bibitem{barrenechea2016}
G.~R. Barrenechea, V.~John, and P.~Knobloch.
\newblock {Analysis of algebraic flux correction schemes}.
\newblock {\em SIAM J. Numer. Anal.}, 54:2427--2451, 2016.
\newblock \url{https://doi.org/10.1137/15M1018216}.

\bibitem{barrenechea2017b}
G.~R. Barrenechea, V.~John, and P.~Knobloch.
\newblock {An algebraic flux correction scheme satisfying the discrete maximum
  principle and linearity preservation on general meshes}.
\newblock {\em Math. Models Methods Appl. Sci.}, 27:525--548, 2017.
\newblock \url{https://doi.org/10.1142/S0218202517500087}.

\bibitem{barrenechea2024}
G.~R. Barrenechea, V.~John, and P.~Knobloch.
\newblock {Finite element methods respecting the discrete maximum principle for
  convection--diffusion equations}.
\newblock {\em SIAM Rev.}, 66(1):3--88, 2024.
\newblock \url{https://doi.org/10.1137/22M1488934}.

\bibitem{barrenechea2025}
G.~R. Barrenechea, V.~John, and P.~Knobloch.
\newblock {\em {Monotone Discretizations for Elliptic Second Order Partial
  Differential Equations}}.
\newblock Springer Series in Computational Mathematics. Springer Cham, 1
  edition, 2025.
\newblock \url{https://doi.org/10.1007/978-3-031-80684-1}.

\bibitem{barrenechea2018}
G.~R. Barrenechea, V.~John, P.~Knobloch, and R.~Rankin.
\newblock {A unified analysis of algebraic flux correction schemes for
  convection--diffusion equations}.
\newblock {\em SeMA}, 75:655--685, 2018.
\newblock \url{https://doi.org/10.1007/s40324-018-0160-6}.

\bibitem{bartels2022}
S.~Bartels, C.~Palus, and Z.~Wang.
\newblock {Quasi-optimal error estimates for the approximation of stable
  harmonic maps}.
\newblock {\em arXiv preprint \url{https://arxiv.org/abs/2209.11985v1}}, 2022.

\bibitem{baumgartner2022}
J.~Baumgartner.
\newblock {\em {Optimal control problems and algebraic flux correction
  schemes}}.
\newblock PhD thesis, Universit{\"a}t Duisburg--Essen, 2022.
\newblock \url{https://doi.org/10.17185/duepublico/75439}.

\bibitem{baumgartner2025}
J.~Baumgartner and A.~R{\"o}sch.
\newblock {AFC Stabilization of a State Constrained Optimal Control Problem}.
\newblock {\em Math. Control Relat. Fields}, 2025.
\newblock \url{https://doi.org/10.3934/mcrf.2025033}.

\bibitem{becker2007}
R.~Becker, D.~Meidner, and B.~Vexler.
\newblock {Efficient numerical solution of parabolic optimization problems by
  finite element methods}.
\newblock {\em Optim. Methods Softw.}, 22(5):813--833, 2007.
\newblock \url{https://doi.org/10.1080/10556780701228532}.

\bibitem{becker2007b}
R.~Becker and B.~Vexler.
\newblock {Optimal control of the convection--diffusion equation using
  stabilized finite element methods}.
\newblock {\em Numer. Math.}, 106:349--367, 2007.
\newblock \url{https://doi.org/10.1007/s00211-007-0067-0}.

\bibitem{becker2007a}
R.~Becker and B.~Vexler.
\newblock {Optimal Control of the Convection-Diffusion Equation Using
  Stabilized Finite Element Methods}.
\newblock {\em Numer. Math.}, 106(3):349--367, 2007.
\newblock \url{http://dx.doi.org/10.1007/s00211-007-0067-0}.

\bibitem{brenner2008}
S.~C. Brenner and L.~R. Scott.
\newblock {\em {The mathematical theory of finite element methods}}.
\newblock Springer, New York, 2 edition, 2008.
\newblock \url{https://doi.org/10.1007/978-0-387-75934-0}.

\bibitem{brooks1982}
A.~N. Brooks and T.~J.~R. Hughes.
\newblock {Streamline upwind/Petrov--Galerkin formulations for convection
  dominated flows with particular emphasis on the incompressible Navier--Stokes
  equations}.
\newblock {\em Comput. Methods Appl. Mech. Eng.}, 32(1--3):199--259, 1982.
\newblock \url{https://doi.org/10.1016/0045-7825(82)90071-8}.

\bibitem{burman2004}
E.~Burman and P.~Hansbo.
\newblock {Edge Stabilization for Galerkin Approximations of
  Convection-Diffusion-Reaction Problems}.
\newblock {\em Comput. Methods Appl. Mech. Eng.}, 193:1437--1453, 2004.
\newblock \url{https://doi.org/10.1016/j.cma.2003.12.032}.

\bibitem{casas2005}
E.~Casas, M.~Mateos, and F.~Tr{\"o}ltzsch.
\newblock {Error estimates for the numerical approximation of boundary
  semilinear elliptic control problems}.
\newblock {\em Comput. Optim. Appl.}, 31:193--220, 2005.
\newblock \url{https://doi.org/10.1007/s10589-005-2180-2}.

\bibitem{chatzipantelidis2012}
P.~Chatzipantelidis, R.~Lazarov, and V.~Thom{\'e}e.
\newblock {Some error estimates for the lumped mass finite element method for a
  parabolic problem}.
\newblock {\em Math. Comp.}, 81:1--20, 2012.
\newblock \url{https://doi.org/10.1090/S0025-5718-2011-02503-2}.

\bibitem{chatzipantelidis2022}
P.~Chatzipantelidis and C.~Pervolianakis.
\newblock {Error analysis of a backward Euler positive preserving stabilized
  scheme for a Chemotaxis system}.
\newblock {\em arXiv preprint (v5 version)}, 2022.
\newblock \url{https://arxiv.org/pdf/2210.04709v5}.

\bibitem{collis2002}
S.~S. Collis and M.~Heinkenschloss.
\newblock {Analysis of the Streamline Upwind/Petrov Galerkin method applied to
  the solution of optimal control problems}.
\newblock Technical Report TR02-01, CAAM, 2002.
\newblock Also available at \url{https://arxiv.org/abs/2411.09828}.

\bibitem{draganescu2004}
A.~Dr\u{a}g\u{a}nescu, T.~F. Dupont, and L.~R. Scott.
\newblock {Failure of the discrete maximum principle for an elliptic finite
  element problem}.
\newblock {\em Math. Comp.}, 74:1--23, 2004.
\newblock \url{https://doi.org/10.1090/s0025-5718-04-01651-5}.

\bibitem{evans2010}
L.~Evans.
\newblock {\em {Partial differential equations}}.
\newblock American Mathematical Society, second edition, 2010.

\bibitem{fu2009}
H.~Fu and H.~Rui.
\newblock {A priori error estimates for optimal control problems governed by
  transient advection--diffusion equations}.
\newblock {\em J. Sci. Comput.}, 38:290--315, 2009.
\newblock \url{https://doi.org/10.1007/s10915-008-9224-6}.

\bibitem{guermond1999}
J.~L. Guermond.
\newblock Stabilization of galerkin approximations of transport equations by
  subgrid modeling.
\newblock {\em Modél. Math. Anal. Numér.}, 36(6):1293--1316, 1999.
\newblock \url{https://doi.org/10.1051/m2an:1999145}.

\bibitem{hinze2005}
M.~Hinze.
\newblock {A variational discretization concept in control constrained
  optimization: the linear--quadratic case}.
\newblock {\em Comput. Optim. Appl.}, 30:45--61, 2005.
\newblock \url{https://doi.org/10.1007/s10589-005-4559-5}.

\bibitem{hinze2009book}
M.~Hinze, R.~Pinnau, M.~Ulbrich, and S.~Ulbrich.
\newblock {\em {Optimization with PDE constraints}}.
\newblock Mathematical Modelling: Theory and Applications. Springer Dordrecht,
  2009.
\newblock \url{https://doi.org/10.1007/978-1-4020-8839-1}.

\bibitem{hinze2009}
M.~Hinze, N.~Yan, and Z.~Zhou.
\newblock {Variational discretization for optimal control governed by
  convection dominated diffusion equations}.
\newblock {\em J. Comput. Math.}, 27(2--3):237--253, 2009.
\newblock \url{http://www.jstor.org/stable/43693504}.

\bibitem{jha2021}
A.~Jha and N.~Ahmed.
\newblock {Analysis of flux corrected transport schemes for evolutionary
  convection--diffusion--reaction equations}.
\newblock {\em arXiv preprint}, 2021.
\newblock \url{https://doi.org/10.48550/arXiv.2103.04776}.

\bibitem{jin2020}
B.~Jin, B.~Li, and Z.~Zhou.
\newblock Pointwise-in--time error estimates for an optimal control problem
  with subdiffusion constraint.
\newblock {\em IMA J. Numer. Anal.}, 40:377--404, 2020.
\newblock \url{https://doi.org/10.1093/imanum/dry064}.

\bibitem{john2007}
V.~John and P.~Knobloch.
\newblock {On spurious oscillations at layers diminishing SOLD methods for
  convection--diffusion equations. I. A review}.
\newblock {\em Comput. Methods Appl. Mech. Eng.}, 196(17--20):2197--2215, 2007.
\newblock \url{https://doi.org/10.1016/j.cma.2006.11.013}.

\bibitem{john2021}
V.~John and P.~Knobloch.
\newblock {Existence of solutions of a finite element
  flux--corrected--transport scheme}.
\newblock {\em Appl. Math. Lett.}, 115:106932, 2021.
\newblock \url{https://doi.org/10.1016/j.aml.2020.106932}.

\bibitem{john2008}
V.~John and E.~Schmeyer.
\newblock {Finite element methods for time-dependent
  convection--diffusion--reaction equations with small diffusion}.
\newblock {\em Comput. Methods Appl. Mech. Engrg.}, 198(3):475--494, 2008.
\newblock \url{https://doi.org/10.1016/j.cma.2008.08.016}.

\bibitem{knobloch2017}
P.~Knobloch.
\newblock {On the discrete maximum principle for algebraic flux correction
  schemes with limiters of upwind type}.
\newblock In Z.~Huang, M.~Stynes, and Z.~Zhang, editors, {\em Boundary and
  Interior Layers, Computational and Asymptotic Methods BAIL 2016}, volume 120
  of {\em Lect. Notes Comput. Sci. Eng.}, pages 129--139. Springer, Berlin,
  2017.
\newblock \url{https://doi.org/10.1007/978-3-319-67202-1_10}.

\bibitem{kuzmin2010book}
D.~Kuzmin.
\newblock {\em {A guide to numerical methods for transport equations}}.
\newblock University Erlangen--Nuremberg, Nuremberg, 2010.

\bibitem{kuzmin2005}
D.~Kuzmin and M.~M{\"o}ller.
\newblock {Algebraic flux correction I. Scalar conservation laws}.
\newblock In {\em Flux--corrected transport: principles, algorithms, and
  applications}, pages 155--206. Springer Berlin Heidelberg, Berlin,
  Heidelberg, 2005.
\newblock \url{https://doi.org/10.1007/3-540-27206-2_6}.

\bibitem{kuzmin2002}
D.~Kuzmin and S.~Turek.
\newblock {Flux correction tools for finite elements}.
\newblock {\em J. Comput. Phys.}, 175:525--558, 2002.
\newblock \url{https://doi.org/10.1006/jcph.2001.6955}.

\bibitem{kuzmin2004}
D.~Kuzmin and S.~Turek.
\newblock {High--resolution FEM--TVD schemes based on a fully multidimensional
  flux limiter}.
\newblock {\em J. Comput. Phys.}, 198:131--158, 2004.
\newblock \url{https://doi.org/10.1016/j.jcp.2004.01.015}.

\bibitem{lions1971}
J.~L. Lions.
\newblock {\em {Optimal Control of Systems Governed by Partial Differential
  Equations}}, volume 170 of {\em Grundlehren der mathematischen
  Wissenschaften}.
\newblock Springer, Berlin, Heidelberg, 1971.
\newblock Softcover reprint of the original 1st edition.

\bibitem{lohmann2019}
C.~Lohmann.
\newblock {\em {Physics-Compatible Finite Element Methods for Scalar and
  Tensorial Advection Problems}}.
\newblock Springer Spektrum Wiesbaden, 1 edition, 2019.
\newblock \url{https://doi.org/10.1007/978-3-658-27737-6}.

\bibitem{meidner2008}
D.~Meidner and B.~Vexler.
\newblock {A priori error estimates for space--time finite element
  discretization of parabolic optimal control problems. Part I: problems
  without control constraints}.
\newblock {\em SIAM J. Control Optim.}, 47(3):1150--1177, 2008.
\newblock \url{https://doi.org/10.1137/070694016}.

\bibitem{meidner2008b}
D.~Meidner and B.~Vexler.
\newblock {A priori error estimates for space--time finite element
  discretization of parabolic optimal control problems. Part II: problems with
  control constraints}.
\newblock {\em SIAM J. Control Optim.}, 47(3):1301--1329, 2008.
\newblock \url{https://doi.org/10.1137/070694028}.

\bibitem{meidner2011}
D.~Meidner and B.~Vexler.
\newblock {A priori error analysis of the Petrov–Galerkin Crank–Nicolson
  scheme for parabolic optimal control problems}.
\newblock {\em SIAM J. Control Optim.}, 49(5):2183--2211, 2011.
\newblock \url{https://doi.org/10.1137/100809611}.

\bibitem{meyer2004}
C.~Meyer and A.~R{\"o}sch.
\newblock {Superconvergence properties of optimal control problems}.
\newblock {\em SIAM J. Control Optim.}, 43:970--985, 2004.
\newblock \url{https://doi.org/10.1137/S0363012903431608}.

\bibitem{neitzel2012}
I.~Neitzel and B.~Vexler.
\newblock {A priori error estimates for space--time finite element
  discretization of semilinear parabolic optimal control problems}.
\newblock {\em Numer. Math.}, 120(2):345--386, 2012.
\newblock \url{https://doi.org/10.1007/s00211-011-0409-9}.

\bibitem{ryu2001}
S.~U. Ryu and A.~Yagi.
\newblock {Optimal control of Keller--Segel equations}.
\newblock {\em J. Math. Anal. Appl.}, 256(1):45--66, 2001.
\newblock \url{https://doi.org/10.1006/jmaa.2000.7254}.

\bibitem{springer2014}
A.~Springer and B.~Vexler.
\newblock {Third order convergent time discretization for parabolic optimal
  control problems with control constraints}.
\newblock {\em Comput. Optim. Appl.}, 57:205--240, 2014.
\newblock \url{https://doi.org/10.1007/s10589-013-9580-5}.

\bibitem{thomee2006}
V.~Thom{\'e}e.
\newblock {\em {Galerkin finite element methods for parabolic problems}}.
\newblock Springer, 2 edition, 2006.
\newblock \url{https://doi.org/10.1007/3-540-33122-0}.

\bibitem{vondaniels2015}
N.~von Daniels, M.~Hinze, and M.~Vierling.
\newblock {Crank--Nicolson time stepping and variational discretization of
  control-constrained parabolic optimal control problems}.
\newblock {\em SIAM J. Control Optim.}, 53(3):1500--1523, 2015.
\newblock \url{https://doi.org/10.1137/14099680X}.

\bibitem{weng2015}
Z.~Weng, J.~Z. Yang, and X.~Lu.
\newblock {A Stabilized Finite Element Method for the Convection Dominated
  Diffusion Optimal Control Problem}.
\newblock {\em Appl. Anal.}, pages 1--17, 2015.
\newblock \url{https://doi.org/10.1080/00036811.2015.1114606}.

\bibitem{yan2009}
N.~Yan and Z.~Zhou.
\newblock {A priori and a posteriori error analysis of edge stabilization
  {G}alerkin method for the optimal control problem governed by
  convection-dominated diffusion equation}.
\newblock {\em J. Comput. Appl. Math.}, 223(1):198--217, 2009.
\newblock \url{https://doi.org/10.1016/j.cam.2008.01.006}.

\bibitem{zhang2023}
X.~Zhang, J.~Zhao, and Y.~Hou.
\newblock {A Priori Error Estimates of Crank--Nicolson Finite Element Method
  for Parabolic Optimal Control Problems}.
\newblock {\em Comput. Math. Appl.}, 144:274--289, 2023.
\newblock \url{https://doi.org/10.1016/j.camwa.2023.06.017}.

\bibitem{zhou2010}
Z.~Zhou and N.~Yan.
\newblock {The Local Discontinuous Galerkin Method for Optimal Control Problem
  Governed by Convection Diffusion Equations}.
\newblock {\em Int. J. Numer. Anal. Model.}, 7(4):681--699, 2010.

\end{thebibliography}

\end{document}